\documentclass[12pt]{amsart}
\usepackage{amssymb}
\usepackage{mathrsfs} 
\usepackage{bm}
\usepackage{graphicx}
\usepackage{color}
\usepackage[all]{xy}

\allowdisplaybreaks[4]

\newtheorem{thm}{Theorem}[section]
\newtheorem{lem}[thm]{Lemma}
\newtheorem{cor}[thm]{Corollary}
\newtheorem{prop}[thm]{Proposition}
\newtheorem{definition}[thm]{Definition}

\def \para{\refstepcounter{thm} \par\medskip\noindent
                \textbf{\thethm .} }

\def \remark{\refstepcounter{thm} \par\medskip\noindent
                \textbf{Remark \thethm .} }

\numberwithin{equation}{thm}

\renewcommand\AA{\mathbb A} 
\newcommand\BB{\mathbb B}
\newcommand\CC{\mathbb C}

\newcommand\KK{\mathbb K}

\newcommand\XX{\mathbb X}

\newcommand\ZZ{\mathbb Z}

\newcommand\bQ{\mathbf Q}

\newcommand\bb{\mathbf b}

\newcommand\bk{\mathbf k}

\newcommand\bn{\mathbf n}

\newcommand\bu{\mathbf u}

\newcommand\bx{\mathbf x}

\newcommand\cA{\mathcal{A}}
\newcommand\cB{\mathcal{B}}

\newcommand\cD{\mathcal{D}}

\newcommand\cI{\mathcal{I}}
\newcommand\cJ{\mathcal{J}}
\newcommand\cK{\mathcal{K}}

\newcommand\cU{\mathcal{U}}

\newcommand\cX{\mathcal{X}}

\newcommand\fI{\mathfrak I}

\newcommand\fS{\mathfrak S}

\newcommand\fX{\mathfrak X}

\newcommand\sH{\mathscr H}

\newcommand\sS{\mathscr S}

\renewcommand\a{\alpha}  
\renewcommand\b{\beta}   
\newcommand\g{\gamma}  
\renewcommand\d{\delta}  

\newcommand\la{\lambda}

\newcommand\s{\sigma}

\newcommand\w{\omega}

\newcommand\ve{\varepsilon}

\newcommand\vf{\varphi}

\newcommand\D{\Delta}

\renewcommand\Xi{\Xi}
\renewcommand\Pi{\Pi} 

\newcommand\Up{\Upsilon}

\newcommand\Om{\Omega}

\newcommand\vG{\varGamma}

\newcommand\vL{\varLambda}

\newcommand\Bb{\boldsymbol\beta}

\newcommand{\dis}{\displaystyle}

\newcommand\wh{\widehat}
\newcommand\wt{\widetilde}

\newcommand\ra{\rightarrow} 
\newcommand\lar{\leftarrow}

\newcommand\LRa{\Leftrightarrow}

\newcommand\op{\oplus}

\newcommand\lan{\langle}
\newcommand\ran{\rangle}

\newcommand\Id{\operatorname{Id}}

\newcommand\Top{\operatorname{Top}}

\newcommand\rad{\operatorname{rad}}

\newcommand\opp{\operatorname{opp}}

\newcommand\Fgl{\mathfrak{gl}}
\newcommand\Fsl{\mathfrak{sl}}

\newcommand\ev{\mathbf{{ev}}}

\newcommand{\isom}{\,\raise2pt\hbox{$\underrightarrow{\sim}$}\,}

\newcounter{ichi}
\newcommand{\roi}{\roman{ichi}}
\newcounter{ni}
\newcommand{\roii}{\roman{ni}}
\newcounter{san}
\newcommand{\roiii}{\roman{san}}
\newcounter{yon}
\newcommand{\roiv}{\roman{yon}}
\newcounter{go}
\newcommand{\rov}{\roman{go}}
\newcounter{roku}
\newcounter{nana}
\newcounter{hachi}
\newcounter{kyu}
\begin{document}

\setlength{\baselineskip}{4.9mm}
\setlength{\abovedisplayskip}{4.5mm}
\setlength{\belowdisplayskip}{4.5mm}


\renewcommand{\theenumi}{\roman{enumi}}
\renewcommand{\labelenumi}{(\theenumi)}
\renewcommand{\thefootnote}{\fnsymbol{footnote}}
\renewcommand{\thefootnote}{\fnsymbol{footnote}}
\parindent=20pt


\setcounter{section}{-1}




\address{R. Kodera : Department of Mathematics, Graduate School of Science,  Kobe University, 
	1-1, Rokkodai, Nada-ku, Kobe 657-8501, Japan}
\email{kodera@math.kobe-u.ac.jp}
		
\address{K. Wada : Department of Mathematics, Faculty of Science, Shinshu University, 
		Asahi 3-1-1, Matsumoto 390-8621, Japan}
		
\email{wada@math.shinshu-u.ac.jp}



\medskip
\begin{center}
{\large \textbf{Finite dimensional simple modules of $(q,\bQ)$-current algebras}}  
\\
\vspace{1cm}
Ryosuke Kodera and Kentaro Wada 
\\[1em]
\end{center}


\title{} 
\maketitle 

\markboth{R. Kodera and K. Wada}{ Finite dimensional simple modules of $(q,\bQ)$-current algebras}



\begin{abstract}
The $(q, \mathbf{Q})$-current algebra associated with the general linear Lie algebra 
was introduced by the second author in the study of representation theory of cyclotomic $q$-Schur algebras. 
In this paper, we study the  $(q, \mathbf{Q})$-current algebra $U_q(\mathfrak{sl}_n^{\langle \mathbf{Q} \rangle}[x])$ 
associated with the special linear Lie algebra $\mathfrak{sl}_n$. 
In particular, we classify  finite dimensional simple 
$U_q(\mathfrak{sl}_n^{\langle \mathbf{Q} \rangle}[x])$-modules. 
\end{abstract}

\tableofcontents 

\section{Introduction} 
\para 
The $(q, \bQ)$-current algebra associated with the general linear Lie algebra 
was introduced in \cite{W16} to study the representation theory of cyclotomic $q$-Schur algebras. 
(In fact, the algebra introduced in \cite{W16} is isomorphic to a 
$(q, \bQ)$-current algebra, considered in this paper, with special parameters. 
See Appendix \ref{section gl} for these connections.) 
We expected that the $(q, \bQ)$-current algebra has good properties like  quantum groups. 

\para 
In this paper, we study the  $(q, \bQ)$-current algebra $U_q(\Fsl_n^{\lan \bQ \ran}[x])$ 
associated with the special linear Lie algebra $\Fsl_n$. 
The $(q, \bQ)$-current algebra $U_q(\Fsl_n^{\lan \bQ \ran}[x])$ has parameters 
$q \in \CC^{\times}$ and $\bQ=(Q_1,Q_2,\dots, Q_{n-1}) \in \CC^{n-1}$. 

In the case where $q=1$, 
the algebra $U_1(\Fsl_n^{\lan \bQ \ran}[x])$ is isomorphic to the universal enveloping algebra of the deformed current Lie algebra 
$\Fsl_n^{\lan \bQ \ran}[x]$ given in \cite{W18} under avoiding some ambiguities of signs  
(see Remark \ref{Remark q=1} (\roii)). 
We remark that the deformed current Lie algebra $\Fsl_n^{\lan \bQ \ran}[x]$ 
is isomorphic to the polynomial current Lie algebra $\Fsl_n [x]$ if $\bQ =(0, \dots, 0)$. 

On the other hand, 
in the case where $\bQ= \mathbf{0} = (0,\dots,0)$, 
the algebra $U_q (\Fsl_n^{\lan \mathbf{0} \ran}[x])$ is a subalgebra 
of the quantum loop algebra $U_q(L \Fsl_n)$. 
This connection corresponds to the fact that 
the polynomial current Lie algebra $\Fsl_n[x] = \Fsl_n \otimes \CC[x]$ 
is a subalgebra of the loop Lie algebra $L \Fsl_n=\Fsl_n \otimes \CC [x,x^{-1}]$ in the natural way 
(see Remark \ref{Remark shift inj} (\roi)). 
By using the explicit description in \cite{FT} for the coproduct of $U_q(L\Fsl_n)$ under Drinfeld's new generators, 
we see that the coproduct of $U_q(L\Fsl_n)$ induces the coproduct of $U_q (\Fsl_n^{\lan \mathbf{0} \ran}[x])$ by the restriction 
(see Proposition \ref{Prop coprod Uq0}).  

In general, 
we prove that the algebra $U_q (\Fsl_n^{\lan \bQ \ran}[x])$ 
is a subalgebra of a quotient of a shifted quantum affine algebra $\cU_{\bb,0}$ introduced in \cite{FT}. 
The quotient is obtained by regarding some central elements in $\cU_{\bb,0}$ as scalars 
depending on the parameters $Q_1, \dots, Q_{n-1}$. 
(see Proposition \ref{Prop ThetaQ} for details).  
Then, by applying an analogy of the argument in \cite{FT}, 
we have the following theorem. 

\begin{thm}[{Proposition \ref{Prop inj Q to 0}, Theorem \ref{Thm DrQ} and Proposition \ref{Prop coass}}]\ 
\begin{enumerate}
\item 
There exist  injective algebra homomorphisms 
$\iota_+^{\lan \bQ \ran}$ and $\iota_-^{\lan \bQ \ran}$ from  
$U_q(\Fsl_n^{\lan \bQ \ran}[x])$ to $U_q(\Fsl_n^{\lan \mathbf{0} \ran}[x])$. 

\item 
The algebra $U_q (\Fsl_n^{\lan \bQ \ran}[x])$ is a right (resp. left) coideal subalgebra of 
$U_q (\Fsl_n^{\lan \mathbf{0} \ran}[x])$ through the injection $\iota_-^{\lan \bQ \ran}$ (resp. $\iota_+^{\lan \bQ \ran}$). 
\end{enumerate}
\end{thm}

\para 
The goal of this paper is to classify the finite dimensional simple modules of $U_q (\Fsl_n^{\lan \bQ \ran}[x])$ 
in the case where $q$ is not a root of unity. 
The algebra $U_q (\Fsl_n^{\lan \bQ \ran}[x])$ has a triangular decomposition (Theorem \ref{Theorem PBW}), 
and we see that every finite dimensional simple $U_q (\Fsl_n^{\lan \bQ \ran}[x])$-module is a highest weight module 
in the usual manner. 
Thus, it is enough to classify the highest weights such that the corresponding simple highest weight modules are finite dimensional. 
The highest weight for a highest weight $U_q(\Fsl_n^{\lan \bQ \ran}[x])$-module 
is described by an element of $(\CC^{\times} \times \prod_{t >0} \CC)^{n-1}$. 
For $\bu=((\la_i, (u_{i,t})_{t>0}))_{1\leq i \leq n-1} \in (\CC^{\times} \times \prod_{t>0} \CC)^{n-1}$, 
we denote the simple highest weight $U_q(\Fsl_n^{\lan \bQ \ran}[x])$-module of the highest weight $\bu$ by $L(\bu)$ 
(see \S \ref{section hw} for details). 

In order to describe the highest weights for finite dimensional simple $U_q (\Fsl_n^{\lan \bQ \ran}[x])$-modules, 
we prepare some combinatorics as follows. 

For $t,k \in \ZZ_{>0}$, we define the symmetric polynomial $p_t(q)(x_1,x_2,\dots, x_k)$ with variables $x_1,x_2,\dots, x_k$ by 
\begin{align*}
p_t(q) (x_1,x_2,\dots, x_k) := \sum_{\la \vdash t \atop \ell(\la) \leq k} q^{- \ell (\la)} (q-q^{-1})^{\ell (\la) -1} m_{\la} (x_1,x_2,\dots, x_k),  
\end{align*}
where we denote by $\la \vdash t$ if $\la$ is a partition of $t$, 
denote by $\ell(\la)$  the length of $\la$ and 
denote by $m_{\la}(x_1,x_2, \dots, x_k)$  the monomial symmetric polynomial associated with $\la$. 
For $t,k \in \ZZ_{>0}$ and $Q, \b \in \CC^{\times}$, 
we also define the symmetric polynomial $p_t^{\lan Q \ran} (q; \b)(x_1,x_2,\dots, x_k)$ by 
\begin{align*}
&p_t^{\lan Q \ran} (q;\b) (x_1,x_2, \dots, x_k)  
\\
&:=p_t(q) (x_1,x_2, \dots,x_k) + \wt{\b} Q^{-t} + (q-q^{-1}) \sum_{z=1}^{t-1} \wt{\b} Q^{-t+z} p_z (q) (x_1,x_2, \dots,x_k), 
\end{align*}
where we put $\wt{\b} = (q-q^{-1})^{-1} (1 - \b^{-2})$. 
We remark that, 
in the case where $q=1$, 
the polynomial $p_t(1)(x_1,x_2,\dots, x_k)$ coincides with 
the power sum symmetric polynomial of degree $t$. 
We also remark that, 
in the case where $\b = \pm 1$, 
we have 
$p_t^{\lan Q \ran}(q; \pm 1)(x_1,x_2,\dots, x_k) = p_t(q)(x_1,x_2,\dots, x_k)$. 

Let $\CC[x]$ be the polynomial ring over $\CC$ with an indeterminate variable $x$. 
For $\vf \in \CC[x]$, we denote the leading coefficient of $\vf$ by $\b_{\vf}$. 
Then we define a map 
$\bu^{\lan Q \ran} : \CC[x] \setminus \{0\} \to \CC^{\times} \times \prod_{t >0} \CC$ by 
\begin{align*}
\bu^{\lan Q \ran} (\vf ) = \begin{cases} 
		(\b_{\vf}, (0)_{t>0}) & \text{ if } Q=0 \text{ and } \deg \vf=0, 
		\\
		(\b_{\vf} q^{\deg \vf}, (p_t(q)(\g_1,\g_2,\dots, \g_k))_{t>0}) & \text{ if } Q=0 \text{ and } \deg \vf >0, 
		\\
		(\b_{\vf}, (\wt{\b}_{\vf} Q^{-t})_{t>0}) & \text{ if } Q \not=0 \text{ and } \deg \vf=0, 
		\\
		(\b_{\vf} q^{\deg \vf}, (p_t^{\lan Q \ran} (q; \b_{\vf}) (\g_1,\g_2,\dots, \g_k))_{t>0}) 
			& \text{ if } Q \not=0 \text{ and } \deg \vf >0 
	\end{cases}
\end{align*}
for $\vf = \b_{\vf} (x-\g_1)(x-\g_2) \dots (x-\g_k) \in \CC[x] \setminus \{0\}$. 
For $Q \in \CC$, set 
\begin{align*}
\CC[x]^{\lan Q \ran} = \begin{cases}
		 \{ \vf \in \CC[x] \setminus \{0 \} \mid \b_{\vf} = \pm 1 \} & \text{ if } Q=0, 
		 \\
		 \{ \vf \in \CC[x] \setminus \{0\} \mid \b_{\vf}^{-2} Q^{-1} \text{ is not a root of } \vf\} 
		 &\text{ if } Q \not=0.
	\end{cases}
\end{align*}
Then we have the following classification of the isomorphism classes of finite dimensional simple 
$U_q (\Fsl_n^{\lan \bQ \ran}[x])$-modules. 
\begin{thm}[{Theorem \ref{Thm  slnQ}}]  
There exists the bijection between $\prod_{1 \leq i \leq n-1} \CC[x]^{\lan Q_i \ran}$ 
and the isomorphism classes of finite dimensional simple 
$U_q( \Fsl_n^{\lan \bQ \ran}[x])$-modules 
given by $(\vf_i)_{1\leq i \leq n-1} \mapsto L((\bu^{\lan Q_i \ran}(\vf_i))_{1 \leq i \leq n-1})$. 
\end{thm}

We remark that 
the simple highest weight module 
$L((\bu^{\lan Q_i \ran}(\vf_i))_{1 \leq i \leq n-1})$ is finite dimensional 
even if $\vf_i \not\in \CC[x]^{\lan Q_i \ran}$ for some $i$ such that $Q_i \not=0$ 
although it is infinite dimensional if $\vf_i \not\in \CC[x]^{\lan 0 \ran}$ for some $i$ such that $Q_i=0$. 
In the case where $Q \not=0$, 
the map $\bu^{\lan Q \ran} : \CC[x] \setminus \{0\} \ra \CC^{\times} \times \prod_{t>0} \CC$ is not injective, 
and we have the following proposition. 

\begin{prop}[{Proposition \ref{Lemma uQ vf = UQ vf'}}]
For $\vf, \vf' \in \CC[x] \setminus \{0\}$ such that $\deg \vf \geq \deg \vf'$, 
we have that 
$\bu^{\lan Q \ran}(\vf) = \bu^{\lan Q \ran}(\vf')$ 
if and only if 
\begin{align*}
\vf = q^{- (\deg \vf - \deg \vf')} \vf' \prod_{z=1}^{\deg \vf - \deg \vf'} ( x - q^{- 2 (z-1)} \b_{\vf}^{-2} Q^{-1}).
\end{align*}
\end{prop}
Thanks to this proposition, 
we can take the set $\prod_{1 \leq i \leq n-1} \CC[x]^{\lan Q _i \ran}$ 
as an index set for the  isomorphism classes of finite dimensional simple 
$U_q( \Fsl_n^{\lan \bQ \ran}[x])$-modules. 

We also remark that, in the case where $\bQ=\mathbf{0}=(0,\dots,0)$, 
the algebra $U_q (\Fsl_n^{\lan \mathbf{0} \ran}[x])$ is a subalgebra of the quantum loop algebra 
$U_q(L\Fsl_n)$, 
and the argument to classify finite dimensional simple $U_q (\Fsl_n^{\lan \mathbf{0} \ran}[x])$-modules 
is essentially same as the argument for $U_q (L\Fsl_n)$ given in \cite{CP91} and \cite{CP-book}. 
However, in the case where $\bQ \not=(0,\dots,0)$, we need more careful treatments.   

\para 
In the theory of quantum loop algebras and shifted quantum affine algebras, 
we usually use generating functions for generators. 
In order to describe the corresponding statements for $U_q (\Fsl_n^{\lan \bQ \ran}[x])$, we need other generators 
$\Psi_{i,t}^+ \in U_q (\Fsl_n^{\lan \bQ \ran}[x])$  ($ 1 \leq i \leq n-1, t \geq - b_i$)
defined by \eqref{Def Psi Q=0} and \eqref{Def Psi Q not=0}. 
We consider the generating function  $\Psi_i^+ (\w) = \sum_{t \geq - b_i} \Psi_{i,t}^+ \w^t$. 
We also define a map $\flat : \CC[x] \ra \CC [\w]$ ($\vf \mapsto \vf^{\flat} (\w)$) by 
\begin{align*} 
\vf^{\flat}(\w) =  (1 - \g_1 \w) (1- \g_2 \w) \dots (1 - \g_k \w) 
\end{align*} 
if $\vf = \b_{\vf} (x - \g_1) (x-\g_2) \dots (x-\g_k)$.  
Then we have the following corollary. 
\begin{cor}[{Corollary \ref{Cor Psi+ Q=0} and Corollary \ref{Cor Psi+ Q not=0}}] 
\label{Cor Intro Psi+}
For $(\vf_i)_{1\leq i \leq n-1} \in \prod_{1\leq i \leq n-1} \CC[x]^{\lan Q_i \ran}$, 
let $v_0$ be a highest weight vector 
of $L((\bu^{\lan Q_i \ran}(\vf_i))_{1\leq i \leq n-1})$. 
Then we have 
\begin{align*}
\Psi_i^+(\w) \cdot v_0 
= \begin{cases} 
	\dis \b_{\vf_i} q^{\deg \vf_i} \frac{\vf_i^{\flat} (q^{-2} \w)}{\vf_i^{\flat} (\w)}  v_0 
	& \text{ if } Q_i=0, 
	\\
	\dis 
	q^{\deg \vf_i} \frac{\vf_i^{\flat} (q^{-2} \w)}{\vf_i^{\flat} (\w)} (\b_{\vf_i}^{-1} - Q_i \b_{\vf_i} \w^{-1} ) v_0 & \text{ if } Q_i \not=0
\end{cases}
\end{align*}
for $i=1,2, \dots, n-1$. 
\end{cor}

We remark that Corollary \ref{Cor Intro Psi+} is an analogue of the statement 
for shifted Yangians given in \cite[Corollary 7.10]{BK} and \cite[Theorem 3.5]{{KTWWY}}.

\para 
After writing the first version of this paper, 
Alexander Tsymbaliuk informed us that 
he and Michael Finkelberg obtained a classification of finite dimensional simple modules of shifted quantum affine algebras of type $A$. 
Unfortunately, 
their work is unpublished. 
\\

{\bf Acknowledgements:} 
The authors are grateful to Hiraku Nakajima and Yoshihisa Saito 
for their useful suggestions and discussions. 
The authors are also grateful to Alexander Tsymbaliuk 
for his valuable comments on the first version of the paper. 
The first author was supported by JSPS KAKENHI Grant Number JP17H06127 and JP18K13390. 
The second author was supported by JSPS KAKENHI Grant Number JP16K17565.


\section{The  $(q, \bQ)$-current algebra $U_q(\Fsl_n^{\lan \bQ \ran}[x])$}
\label{section qQCA}

In this section, we give a definition of the $(q, \bQ)$-current algebra $U_q (\Fsl_n^{\lan \bQ \ran}[x])$ associated with 
the special linear Lie algebra $\Fsl_n$. 
We also give some basic properties of $U_q(\Fsl_n^{\lan \bQ \ran}[x])$. 
\para 
For $v \in \CC^{\times}$ and  any elements $x,y$ of an associative algebra over $\CC$, we put 
$[x,y]_v = xy - v yx$. 
In the case where $v=1$, we denote $[x,y]_1 = xy-yx$ by $[x,y] $ simply. 

Put $I=\{1,2, \dots, n-1\}$. 
Let $A=(a_{ij})_{i,j \in I}$ 
be the Cartan matrix of type $A_{n-1}$, 
namely we have $a_{ii}=2$, $a_{i, i\pm1}=-1$ and $a_{ij}=0$ if $j \not=i, i \pm1$. 

Take $q \in \CC^{\times}$.  
Put $[k] = (q-q^{-1})^{-1} (q^k - q^{-k})$ for $k \in \ZZ$, 
and 
$[k] ! = [k] [k-1] \dots [1]$ for $k \in \ZZ_{>0}$ with $[0]!=1$.

We define the $(q, \bQ)$-current algebra $U_q (\Fsl_n^{\lan \bQ \ran}[x])$ 
associated with 
the special linear Lie algebra $\Fsl_n$ as follows. 
\begin{definition}
For $ q \in \CC^{\times}$ and $\bQ=(Q_1,Q_2, \dots, Q_{n-1}) \in \CC^{I}$, 
we define an associative algebra $U_q( \Fsl_n^{\lan \bQ \ran}[x])$ 
over $\CC$ by the following generators and defining relations: 
\begin{description}
\item[Generators] 
$X_{i,t}^{\pm}$, $J_{i,t}$, $K_i^{\pm}$ ($i \in I$, $t \in \ZZ_{\geq 0}$), 

\item[Defining relations] 
\begin{align*}
&\tag{Q1-1} 
	[K_i^+, K_j^+] = [K_i^+, J_{j,t}] = [J_{i,s}, J_{j,t}]=0, 
	\\ &
	\tag{Q1-2}
	K_i^+ K_i^- = 1 = K_i^- K_i^+, 
	\quad 
	(K_i^-)^2 = 1 - (q-q^{-1}) J_{i,0},  
\\
&\tag{Q2} 
	X_{i,t+1}^{+} X_{j,s}^{+} - q^{a_{ij}} X_{j,s}^+ X_{i,t+1}^+ = q^{a_{ij}} X_{i,t}^+ X_{j,s+1}^+ - X_{j,s+1}^+ X_{i,t}^+, 
\\
&\tag{Q3} 
	X_{i,t+1}^{-} X_{j,s}^{-} - q^{-a_{ij}} X_{j,s}^- X_{i,t+1}^- = q^{-a_{ij}} X_{i,t}^- X_{j,s+1}^- - X_{j,s+1}^- X_{i,t}^-, 
\\
&\tag{Q4-1} 
	K_i^+ X_{j,t}^+ K_i^- = q^{a_{ij}} X_{j,t}^+, 
	\\ 
	& \tag{Q4-2} 
	q^{a_{ij}} J_{i,0} X_{j,t}^+ - q^{- a_{ij}} X_{j,t}^+ J_{i,0} = [a_{ij}] X_{j,t}^+, 
	\\
	& \tag{Q4-3} 
	[J_{i,s+1}, X_{j,t}^+] = q^{a_{ij}} J_{i,s} X_{j,t+1}^+ - q^{-a_{ij}} X_{j,t+1}^+ J_{i,s},  
\\
&\tag{Q5-1} 
	K_i^+ X_{j,t}^- K_i^- = q^{-a_{ij}} X_{j,t}^-, 
	\\
	& \tag{Q5-2}
	q^{-a_{ij}} J_{i,0} X_{j,t}^- - q^{ a_{ij}} X_{j,t}^- J_{i,0} = [-a_{ij}] X_{j,t}^-, 
	\\
	& \tag{Q5-3} 
	[J_{i,s+1}, X_{j,t}^-] = q^{- a_{ij}} J_{i,s} X_{j,t+1}^- - q^{a_{ij}} X_{j,t+1}^- J_{i,s},  
\\
&\tag{Q6}
	[X_{i,t}^+, X_{j,s}^-] = \d_{i,j} K_i^+ (  J_{i,s+t} - Q_i  J_{i,s+t+1}), 
\\
&\tag{Q7} 
	[X_{i,t}^+, X_{j,s}^+] =0  \text{ if } j \not=i,i \pm 1, 
	\\& 
	X_{i\pm 1,u}^+ (X_{i,s}^+ X_{i,t}^+ + X_{i,t}^+ X_{i,s}^+) + (X_{i,s}^+ X_{i,t}^+ + X_{i,t}^+ X_{i,s}^+) X_{i \pm 1, u}^+  
	\\ & = (q+q^{-1}) (X_{i,s}^+ X_{i \pm 1,u}^+ X_{i,t}^+ + X_{i,t}^+ X_{i \pm 1, u}^+ X_{i,s}^+), 
\\
&\tag{Q8} 
	[X_{i,t}^-, X_{j,s}^-] =0  \text{ if } j \not=i,i \pm 1, 
	\\& 
	X_{i\pm 1,u}^- (X_{i,s}^- X_{i,t}^- + X_{i,t}^- X_{i,s}^-) + (X_{i,s}^- X_{i,t}^- + X_{i,t}^- X_{i,s}^-) X_{i \pm 1, u}^-  
	\\ & = (q+q^{-1}) (X_{i,s}^- X_{i \pm 1,u}^- X_{i,t}^- + X_{i,t}^- X_{i \pm 1, u}^- X_{i,s}^-). 
\end{align*}
\end{description}
\end{definition}
We call $U_q(\Fsl_n^{\lan \bQ \ran}[x])$ the $(q, \bQ)$-current algebra associated with $\Fsl_n$. 
We denote $U_q (\Fsl_n^{\lan \bQ \ran}[x])$ by $U_q^{\lan \bQ \ran}$ simply 
unless there is any confusion. 

\begin{remarks}
\label{Remark q=1}
\begin{enumerate}
\item 
If $q \not=1$, 
the relation (Q4-2) (resp. (Q5-2)) follows from the relations (Q1-2) and (Q4-1) 
(resp. (Q1-2) and (Q5-1)).  

\item 
In the case where $q=1$, 
we see easily that 
$U_1 (\Fsl_n^{\lan \bQ \ran}[x])/\lan K_i^+ -1 \mid i \in I\ran$ is isomorphic to the universal enveloping algebra 
of the deformed current Lie algebra 
$\Fsl_n^{\lan \bQ \ran}[x]$ given in \cite[Definition 1.1]{W18}, 
where $\lan K_i^+  -1 \mid i \in I \ran$ is the two-sided ideal of $U_1 (\Fsl_n^{\lan \bQ \ran}[x])$ 
generated by 
$\{ K_i^+ -1 \mid i \in I\}$. 
Under this isomorphism, 
the generators $X_{i,t}^{\pm}$ and $J_{i,t}$ of $U_1 (\Fsl_n^{\lan \bQ \ran}[x])$ 
correspond to the generators of the enveloping algebra of $\Fsl_n^{\lan \bQ \ran}[x]$ denoted by the same symbols respectively.  
We note that $\Fsl_n^{\lan \bQ \ran}[x]$ is isomorphic to the polynomial current Lie algebra $\Fsl_n [x]$ 
if $\bQ=(0,\dots,0)$. 
\end{enumerate}
\end{remarks}

>From the defining relations, 
we can easily check the following lemma. 

\begin{lem}
\label{Lemma dag}
There exists the algebra anti-involution 
$\dag : U_q^{\lan \bQ \ran} \ra U_q^{\lan \bQ \ran}$  
such that 
$\dag(X_{i,t}^{\pm}) = X_{i,t}^{\mp}$, 
$\dag(J_{i,t}) = J_{i,t}$ 
and $\dag(K_i^{\pm})  = K_i^{\pm}$ 
for $i \in I$ and $t \in \ZZ_{\geq 0}$. 
\end{lem}

\para 
The relation (Q1-2) implies that  
\begin{align}
\label{Ji0}
J_{i,0} = \frac{1- (K_i^-)^2}{q-q^{-1}}
\end{align}
if $q^2 \not= 1$. 
By the relations (Q4-2), (Q4-3), (Q5-2) and (Q5-3), 
we have 
\begin{align}
\label{Ji1 Xitpm}
[J_{i,1},  X_{i,t}^{\pm} ]  = \pm [2] X_{i,t+1}^{\pm}. 
\end{align}
This implies that   
\begin{align}
\label{Xit+1pm} 
X_{i,t+1}^{\pm} = \pm \frac{1}{[2]}  [J_{i,1}, X_{i,t}^{\pm}] 
\end{align}
if $q^2 \not= -1$. 
The relations (Q1-2) and (Q6) imply that 
\begin{align}
\label{Jit+1}
J_{i,t+1} = \begin{cases}
		K_i^- [X_{i,t+1}^+, X_{i,0}^-] & \text{ if } Q_i=0,  
		\\
		Q_i^{-1} J_{i,t} - Q_i^{-1} K_i^- [X_{i,t}^+, X_{i,0}^-] & \text{ if } Q_i \not=0.
	\end{cases}
\end{align}
Thanks to the relations \eqref{Ji0}, \eqref{Xit+1pm} and \eqref{Jit+1}, 
we have the following lemma.

\begin{lem}
\label{Lemma gen}
Assume that $q^2 \not=\pm1$. 
The algebra $U_q (\Fsl_n^{\lan \bQ \ran}[x])$ is generated by 
$X_{i,0}^{\pm}$ $J_{i,1}$ and $K_i^{\pm}$ for $i \in I$. 
\end{lem}


\para
Let $U_{q,\bQ}^{\pm}$  be a subalgebra of $U_q^{\lan \bQ \ran}$ generated by 
$X_{i,t}^{\pm}$ ($(i,t) \in I \times \ZZ_{\geq 0}$), 
and $U_{q, \bQ}^0$ be a subalgebra of $U_q^{\lan \bQ \ran}$ generated by 
$J_{i,t}$ ($(i,t) \in I \times \ZZ_{\geq 0}$) and $K_i^{\pm}$ ($i \in I$). 
>From the defining relations, 
we see that 
\begin{align}
\label{weak tri decom}
U_q^{\lan \bQ \ran} = U_{q,\bQ}^- \cdot U_{q,\bQ}^0 \cdot U_{q, \bQ}^+.
\end{align}
\para 
Through the connection with the shifted quantum affine algebra given in the next section, 
and using the PBW theorem for the quantum loop algebra in \cite{T}, 
we can obtain the PBW theorem for $U_q^{\lan \bQ \ran}$. 
In this section, 
we give only the statement of PBW theorem 
for $U_q^{\lan \bQ \ran}$, 
and a proof is given in Appendix \ref{Proof PBW}. 

\para
Let $\{\a_1, \a_2,\dots, \a_{n-1}\}$ be the 
set of simple roots of $\Fsl_n$, 
and  
\begin{align*}
\D^+ = \{\a_{i,j} := \a_i + \a_{i+1} + \dots + \a_{j-1} \mid 1 \leq i <  j \leq n\}
\end{align*} 
be the set of positive roots. 
We define a total order on $\D^+$ by 
\begin{align*}
\a_{i,j} \leq \a_{i',j'} 
\text{ if }
i <i' \text{ or } i=i', j \leq j'.
\end{align*}
We also define a  total order on $\D^+ \times \ZZ$ by 
\begin{align*}
(\b,t) \leq (\b', t') \text{ if } \b < \b' \text{ or } \b=\b', t \leq t'. 
\end{align*}
Let $H_{\geq 0}$ denote the set of all functions 
$h : \D^+ \times \ZZ_{\geq 0} \ra \ZZ_{\geq 0}$ with finite support.  

For $(\a_{i,j}, t) \in \D^+ \times \ZZ_{\geq 0}$, put 
\begin{align*}
& X_{\a_{i,j}}^+(t) 
	:= [[ \dots [[X_{j-1,0}^+, X_{j-2,0}^+]_{q}, X_{j-3,0}^+]_{q}, \dots, X_{i+1,0}^+]_{q},X_{i,t}^+]_{q}, 
\\
& X_{\a_{i,j}}^-(t) 
	:= [X_{i,t}^-, [X_{i+1,0}^-, \dots, [X_{j-3,0}^- [X_{j-2,0}^-, X_{j-1,0}^-]_q ]_q\dots ]_q ]_q. 
\end{align*}


For $h \in H_{\geq 0}$, put 
\begin{align*}
X_h^+ := \prod_{(\b,t) \in \D^+ \times \ZZ_{\geq 0}}^{\ra} X_{\b}^+ (t)^{h(\b,t)}, 
\quad 
X_h^- := \prod_{(\b,t) \in \D^+ \times \ZZ_{\geq 0}}^{\lar} X_{\b}^-(t)^{h(\b,t)}.
\end{align*}
We define a total order on $I \times \ZZ_{\geq 0}$ by 
$(i,t) \leq (i',t')$ if $i < i'$ or $i=i'$, $t \leq t'$. 
Let $H_0$ denote the set of all functions $h_0 : I \times \ZZ_{> 0} \ra \ZZ_{\geq 0}$ with finite support, 
and put 
\begin{align*}
J_{h_0} := \prod_{(i,t)\in I \times \ZZ_{ > 0}}^{\ra} J_{i,t}^{h_0(i,t)} 
\end{align*}
for $h_0 \in H_0$. 
For $\bk =(k_1,k_2,\dots,k_{n-1}) \in \ZZ^{I}$, 
put 
\begin{align*}
K^{\bk} = K_1^{k_1} K_2^{k_2} \dots K_{n-1}^{k_{n-1}}.
\end{align*}
(Note that $J_{i,0} = (q-q^{-1})^{-1} ( 1 - (K_i^-)^2)$ by the relation (Q1-1).)  
Then we have the following theorem. 


\begin{thm}
\label{Theorem PBW}
Assume that $q \not= \pm 1$, then we have the following.  
\begin{enumerate} 
\item 
The multiplication map 
\begin{align*}
U_{q, \bQ}^- \otimes U_{q,\bQ}^0 \otimes U_{q,\bQ}^+ \ra U_q^{\lan \bQ \ran} 
\end{align*}
gives an isomorphism of vector spaces. 

\item 
\begin{enumerate}
\item 
$\{X_h^+ \mid h \in H_{\geq 0}\}$  gives a $\CC$-basis of $U_{q,\bQ}^+$. 
\item 
$\{ X_h^- \mid h \in H_{\geq 0} \}$ gives a $\CC$-basis of $U_{q,\bQ}^-$. 
\item 
$\{ K^{\bk} J_{h_0} \mid \bk \in \ZZ^{I}, h_0 \in H_0\}$ gives a $\CC$-basis of $U_{q,\bQ}^0$. 
\item  
$\{X_h^- K^\bk J_{h_0} K_{h'}^+ \mid h,h' \in H_{\geq 0}, h_0 \in H_0, \bk \in \ZZ^I \}$ 
gives a $\CC$-basis of $U_q^{\lan \bQ \ran}$. 
\end{enumerate} 

\item 
\begin{enumerate} 
\item 
The algebra $U_{q,\bQ}^+$ is generated by $\{X_{i,t}^+ \mid (i,t) \in I\times \ZZ_{\geq 0}\}$
	subject to the defining relations (Q2) and (Q7). 
\item 
The algebra $U_{q,\bQ}^-$ is generated by $ \{ X_{i,t}^- \mid(i,t) \in I\times \ZZ_{\geq 0}\}$
	subject to the defining relations (Q3) and (Q8).  
\item 
The algebra $U_{q,\bQ}^0$ is generated by $\{ J_{i,t}, K_i^{\pm} \mid i \in I, t \in \ZZ_{\geq 0}\}$ 
	subject to the defining relations (Q1-1) and (Q1-2). 
\end{enumerate}
\end{enumerate} 
\begin{proof}
See Appendix \ref{Proof PBW}. 
\end{proof}
\end{thm}

\para 
In the next section, 
we give a connection with the shifted quantum affine algebras introduced in \cite{FT}. 
In particular, 
we see that $U_q^{\lan \mathbf{0} \ran}$, where $\mathbf{0}=(0,\dots,0)$,  
turns out to be a Hopf subalgebra of the quantum loop algebra $U_q(L\Fsl_n)$ associated with $\Fsl_n$. 
Then the  injective algebra homomorphisms 
$\iota^{\lan \bQ \ran}_{\pm} : U_q^{\lan \bQ \ran} \ra U_q^{\lan \mathbf{0} \ran}$ 
given in the following proposition have an important role in this paper. 

\begin{prop}
\label{Prop inj Q to 0}
Assume that $q^2 \not= \pm1$. We have the followings. 
\begin{enumerate}
\item 
There exists an injective algebra homomorphism 
$\iota^{\lan \bQ \ran}_+ : U_q^{\lan \bQ \ran} \ra U_q^{\lan \mathbf{0} \ran}$ such that 
$X_{i,t}^+ \mapsto X_{i,t}^+ - Q_i X_{i,t+1}^+$, \, 
$X_{i,t}^- \mapsto X_{i,t}^-$, \, 
$K_i^{\pm} \mapsto K_i^{\pm}$, \,  
$J_{i,t} \mapsto J_{i,t}$. 

\item 
There exists an injective algebra homomorphism 
$\iota^{\lan \bQ \ran}_- : U_q^{\lan \bQ \ran} \ra U_q^{\lan \mathbf{0} \ran}$ such that 
$X_{i,t}^+ \mapsto X_{i,t}^+ $, \, 
$X_{i,t}^- \mapsto X_{i,t}^- - Q_i X_{i,t+1}^-$, \, 
$K_i^{\pm} \mapsto K_i^{\pm}$, \, 
$J_{i,t} \mapsto J_{i,t}$. 
\end{enumerate}

\begin{proof}
We can prove the well-defindness of the homomorphisms $\iota_{\pm}^{\lan \bQ \ran}$ 
by checking the defining relations directly. 

In order to show the injectivity, 
it is enough to show that the restrictions of $\iota_{\pm}^{\lan \bQ \ran}$ to each subalgebras 
$U_{q,\bQ}^{+}$, $U_{q,\bQ}^-$ and $U_{q,\bQ}^{0}$ are injective 
thanks to Theorem \ref{Theorem PBW}. 
By Theorem \ref{Theorem PBW} (\roiii) and the definitions of $\iota_{\pm}^{\lan \bQ \ran}$, 
it is clear that the restrictions 
$\iota_+^{\lan \bQ \ran}|_{U_{q,\bQ}^-}$, 
$\iota_-^{\lan \bQ \ran}|_{U_{q,\bQ}^+}$ 
and $\iota_{\pm}^{\lan \bQ \ran}|_{U_{q,\bQ}^0}$ are injective. 
We prove the restriction $\iota_+^{\lan \bQ \ran}|_{U_{q,\bQ}^+}$ is injective. 
The injectivity of $\iota_-^{\lan \bQ \ran}|_{U_{q,\bQ}^-}$ is similar. 

Let $U^+$ be the associative algebra generated by 
$X_{i,t}^+$ ($(i,t) \in I \times \ZZ_{\geq 0}$) 
subject to the defining relations (Q2) and (Q7). 
Then both $U_{q,\bQ}^+$ and $U_{q, \mathbf{0}}^+$ are isomorphic to the algebra $U^+$ 
in natural way by Theorem \ref{Theorem PBW} (\roiii)-(a), 
and the homomorphism $\iota_+^{\lan \bQ \ran}|_{U_{q,\bQ}^+}$ 
coincides with the endomorphism 
$\iota^+ : U^+ \ra U^+$ ($X_{i,t}^+ \mapsto X_{i,t}^+ - Q_i X_{i,t+1}^+$). 
We prove that the endomorphism $\iota^+$ is injective. 
We see that the algebra $U^+$ becomes a $\ZZ$-graded algebra by putting 
$\deg (X_{i,t}^+)=t$ for $(i,t) \in I \times \ZZ_{\geq 0}$. 
Then, for $x \in U^+$, we can write $x= \sum_{t \geq s} x_t$ with $\deg (x_t)=t$ 
for some $s \in \ZZ_{\geq 0}$. 
>From the definition of $\iota^+$, 
we see that $\iota (x) = x_s + \sum_{t >s} x'_t$ with $\deg(x'_t)=t$. 
Thus, we have $\iota^+(x) \not=0$ if $x \not=0$, and $\iota^+$ is injective. 
\end{proof}
\end{prop} 


\begin{remark}
The injections in Proposition \ref{Prop inj Q to 0} are certain modifications of ones 
in \cite[Lemma 10.18]{FT} 
(see Remarks \ref{Remark shift inj} (\roii)). 
\end{remark}

\section{A connection with the shifted quantum affine algebras} 

In this section, we give a connection between the $(q, \bQ)$-current algebra and the shifted quantum affine algebra 
introduced in \cite{FT}. 
In fact, the $(q, \bQ)$-current algebra turns out to be a subalgebra of a shifted quantum affine algebra with a suitable shift. 
We recall the definition of the shifted quantum affine algebras in \cite{FT} whose shifts are at most $1$ 
since we need only these shifts.  

\begin{definition}[{\cite{FT}}]
For $q \in \CC \setminus \{0, \pm 1\}$ and $\bb=(b_1,b_2, \dots, b_{n-1}) \in \{0,1\}^{I}$, 
the shifted quantum affine algebra $\cU_{\bb,0}$ is an associative algebra over $\CC$ generated by 
$e_{i,t}$, $f_{i,t}$ ($(i,t) \in I \times \ZZ$), 
$\psi_{i,s_i}^+$ ($(i,s_i) \in I \times \ZZ_{\geq - b_i}$), $(\psi_{i,-b_i}^+)^{-1}$, 
$\psi_{i,s}^-$ ($(i,s) \in I \times \ZZ_{\leq 0}$) and $(\psi_{i,0}^-)^{-1}$ 
subject to the following defining relations: 
\begin{align*}
\tag{U1}
&[ \psi_{i, s_i}^+, \psi_{j, t_j}^+] = [ \psi_{i, s_i}^+, \psi_{j, t}^-] 
	 = [ \psi_{i, s}^-, \psi_{j, t}^-] =0 
	\quad  (s_i \geq - b_i, \, t_j \geq - b_j, \, s,t \leq 0), 
	\\
	& \psi_{i, - b_i}^+ (\psi_{i, -b_i}^+)^{-1} = (\psi_{i, -b_i}^+)^{-1} \psi_{i, -b_i}^+ =1, 
	\quad 
	\psi_{i,  0}^- (\psi_{i,  0}^-)^{-1} = (\psi_{i,  0}^-)^{-1} \psi_{i, 0}^- =1, 
\\
\tag{U2}
&e_{i, t+1} e_{j,s} - q^{a_{ij}} e_{j,s} e_{i,t+1}  =   q^{a_{ij}} e_{i,t} e_{j, s+1} - e_{j,s+1} e_{i,t} 
	\quad (s, t \in \ZZ), 
\\
\tag{U3}
&f_{i,t+1} f_{j,s}  - q^{- a_{ij}} f_{j,s} f_{i,t+1} = q^{- a_{ij}} f_{i,t} f_{j,s+1}  -  f_{j,s+1} f_{i,t} 
	\quad (s,t \in \ZZ),
\\
\tag{U4}
&\psi_{i, -b_i}^+ e_{j,s} (\psi_{i, -b_i}^+)^{-1} = q^{a_{ij}} e_{j,s}, 
	\quad 
	\psi_{i, 0}^- e_{j,s} (\psi_{i,0}^-)^{-1}= q^{- a_{ij}} e_{j,s},
	\\ & 
	\psi_{i, t+1}^+ e_{j,s} - q^{a_{ij}} e_{j,s} \psi_{i,t+1}^+  = q^{a_{ij}} \psi_{i,t}^+ e_{j,s+1}  - e_{j,s+1} \psi_{i,t}^+ 
	\quad ( s \in \ZZ, \, t \geq - b_i), 
	\\ & 
	\psi_{i, t}^- e_{j,s-1} - q^{a_{ij}} e_{j,s-1} \psi_{i, t}^-  
		= q^{a_{ij}} \psi_{i, t-1}^- e_{j, s}  - e_{j,s} \psi_{i,t-1}^-
		\quad (s \in \ZZ, \, t \leq  0),
\\
\tag{U5}
&\psi_{i, -b_i}^+ f_{j,s} (\psi_{i, -b_i}^+)^{-1} = q^{ - a_{ij}} f _{j,s}, 
	\quad 
	\psi_{i, 0}^- f_{j,s} (\psi_{i,0}^-)^{-1}= q^{a_{ij}} f_{j,s},
	\\ & 
	\psi_{i, t+1}^+ f_{j,s} - q^{ - a_{ij}} f_{j,s} \psi_{i,t+1}^+  = q^{ - a_{ij}} \psi_{i,t}^+ f_{j,s+1}  - f_{j,s+1} \psi_{i,t}^+ 
	\quad ( s \in \ZZ, \, t \geq - b_i), 
	\\ & 
	\psi_{i, t}^- f_{j,s-1} - q^{ - a_{ij}} f_{j,s-1} \psi_{i, t}^-  
		= q^{ - a_{ij}} \psi_{i, t-1}^- f_{j, s}  - f_{j,s} \psi_{i,t-1}^-
		\quad (s \in \ZZ, \, t \leq  0),
\\
\tag{U6}
& [e_{i,t}, f _{j,s}] 
	= \d_{i,j} \begin{cases}
		\dis \frac{ \psi_{i,t+s}^+}{q-q^{-1}} & \text{ if }  s+t > 0, 
		\\[1em]
		\dis \frac{\psi^+_{i,0} - \psi_{i,0}^-}{q-q^{-1}}  & \text{ if } s+t=0,
		\\[1em]
		\dis \frac{\psi^+_{i, -1} - \psi_{i, -1}^-}{q-q^{-1}}  & \text{ if } s+t= -1 \text{ and } b_i=1, 
		\\[1em]
		\dis \frac{ - \psi_{i,t+s}^- }{q-q^{-1}} & \text{ if }  s+ t < - b_i 
		\end{cases}
	\qquad (s,t \in \ZZ),
\\
\tag{U7}
& [e_{i,t}, e_{j,s}] =0 \quad  \text{ if } j \not=i,i\pm1 \quad (s,t \in \ZZ), 
\\
&e_{i \pm 1,u} (e_{i,s} e_{i,t} + e_{i,t} e_{i,s} ) + (e_{i,s} e_{i,t} + e_{i,t} e_{i,s} ) e_{i \pm 1, u} 
	\\ & \quad 
	= (q+q^{-1}) (e_{i,s} e_{i\pm1 ,u} e_{i,t} + e_{i,t} e_{i \pm 1, u} e_{i,s}) 
		\quad (s,t,u \in \ZZ), 
\\
\tag{U8}
& [f_{i,t}, f_{j,s}] =0 \quad  \text{ if } j \not=i,i\pm1 \quad (s,t \in \ZZ), 
\\
&f_{i \pm 1,u} (f_{i,s} f_{i,t} + f_{i,t} f_{i,s} ) + (f_{i,s} f_{i,t} + f_{i,t} f_{i,s} ) f_{i \pm 1, u} 
	\\ & \quad 
	= (q+q^{-1}) (f_{i,s} f_{i\pm1 ,u} f_{i,t} + f_{i,t} f_{i \pm 1, u} f_{i,s}) 
	\quad (s,t,u \in \ZZ). 
\end{align*}
\end{definition}
We define the elements 
$\{h_{i, t}\}_{i\in I, t>0}$ 
by 
\begin{align*}
&(\psi_{i,- b_i}^+ z^{b_i})^{-1} ( \sum_{t \geq - b_i} \psi_{i,t}^+ z^{-t}) 
	= \exp ( (q-q^{-1}) \sum_{t>0} h_{i,t} z^{-t} ). 
\end{align*}
In particular, 
we have 
\begin{align*}
h_{i,1} = (q-q^{-1})^{-1}  (\psi_{i,-b_i}^+)^{-1} \psi_{i,1-b_i}^+.
\end{align*}
\begin{remarks}
\begin{enumerate}
\item 
For each $i \in I$, the element $\psi_{i,-b_i}^+ \psi_{i,0}^-$ is  a central element of $\cU_{\bb,0}$. 

\item 
In the case where $\bb=(0,\dots,0)$, 
the algebra $\cU_{0,0}/\lan \psi_{i,0}^+ \psi_{i,0}^- -1 \mid i \in I \ran$ 
is isomorphic to the quantum loop algebra $U_q (L \Fsl_n)$ associated with $\Fsl_n$. 
\end{enumerate}
\end{remarks}

\para
For $\bQ =(Q_1, Q_2,\dots,Q_{n-1}) \in \CC^{I}$, 
put $\bb_{\bQ} = (b_1,b_2,\dots, b_{n-1}) \in \{0,1\}^I$ with 
\begin{align*}
b_i = \begin{cases} 0 & \text{ if } Q_i=0, \\ 1 & \text{ if } Q_i \not=0. \end{cases}
\end{align*}
Let $\fI^{\lan \bQ\ran}$ be the two-sided ideal of 
$\cU_{\bb_{ \bQ },0}$ 
generated by 
$\{ \psi_{i,-b_i}^+ \psi_{i,0}^- + Q_i + b_i -1 \mid i \in I\}$, 
and we denote the quotient algebra $\cU_{\bb_{\bQ},0}/ \fI^{\lan \bQ \ran}$ by $\cU_{\bb_{\bQ},0}^{\lan \bQ \ran}$. 
Then we have 
\begin{align}
\label{psi 0}
\psi_{i,0}^+ = (\psi_{i,0}^-)^{- 1} \text{ if } Q_i=0, 
\text{ and }
\psi_{i,-1}^+ = - Q_i ( \psi_{i,0}^-)^{-1} \text{ if } Q_i \not=0 
\end{align}
in $\cU_{\bb_{\bQ},0}^{\lan \bQ \ran}$. 
In particular, we have $\cU_{0,0}^{\lan \mathbf 0 \ran} \cong U_q(L \Fsl_n)$ 
if $\bQ = \mathbf{0}=(0,\dots,0)$. 

\begin{prop}
\label{Prop ThetaQ}
Assume that $q \not=\pm 1$. 
There exists an injective algebra homomorphism 
\begin{align*}
&\Theta^{\lan \bQ \ran} : U_q^{\lan \bQ \ran} \ra \cU_{\bb_{\bQ},0}^{\lan \bQ \ran}, 
\\
&X_{i,t}^+ \mapsto e_{i,t}, 
\quad 
X_{i,t}^- \mapsto f_{i,t}, 
\quad 
K_i^+ \mapsto (\psi_{i,0}^-)^{-1}, 
\quad 
K_i^- \mapsto \psi_{i,0}^-, 
\\
& J_{i,t} \mapsto 
	\begin{cases} 
	(q-q^{-1})^{-1} (1- (\psi_{i,0}^-)^2) & \text{ if } t=0, 
	\\
	(q-q^{-1})^{-1} \psi_{i,t}^+ \psi_{i,0}^- & \text{ if } t>0 \text{ and } Q_i=0, 
	\\
	(q-q^{-1})^{-1} ( Q_i^{-t} - \sum_{k=1}^t Q_i^{-k} \psi_{i,t-k}^+ \psi_{i,0}^- ) & \text{ if } t >0 \text{ and } Q_i \not=0.
	\end{cases}
\end{align*} 
\begin{proof}
In order to prove the well-defindness of $\Theta^{\lan \bQ \ran}$, 
we check the relations only (Q4-3), (Q5-3) and (Q6) 
since other defining relations of $U_q^{\lan \bQ \ran}$ are clear. 
For the relation (Q4-3), we have 
\begin{align*}
&\Theta^{\lan \bQ \ran}([J_{i,s+1}, X_{j,t}^+] ) 
\\
&= \begin{cases}
	(q-q^{-1})^{-1} \psi_{i,s+1}^+ \psi_{i,0}^- e_{j,t} - e_{j,t} (q-q^{-1})^{-1} \psi_{i,s+1}^+ \psi_{i,0}^-
		& \text{ if } Q_i=0, 
	\\
	(q-q^{-1})^{-1} ( Q_i^{-(s+1)} - \sum_{k=1}^{s+1} Q_i^{-k} \psi_{i,s+1-k}^+ \psi_{i,0}^- ) e_{j,t} 
		\\ \quad 
		- e_{j,t} (q-q^{-1})^{-1} ( Q_i^{-(s+1)} - \sum_{k=1}^{s+1} Q_i^{-k} \psi_{i,s+1-k}^+ \psi_{i,0}^- ) 
		& \text{ if } Q_i \not=0
	\end{cases}
\\
&= \begin{cases}
		(q-q^{-1})^{-1} ( q^{-a_{ij}} \psi_{i,s+1}^+ e_{j,t} - e_{j,t} \psi_{i,s+1}^+) \psi_{i,0}^- 
			& \text{ if } Q_i =0, 
		\\ - (q-q^{-1})^{-1} \sum_{k=1}^{s+1} Q_i^{-k} ( q^{-a_{i j}} \psi_{i,s+1-k}^+ e_{j,t} 
			- e_{j,t} \psi_{i,s+1-k}^+ ) \psi_{i,0}^- 
			& \text{ if } Q_i \not=0
	\end{cases}
\\
&= \begin{cases}
		(q-q^{-1})^{-1} (  \psi_{i,s}^+ e_{j,t+1} - q^{-a_{ij}} e_{j,t+1} \psi_{i,s}^+) \psi_{i,0}^- 
			& \text{ if } Q_i =0, 
		\\ - (q-q^{-1})^{-1} \sum_{k=1}^{s+1} Q_i^{-k} (  \psi_{i,s-k}^+ e_{j,t+1} 
			- q^{-a_{i j}} e_{j,t+1} \psi_{i,s-k}^+ ) \psi_{i,0}^- 
			& \text{ if } Q_i \not=0
	\end{cases} 
\\
&= \begin{cases}
		q^{a_{ij}} (q-q^{-1})^{-1}  \psi_{i,s}^+ \psi_{i,0}^- e_{j,t+1}  
			- q^{-a_{ij}} e_{j,t+1} (q-q^{-1})^{-1}  \psi_{i,s}^+ \psi_{i,0}^- 
			& \text{ if } Q_i =0, 
		\\ q^{a_{ij}} (q-q^{-1})^{-1} 
			( - Q_i^{-(s+1)} \psi_{i,-1}^+ \psi_{i,0}^- -  \sum_{k=1}^{s} Q_i^{-k} \psi_{i,s-k}^+ \psi_{i,0}^- ) e_{j,t+1}
			\\  \quad 
			- q^{-a_{i j}} e_{j,t+1} (q-q^{-1})^{-1}  ( - Q_i^{-(s+1)} \psi_{i,-1}^+ \psi_{i,0}^- 
				-   \sum_{k=1}^{s} Q_i^{-k} \psi_{i,s-k}^+  \psi_{i,0}^- )
			& \text{ if } Q_i \not=0
	\end{cases}
\\
&= \Theta^{\lan \bQ \ran} ( q^{a_{ij}} J_{i,s} X_{j,t+1}^+ - q^{-a_{ij}} X_{j,t+1} J_{i,s} ), 
\end{align*}
where we note that $\psi_{i,-1}^+ \psi_{i,0}^- = - Q_i \in \cU_{\bb_{\bQ},0}^{\lan \bQ \ran}$ 
	if $Q_i \not=0$ by \eqref{psi 0}. 
The relation (Q5-3) is similar. 

We check the relation (Q6). 
If $s=t=0$, we have 
\begin{align*}
&\Theta^{\lan \bQ \ran} ([X_{i,t}^+, X_{j,s}^-]) 
\\
&= [e_{i,0}, f_{j,0}]
\\
&= \d_{i,j} (q-q^{-1})^{-1} (\psi_{i,0}^+ - \psi_{i,0}^-) 
\\
&= \d_{i,j} 
	\begin{cases}
	(\psi_{i,0}^-)^{-1}  (q-q^{-1})^{-1} (1-(\psi_{i,0}^-)^2) & \text{ if } Q_i=0, 
	\\
	(\psi_{i,0}^-)^{-1}  \big( (q-q^{-1})^{-1} (1-(\psi_{i,0}^-)^2) - Q_i (q-q^{-1})^{-1} (Q_i^{-1} - Q_i^{-1} \psi_{i,0}^+ \psi_{i,0}^- ) \big) 
	&\text{ if } Q_i \not=0
	\end{cases}
\\
&= \Theta^{\lan \bQ \ran} (\d_{i,j} K_i^+ ( J_{i,0} - Q_i J_{i,1}) )
\end{align*}
where we note that $\psi_{i,0}^+ = (\psi_{i,0}^-)^{-1} \in \cU_{\bb_{\bQ},0}^{\lan \bQ \ran}$ 
	if $Q_i=0$ by \eqref{psi 0}. 
If $s+t>0$, we have 
\begin{align*}
&\Theta^{\lan \bQ \ran} ([X_{i,t}^+, X_{j,s}^-]) 
\\
&= [e_{i,t}, f_{j,s}] 
\\
&= \d_{i,j} (q-q^{-1})^{-1} \psi_{i,t+s}^+
\\
&= \d_{i,j} \begin{cases}
	(\psi_{i,0}^-)^{-1} (q-q^{-1})^{-1} \psi_{i,t+s}^+ \psi_{i,0}^- & \text{ if } Q_i=0, 
	\\
	(\psi_{i,0}^-)^{-1} 
	\big\{ (q-q^{-1}) ( Q_i^{-(t+s)} - \sum_{k=1}^{t+s} Q_i^{-k} \psi_{i,t+s-k}^+ \psi_{i,0}^-) 
	\\ \hspace{4em}
	- Q_i (q-q^{-1}) ( Q_i^{-(t+s+1)} - \sum_{k=1}^{t+s+1} Q_i^{-k} \psi_{i,t+s+1-k}^+ \psi_{i,0}^-) 
	\big\} 
	& \text{ if } Q_i \not=0
	\end{cases}
\\
& = \Theta^{\lan \bQ \ran} ( \d_{i,j} K_i^+ (J_{i,s+t} - Q_i J_{s+t+1})).
\end{align*}

The injectivity of $\Theta^{\lan \bQ \ran}$ follows from 
Theorem \ref{Theorem PBW} and \cite[Proposition 5.1]{FT}.
\end{proof}
\end{prop}

\begin{remarks}
\label{Remark shift inj}
\begin{enumerate}
\item 
In the case where $\bQ= \mathbf{0}=(0,\dots,0)$, 
we see that $\bb_{\mathbf{0}}=0$ 
and $\cU_{0,0}^{\lan \mathbf{0} \ran} \cong U_q(L\Fsl_n)$. 
The quantum loop algebra is a $\ZZ$-graded algebra with 
$\deg(e_{i,t})=\deg(f_{i,t}) = t$, 
$\deg (\psi_{i,s}^+)=s$ and $\deg (\psi_{i,-s}^-) =-s$ 
for $i \in I$, $t \in \ZZ$ and $s \in \ZZ_{\geq 0}$. 
By the injection 
$\Theta^{\lan \mathbf{0} \ran} : 
U_q^{\lan \mathbf{0} \ran} \ra \cU_{0,0}^{\lan \mathbf{0} \ran}$, 
we can regard $U_q^{\lan \mathbf{0} \ran}$ as the subalgebra of $U_q(L\Fsl_n)$ 
generated by the elements with nonnegative degree. 
Namely, $U_q^{\lan \mathbf{0} \ran}$ is the counter part of the polynomial current Lie algebra 
$\Fsl_n[x]$ which is a Lie subalgebra of the loop Lie algebra $L\Fsl_n=\Fsl_n[x,x^{-1}]$. 

\item 
There are injective algebra homomorphisms 
\begin{align*}
&\iota'_+ : \cU_{\bb_{\bQ},0} \ra \cU_{0,0}, \quad 
	e_{i,t}  \mapsto e_{i,t} - Q_i e_{i,t+1}, \, 
	f_{i,t} \mapsto f_{i,t}, \,
	\psi_{i,t}^{\pm} \mapsto \psi_{i,t}^{\pm} - Q_i \psi_{i,t+1}^{\pm}, 
\\
&\iota'_- : \cU_{\bb_{\bQ},0} \ra \cU_{0,0}, \quad 
	e_{i,t} \mapsto e_{i,t}, \, 
	f_{i,t} \mapsto f_{i,t} - Q_i f_{i,t+1}, \,
	\psi_{i,t}^{\pm} \mapsto \psi_{i,t}^{\pm} - Q_i \psi_{i,t+1}^{\pm}, 
\end{align*}
where we put $\psi_{i,-1}^+ = \psi_{i,1}^- =0$ in $\cU_{0,0}$. 
We easily see that the injections $\iota'_{\pm}$
induce the injections $ \iota'_{\pm} : \cU_{\bb_{\bQ}}^{\lan \bQ \ran}\ra \cU_{0,0}^{\lan \mathbf{0} \ran} \cong U_q(L\Fsl_n)$. 
The injection $\iota'_{+}$ (resp. $\iota'_-$) is a certain modification of the injection $\iota_{\mu, -\mu,0}$ (resp. $\iota_{\mu,0, - \mu}$) 
given in \cite[Lemma 10.18]{FT} 
for the suitable $\mu$ through the isomorphism $\cU_{\mu,0}^{sc} \cong \cU_{0,\mu}^{sc}$. 
We need this modification 
to obtain  the injections 
from $\cU_{\bb_{\bQ},0}^{\lan \bQ \ran}$ to $\cU_{0,0}^{\lan \mathbf{0} \ran}$. 
Then, we can check  the diagram 
\begin{align*}
\xymatrix{
U_q^{\lan \bQ \ran}  \ar[r]^{\iota_{\pm}^{\lan \bQ \ran}} \ar[d]_{\Theta^{\lan \bQ \ran}}
	&U_q^{\lan \mathbf{0} \ran}  \ar[d]^{\Theta^{\lan \mathbf{0} \ran}}
\\
\cU_{\bb_{\bQ},0}^{\lan \bQ \ran} \ar[r]^{\iota'_{\pm} \hspace{2em} }
	& \cU_{0.0}^{\lan \mathbf{0}\ran} \cong U_q(L\Fsl_n)
}
\end{align*}
commutes.    
\end{enumerate}
\end{remarks}

\para 
In arguments for quantum loop algebras and shifted quantum affine algebras, 
we usually use generating functions for generators. 
In order to compare with such arguments, we prepare generating functions  
for $U_q^{\lan \bQ \ran}$ as follows. 

We define generators $\Psi_{i,t}^+ \in U_q^{\lan \bQ \ran}$  ($i \in I, t \geq - b_i$)
by 
\begin{align}
\label{Def Psi Q=0}
&\Psi_{i,0}^+ = K_i^+, 
\quad 
	\Psi_{i,t}^+ = (q-q^{-1}) K_i^+ J_{i,t} \quad  (t >0) 
\end{align} 
if $Q_i =0$, and by 
\begin{align} 
\label{Def Psi Q not=0}
\begin{split}
& \Psi_{i,-1}^+ = - Q_i K_i^+, 
	\quad 
	\Psi_{i,0}^+ = K_i^+  - (q-q^{-1}) Q_i K_i^+ J_{i,1}, 
\\
& \Psi_{i,t}^+ = (q-q^{-1}) K_i^+ (J_{i,t} - Q_i J_{i,t+1}) 
	\quad (t>0) 
\end{split}
\end{align}
if $Q_i \not=0$. 
Then, Proposition \ref{Prop ThetaQ} implies that 
$\Theta^{\lan \bQ \ran} (\Psi_{i,t}^+) = \psi_{i,t}^+$ 
for $i \in I$ and $t \geq - b_i$. 
Set 
\begin{align*}
X_i^{\pm} (\w) := \sum_{t \geq 0} X_{i,t}^{\pm} \w^t, 
\quad 
\Psi_i^+(\w) := \sum_{t \geq - b_i} \Psi_{i,t}^+ \w^t
\end{align*}
for $i \in I$, 
then we have 
\begin{align*}
\Theta^{\lan \bQ \ran} (X_i^+(\w) ) = \sum_{t \geq 0} e_{i,t} \w^t, 
\quad 
\Theta^{\lan \bQ \ran} (X_i^-(\w) ) = \sum_{t \geq 0} f_{i,t} \w^t, 
\quad 
\Theta^{\lan \bQ \ran} (\Psi_{i}^+(\w) ) 
= \sum_{t \geq - b_i} \psi_{i,t}^+ \w^t. 
\end{align*}


\section{Algebra homomorphisms $\D_r^{\lan \bQ \ran} $ and $\D_l^{\lan \bQ \ran}$}

\para 
In the case where $\bQ=\mathbf{0}=(0,\dots,0)$, 
we recall the injective homomorphism 
$\Theta^{\lan \mathbf{0} \ran} : U_q^{\lan \mathbf{0} \ran} \ra \cU_{0,0}^{\lan \mathbf{0} \ran} \cong U_q(L\Fsl_n)$ in Proposition \ref{Prop ThetaQ}. 
Let $\D : U_q(L \Fsl_n) \ra U_q (L \Fsl_n) \otimes U_q (L \Fsl_n)$ 
be the Drinfeld-Jimbo coproduct on $U_q(L \Fsl_n)$ (see \cite[Theorem 10.13]{FT} 
for the coproduct $\D$). 
Then we denote the composition of $\Theta^{\lan \mathbf{0} \ran}$ and $\D$ by 
\begin{align*}
\D^{\lan \mathbf{0} \ran} = \D \circ \Theta^{\lan \mathbf{0} \ran} 
: U_q^{\lan \mathbf{0} \ran} \ra U_q(L \Fsl_n) \otimes U_q (L\Fsl_n). 
\end{align*}
We regard $U_q^{\lan \mathbf{0} \ran} \otimes U_q^{\lan \mathbf{0} \ran}$ 
as a subalgebra of 
$U_q(L\Fsl_n) \otimes U_q (L\Fsl_n)$ through the injection $\Theta^{\lan \mathbf{0} \ran}\otimes \Theta^{\lan \mathbf{0} \ran}$. 
Then we have the following proposition. 

\begin{prop} 
\label{Prop coprod Uq0}
Assume that $q \not= \pm 1$ and $\bQ=\mathbf{0}=(0,\dots,0)$, 
then we have 
$\D^{\lan \mathbf{0} \ran} (U_q^{\lan \mathbf{0} \ran}) 
	\subset U_q^{\lan \mathbf{0} \ran} \otimes U_q^{\lan \mathbf{0} \ran}$. 
In particular, the homomorphism $\D^{\lan \mathbf{0} \ran}$ 
induces the algebra homomorphism 
\begin{align*}
\D^{\lan \mathbf{0} \ran} : U_q^{\lan \mathbf{0}\ran} 
	\ra U_q^{\lan \mathbf{0}\ran} \otimes U_q^{\lan \mathbf{0}\ran}. 
\end{align*}
Moreover, we have 
\begin{align*}
\D^{\lan \mathbf{0} \ran} (X^+_{i,0}) 
	&= 1 \otimes X_{i,0}^+ + X_{i,0}^+ \otimes K_i^+, 
	\quad 
	\D^{\lan \mathbf{0} \ran} (X^-_{i,0})= X_{i,0}^- \otimes 1 + K_i^- \otimes X_{i,0}^-,
\\
\D^{\lan \mathbf{0} \ran} (K_i^{\pm}) 
	&= K_i^{\pm} \otimes K_i^{\pm} 
\end{align*}
and 
\begin{align}
\label{DJi1}
\begin{split}
\D^{\lan \mathbf{0} \ran} (J_{i,1}) 
	&= J_{i,1} \otimes 1 + 1 \otimes J_{i,1} 
		- (q^2 - q^{-2}) X_{i,0}^+ \otimes X_{i,1}^- 
		\\ & \quad 
		+ (q-q^{-1}) \sum_{l > i+1} \wt{X}_{\a_{i+1,l}}^+(0) \otimes X_{\a_{i+1,l}}^- (1) 
		\\ & \quad 
		+ (q-q^{-1}) \sum_{k <i} q^{k+1-i} X_{\a_{k,i}}^+(0) \otimes X_{\a_{k,i}}^- (1)
		\\ & \quad 
		+ q^{-2} (q-q^{-1}) \sum_{l >i+1} [X_{i,0}^+, \wt{X}_{\a_{i+1,l}}^+(0)]_{q^3} 
			\otimes X_{\a_{i,l}}^-(1) 
		\\ & \quad 
		- (q-q^{-1}) \sum_{k<i} q^{k-i-1} [X_{i,0}^+, X_{\a_{k,i}}^+(0)]_{q^3} 
			\otimes X_{\a_{k,i+1}}^-(1)
		\\ & \quad 
		+ (q-q^{-1})^2 \sum_{l>i+1}^{k <i} q^{k-i} 
			(\wt{X}_{\a_{i,l}}^+ (0) X_{\a_{k,i}}^+(0) 
				- \wt{X}_{\a_{i+1,l}}^+ (0)  X_{\a_{k,i+1}}^+(0) ) 
			\otimes X_{\a_{k,l}}^-(1),
\end{split}
\end{align}
where 
\begin{align*}
& X_{\a_{i,j}}^+ (0) = [[\dots [X_{j-1,0}^+, X_{j-2,0}^+]_q, \dots , X_{i+1,0}^+]_q, X_{i,0}^+]_q, 
\\
& \wt{X}_{\a_{i,j}}^+(0) 
	= [[\dots [X_{j-1,0}^+, X_{j-2,0}^+]_{q^{-1}}, 
		\dots , X_{i+1,0}^+]_{q^{-1}}, X_{i,0}^+]_{q^{-1}}, 
\\
& X_{\a_{i,j}}^-(1) 
	= [X_{i,1}^-, [X_{i+1,0}^-, \dots, [X_{j-2,0}^-, X_{j-1,0}^-]_q \dots ]_q ]_q. 
\end{align*} 

\begin{proof}
By Lemma \ref{Lemma gen}, 
it is enough to check that 
$\D^{\lan \mathbf{0} \ran} (X_{i,0}^{\pm})$, 
$\D^{\lan \mathbf{0} \ran} (J_{i,1})$ 
and 
$\D^{\lan \mathbf{0} \ran} (K_i^{\pm})$ 
belong to $U_q^{\lan \mathbf{0} \ran} \otimes U_q^{\lan \mathbf{0} \ran}$  
for $i \in I$. 
Note that 
\begin{align*}
&\Theta^{\lan \mathbf{0} \ran} (K_i^+) = (\psi_{i,0}^-)^{-1} = \psi_{i,0}^+, 
\quad  
\Theta^{\lan \mathbf{0} \ran} (K_i^-) = \psi_{i,0}^-, 
\\  
& \Theta^{\lan \mathbf{0} \ran} (J_{i,1}) = (q-q^{-1})^{-1} \psi_{i,1}^+ \psi_{i,0}^- =h_{i,1}, 
\\
& \Theta^{\lan \mathbf{0} \ran} (X_{i,0}^+) =e_{i,0}, 
\quad 
\Theta^{\lan \mathbf{0} \ran} (X_{i,0}^-) = f_{i,0}. 
\end{align*} 
Moreover, we see that 
\begin{align*}
& \Theta^{\lan \mathbf{0} \ran} (X_{\a_{i,j}}^+(0)) 
	= [[\dots [e_{j-1,0}, e_{j-2,0}]_q, \dots , e_{i+1,0}]_q, e_{i,0}]_q, 
\\
& \Theta^{\lan \mathbf{0} \ran} (\wt{X}_{\a_{i,j}}^+ (0) ) 
	= [[\dots [e_{j-1,0}, e_{j-2,0}]_{q^{-1}}, \dots , e_{i+1,0}]_{q^{-1}}, e_{i,0}]_{q^{-1}}, 
\\
& \Theta^{\lan \mathbf{0} \ran} (X_{\a_{i,j}}^-(1)) 
	= [f_{i,1}, [f_{i+1,0}, \dots, [f_{j-2,0}, f_{j-1,0}]_q \dots ]_q ]_q,
\end{align*}
Then, the proposition follows from \cite[Theorem 10.13]{FT}. 
(In  \cite[Theorem 10.13]{FT}, 
the elements $\Theta^{\lan \mathbf{0} \ran} (X_{\a_{i,j}}^+(0)) $, 
$\Theta^{\lan \mathbf{0} \ran} (\wt{X}_{\a_{i,j}}^+(0)) $ 
and 
$\Theta^{\lan \mathbf{0} \ran} (X_{\a_{i,j}}^-(1)) $
are denoted by 
$\wt{E}_{ij}^{(0)}$, $E_{ij}^{(0)}$ and $F_{ji}^{(1)}$ respectively.)  
\end{proof}
\end{prop}
\begin{remark}
In fact, the statement $\D^{\lan \mathbf{0} \ran} (U_q^{\lan \mathbf{0} \ran}) \subset U_q^{\lan \mathbf{0} \ran} \otimes U_q^{\lan \mathbf{0} \ran}$ immediately follows from the RTT presentation of the quantum loop algebra and the Drinfeld-Jimbo coproduct.
The RTT presentation is crucially used in \cite[Appendix G]{FT} to derive the formula for $\Delta(h_{i,1})$ which we recalled in the above proof.
These were pointed out by Alexander Tsymbaliuk after writing the first version of this paper. 
\end{remark}

\begin{remark}
In the case where $q=\pm1$ and $\bQ=\mathbf{0}$, 
we can define the algebra homomorphism 
$\D^{\lan \mathbf{0} \ran}  
: U_q^{\lan \mathbf{0} \ran} \ra U_q^{\lan \mathbf{0} \ran} \otimes U_q^{\lan \mathbf{0} \ran}$ by 
\begin{align*}
\D^{\lan \mathbf{0} \ran} (X^+_{i,t}) 
	&= 1 \otimes X_{i,t}^+ + X_{i,t}^+ \otimes K_i^+, 
	\quad 
	\D^{\lan \mathbf{0} \ran} (X^-_{i,t})= X_{i,0}^- \otimes 1 + K_i^- \otimes X_{i,t}^-,
\\
\D^{\lan \mathbf{0} \ran} (K_i^{\pm}) 
	&= K_i^{\pm} \otimes K_i^{\pm}, 
	\quad 
	\D^{\lan \mathbf{0} \ran} (J_{i,t}) 
	= J_{i,t} \otimes 1 + 1 \otimes  J_{i,t}. 
\end{align*} 
In this case, we can check the well-defindness by direct calculations. 
\end{remark}

\para 
By \eqref{Xit+1pm}, we have 
\begin{align*}
&\D^{\lan \mathbf{0} \ran} (X_{i,1}^+) 
= \frac{1}{[2]} ( \D^{\lan \mathbf{0} \ran} (J_{i,1}) \D^{\lan \mathbf{0} \ran} (X_{i,0}^+) 
	- \D^{\lan \mathbf{0} \ran} (X_{i,0}^+) \D^{\lan \mathbf{0} \ran} (J_{i,1}) ), 
\\
&\D^{\lan \mathbf{0} \ran} (X_{i,1}^-) 
=-  \frac{1}{[2]} ( \D^{\lan \mathbf{0} \ran} (J_{i,1}) \D^{\lan \mathbf{0} \ran} (X_{i,0}^-) 
	- \D^{\lan \mathbf{0} \ran} (X_{i,0}^-) \D^{\lan \mathbf{0} \ran} (J_{i,1}) ). 
\end{align*}
Thus, Proposition \ref{Prop coprod Uq0} implies the following corollary. 


\begin{cor} 
\label{Cor D0 Xi1}
We have 
\begin{align*}
\D^{\lan \mathbf{0} \ran} (X_{i,1}^+) 
&=
1 \otimes X_{i,1}^+ 
	+  X_{i,1}^+  \otimes K_i^+ 
	+  (q-q^{-1}) X_{i,0}^+ \otimes  K_i^+ J_{i,1} 
	-  q^{-1} (q - q^{-1})^2 X_{i,0}^+ X_{i,0}^+ \otimes  X_{i,1}^- K_i^+  
	\\ & \quad 
	+ q (q-q^{-1}) \sum_{l >i+1}  \wt{X}_{\a_{i,l}}^+(0) \otimes  X_{\a_{i+1,l}}^-(1) K_i^+  
	\\ & \qquad 
	- (q-q^{-1})^2  \sum_{l >i+1}  X_{i,0}^+ \wt{X}_{\a_{i,l}}^+(0)  \otimes X_{\a_{i,l}}^-(1) K_i^+ 
	\\ & \quad 
	- q (q-q^{-1}) \sum_{k <i} q^{k-i}  X_{\a_{k,i+1}}^+(0) \otimes X_{\a_{k,i}}^- (1) K_i^+ 
	\\ & \qquad 
	-  (q-q^{-1})^2 \sum_{k<i} q^{k-i}  X_{i,0}^+ X_{\a_{k,i+1}}^+(0)   \otimes X_{\a_{k,i+1}}^-(1) K_i^+
	\\ & \quad 
	- (q-q^{-1})^2 \sum_{l>i+1}^{k <i} q^{k-i}  \wt{X}_{\a_{i,l}}^+(0) X_{\a_{k,i+1}}^+(0) \otimes X_{\a_{k,l}}^-(1) K_i^+. 
\end{align*}

\begin{align*}
\D^{\lan \mathbf{0} \ran} (X_{i,1}^-) 
	&= X_{i,1}^- \otimes 1 + K_i^+ \otimes X_{i,1}^- 
	+ q^{-1} (q-q^{-1}) \sum_{l>i+1} \wt{X}_{\a_{i+1,l}}^+ (0) K_i^+ \otimes X_{\a_{i,l}}^- (1) 
	\\
	&\quad - (q-q^{-1}) \sum_{k<i} q^{k-i} X_{\a_{k,i}}^+ (0) K_i^+ \otimes X_{\a_{k,i+1}}^- (1) 
		\\ & \quad 
		- (q-q^{-1})^2 \sum_{l > i+1}^{k < i} q^{k-i-1} 
		\wt{X}_{\a_{i+1,l}}^+(0)  X_{\a_{k,i}}^+ (0) K_i^+ \otimes X_{\a_{k,l}}^- (1).
\end{align*}
\end{cor}


\begin{remark}
The explicit form of $\D^{\lan \mathbf{0} \ran}(X_{i,1}^-)$ in Corollary \ref{Cor D0 Xi1} 
follows directly from one of $\D (f_{i,1})$  given in  \cite[Theorem 10.13]{FT} 
through the injection 
$\Theta^{\lan \mathbf{0} \ran} \otimes \Theta^{\lan \mathbf{0} \ran}$. 
\end{remark}

\para 
We recall the injective homomorphisms 
$\iota_{\pm}^{\lan \bQ \ran} : U_q^{\lan \bQ \ran} \ra U_q^{\lan \mathbf{0} \ran}$ 
in Proposition \ref{Prop inj Q to 0}, 
and we consider the algebra homomorphisms 
\begin{align*}
\D_r^{\lan \bQ \ran} := \D^{\lan \mathbf{0} \ran} \circ \iota_-^{\lan \bQ \ran} : 
	U_q^{\lan \mathbf{Q} \ran} \ra U_q^{\lan \mathbf{0} \ran} \otimes U_q^{\lan \mathbf{0} \ran}, 
\quad 
\D_l^{\lan \bQ \ran} := \D^{\lan \mathbf{0} \ran} \circ \iota_+^{\lan \bQ \ran} : 
	U_q^{\lan \mathbf{Q} \ran} \ra U_q^{\lan \mathbf{0} \ran} \otimes U_q^{\lan \mathbf{0} \ran}. 
\end{align*}
Then we have the following theorem 
by a similar argument as one in \cite[Theorem 10.20]{FT}.  

\begin{thm}
\label{Thm DrQ}
\begin{enumerate}
\item 
We have $\D_r^{\lan \bQ \ran} (U_q^{\lan \bQ \ran}) \subset \iota_-^{\lan \bQ \ran} (U_q^{\lan \bQ \ran}) \otimes U_q^{\lan \mathbf{0} \ran}$. 
In particular, the homomorphism $\D_r^{\lan \bQ \ran}$ 
induces the algebra homomorphism 
\begin{align*}
\D_r^{\lan \bQ \ran} : U_q^{\lan \bQ \ran} \ra U_q^{\lan \bQ \ran} \otimes U_q^{\lan \mathbf{0} \ran}.
\end{align*} 
Moreover, 
we have 
\begin{align*}
\D_r^{\lan \bQ \ran}(X_{i,0}^+) 
	&= 1 \otimes X_{i,0}^+ + X_{i,0}^+ \otimes K_i^+, 
\\
\D_r^{\lan \bQ \ran} (X_{i,0}^-) 
	&= X_{i,0}^- \otimes 1 + K_i^- \otimes X_{i,0}^- 
	\\ & \quad 
	- Q_i \big\{ K_i^+ \otimes X_{i,1}^- 
	+ q^{-1} (q-q^{-1}) \sum_{l> i+1} \wt{X}_{\a_{i+1,l}}^+ (0) K_i^+ \otimes X_{\a_{i,l}}^-(1)
	\\ & \hspace{4em}
	- (q-q^{-1}) \sum_{k<i} q^{k-i} X_{\a_{k,i}}^+ (0) K_i^+ \otimes X_{\a_{k,i+1}}^-(1)
	\\ & \hspace{4em}
	- (q-q^{-1})^2 \sum_{l > i+1}^{k<i} q^{k-i-1} \wt{X}_{\a_{i+1,l}}^+(0) X_{\a_{k,i}}^+(0)  K_i^+ \otimes X_{\a_{k,l}}^- (1) \big\},
\\
\D_r^{\lan \bQ \ran}(K_i^{\pm} ) &= K_i^{\pm} \otimes K_i^{\pm}, 
\end{align*}
and $\D_r^{\lan \bQ \ran} (J_{i,1})$ is given by the right-hand side of \eqref{DJi1}. 


\item 
We have $\D_l^{\lan \bQ \ran} (U_q^{\lan \bQ \ran}) 
	\subset U_q^{\lan \mathbf{0} \ran} \otimes \iota_+^{\lan \bQ \ran} (U_q^{\lan \bQ \ran})$. 
In particular, the homomorphism $\D_l^{\lan \bQ \ran}$ 
induces the algebra homomorphism 
\begin{align*}
\D_l^{\lan \bQ \ran} : U_q^{\lan \bQ \ran} \ra U_q^{\lan \mathbf{0} \ran} \otimes U_q^{\lan \bQ \ran}.
\end{align*} 
Moreover, we have 
\begin{align*}
\D_l^{\lan \bQ \ran}(X_{i,0}^+) 
	&= 1 \otimes X_{i,0}^+ + X_{i,0}^+ \otimes K_i^+ 
	\\ & \quad 
	- Q_i \big\{ 
	X_{i,1}^+  \otimes K_i^+ 
	+  (q-q^{-1}) X_{i,0}^+ \otimes  K_i^+ J_{i,1} 
	-  q^{-1} (q - q^{-1})^2 X_{i,0}^+ X_{i,0}^+ \otimes  X_{i,1}^- K_i^+  
	\\ & \hspace{4em}
	+ q (q-q^{-1}) \sum_{l >i+1}  \wt{X}_{\a_{i,l}}^+(0) \otimes  X_{\a_{i+1,l}}^-(1) K_i^+  
	\\ & \hspace{4em}
	- (q-q^{-1})^2  \sum_{l >i+1}  X_{i,0}^+ \wt{X}_{\a_{i,l}}^+(0)  \otimes X_{\a_{i,l}}^-(1) K_i^+ 
	\\ & \hspace{4em} 
	- q (q-q^{-1}) \sum_{k <i} q^{k-i}  X_{\a_{k,i+1}}^+(0) \otimes X_{\a_{k,i}}^- (1) K_i^+ 
	\\ & \hspace{4em} 
	-  (q-q^{-1})^2 \sum_{k<i} q^{k-i}  X_{i,0}^+ X_{\a_{k,i+1}}^+(0)   \otimes X_{\a_{k,i+1}}^-(1) K_i^+
	\\ & \hspace{4em} 
	- (q-q^{-1})^2 \sum_{l>i+1}^{k <i} q^{k-i}  \wt{X}_{\a_{i,l}}^+(0) X_{\a_{k,i+1}}^+(0) \otimes X_{\a_{k,l}}^-(1) K_i^+ 
	\big\}, 
\\
\D_l^{\lan \bQ \ran} (X_{i,0}^-) 
	&= X_{i,0}^- \otimes 1 + K_i^- \otimes X_{i,0}^-, 
\\
\D_l^{\lan \bQ \ran} (K_i^{\pm}) 
	&= K_i^{\pm} \otimes K_i^{\pm}, 
\end{align*}
and $\D_l^{\lan \bQ \ran} (J_{i,1})$ is given by the right-hand side of \eqref{DJi1}. 
\end{enumerate} 

\begin{proof}
We prove  the statement (\roi). 
By Lemma \ref{Lemma gen}, 
it is enough to check  the relations 
for  the  generators $X_{i,0}^{\pm}$, $J_{i,1}$ and $K_i^\pm$ ($i \in I$). 
By the definition of $\iota_-^{\lan \bQ \ran}$, 
it is clear for the generators $X_{i,0}^+$, $K_i^{\pm}$ and $J_{i,1}$ ($i \in I$).
On the other hand, we have 
\begin{align*}
\D_r^{\lan \bQ \ran}(X_{i,0}^-)
&= \D^{\lan \mathbf{0} \ran} ( X_{i,0}^- - Q_i X_{i,1}^-) 
\\
&=  X_{i,0}^- \otimes 1 + K_i^- \otimes X_{i,0}^- 
	\\ & \quad 
	- Q_i \big\{ 
		X_{i,1}^- \otimes 1 + K_i^+ \otimes X_{i,1}^- 
		+ q^{-1} (q-q^{-1}) \sum_{l>i+1} \wt{X}_{\a_{i+1,l}}^+ (0) K_i^+ \otimes X_{\a_{i,l}}^- (1) 
		\\ &\hspace{4em} 
		- (q-q^{-1}) \sum_{k<i} q^{k-i} X_{\a_{k,i}}^+ (0) K_i^+ \otimes X_{\a_{k,i+1}}^- (1) 
		\\ & \hspace{4em}
		- (q-q^{-1})^2 \sum_{l > i+1}^{k < i} q^{k-i-1} 
			\wt{X}_{\a_{i+1,l}}^+(0)  X_{\a_{k,i}}^+ (0) K_i^+ \otimes X_{\a_{k,l}}^- (1) \big\} 
\\
&= (X_{i,0}^- - Q_i X_{i,1}^-) \otimes 1 
	+ K_i^- \otimes X_{i,0}^- 
	\\ & \quad 
	\\ & \quad 
	- Q_i \big\{ 
		 K_i^+ \otimes X_{i,1}^- 
		+ q^{-1} (q-q^{-1}) \sum_{l>i+1} \wt{X}_{\a_{i+1,l}}^+ (0) K_i^+ \otimes X_{\a_{i,l}}^- (1) 
		\\ &\hspace{4em} 
		- (q-q^{-1}) \sum_{k<i} q^{k-i} X_{\a_{k,i}}^+ (0) K_i^+ \otimes X_{\a_{k,i+1}}^- (1) 
		\\ & \hspace{4em}
		- (q-q^{-1})^2 \sum_{l > i+1}^{k < i} q^{k-i-1} 
			\wt{X}_{\a_{i+1,l}}^+(0)  X_{\a_{k,i}}^+ (0) K_i^+ \otimes X_{\a_{k,l}}^- (1) \big\}. 
\end{align*}
Then we see that 
$\D_r^{\lan \bQ \ran}(X_{i,0}^-) 
\subset \iota_-^{\lan \bQ \ran} (U_q^{\lan \bQ \ran}) \otimes U_q^{\lan \mathbf{0} \ran}$ 
by the definition of $\iota_-^{\lan \bQ \ran}$, 
and we have the statement (\roi). 
The statement (\roii) is proven in a similar way. 
\end{proof}
\end{thm}

The homomorphisms $\D_r^{\lan \bQ \ran}$ and $\D_l^{\lan \bQ \ran}$ 
satisfy the following coassociativity. 

\begin{prop}[{cf. \cite[Proposition 4.14]{FKPRW}}] 
\label{Prop coass}
We have the following commutative diagrams. 
\begin{enumerate}
\item 
$\xymatrix{
U_q^{\lan \bQ \ran} \ar[rr]^{\D_r^{\lan \bQ \ran}\quad }  \ar[d]_{\D_r^{\lan \bQ \ran}}
	&& U_q^{\lan \bQ \ran} \otimes U_q^{\lan \mathbf{0} \ran} 
		\ar[d]^{\Id \otimes \D^{\lan \mathbf{0} \ran}}
\\
U_q^{\lan \bQ \ran} \otimes U_q^{\lan \mathbf{0} \ran} \ar[rr]^{\D_r^{\lan \bQ \ran} \otimes \Id\quad }
	&& U_q^{\lan \bQ \ran} \otimes U_q^{\lan \mathbf{0} \ran}  \otimes U_q^{\lan \mathbf{0} \ran}
}$

\item 
$\xymatrix{
U_q^{\lan \bQ \ran} \ar[rr]^{\D_l^{\lan \bQ \ran}\quad }  \ar[d]_{\D_l^{\lan \bQ \ran}}
	&&  U_q^{\lan \mathbf{0} \ran}  \otimes U_q^{\lan \bQ \ran}
		\ar[d]^{ \D^{\lan \mathbf{0} \ran} \otimes \Id }
\\
 U_q^{\lan \mathbf{0} \ran} \otimes U_q^{\lan \bQ \ran} \ar[rr]^{ \Id \otimes \D_l^{\lan \bQ \ran} \quad }
	&& U_q^{\lan \mathbf{0} \ran}  \otimes U_q^{\lan \mathbf{0} \ran} \otimes  U_q^{\lan \bQ \ran}
}$
\end{enumerate}
\begin{proof}
We note that the coassociativity of the coproduct $\D^{\lan \mathbf{0} \ran}$ on $U_q^{\lan \mathbf{0} \ran}$ 
follows from the coaasociativity of the Drinfeld-Jimbo coproduct $\D$ on $U_q(L\Fsl_n)$.  
By Theorem \ref{Thm DrQ} and  the coassociativity of $\D^{\lan \mathbf{0} \ran}$, 
we see that the diagram 
\begin{align*}
\xymatrix{
U_q^{\lan \bQ \ran}  \ar[rr]^{\iota_-^{\lan \bQ \ran}}  \ar[d]_{\D_r^{\lan \bQ \ran}} 
	&& U_q^{\lan \mathbf{0} \ran} \ar[rr]^{\D^{\lan \mathbf{0} \ran}\quad }  \ar[d]_{\D^{\lan \mathbf{0} \ran}}
	&& U_q^{\lan \mathbf{0} \ran} \otimes U_q^{\lan \mathbf{0} \ran} \ar[d]^{\Id \otimes \D^{\lan \mathbf{0} \ran}}
\\
U_q^{\lan \bQ \ran} \otimes U_q^{\lan \mathbf{0} \ran} \ar[rr]^{\iota_-^{\lan \bQ \ran} \otimes \Id}
	&& U_q^{\lan \mathbf{0} \ran} \otimes U_q^{\lan \mathbf{0} \ran}  \ar[rr]^{\D^{\lan \mathbf{0} \ran} \otimes \Id \quad }
	&& U_q^{\lan \mathbf{0} \ran} \otimes U_q^{\lan \mathbf{0} \ran} \otimes U_q^{\lan \mathbf{0} \ran} 
}
\end{align*} 
commutes,  
and this diagram implies (\roi). 
The commutative diagram (\roii) is proven in a similar way. 
\end{proof}
\end{prop}


\section{Evaluation homomorphisms} 
\label{section ev}
In this section, we recall the evaluation homomorphisms  from 
$U_q (\wh{\Fsl}_n)$ to  $U_q(\Fgl_n)$ given in \cite{J}, 
and we prepare some results on evaluation modules along the calculation in \cite[3.6]{CP94}.
In this section, we assume that $\CC^{\times} \ni q \not=\pm 1$. 

\para Put $\wh{I} = I \cup \{0\}$, 
and 
let $\wh{A} =(a_{ij})_{i,j\in \wh{I}}$ be the Cartan matrix of type $A_{n-1}^{(1)}$. 
Namely, the submatrix $(a_{ij})_{i,j \in I}$ is the Cartan matrix of type $A_{n-1}$, 
and we have $a_{0,0}=2$, $a_{0,1} = a_{1,0}= a_{0,n-1} = a_{n-1,0} =-1$ and $a_{0j}=a_{j0}=0$ if  $j \not=0, n-1$. 
Then the quantum affine algebra $U_q(\wh{\Fsl}_n)$ of type $A_{n-1}^{(1)}$ is an associative algebra over $\CC$ 
generated by 
$e_i$, $f_i$, $k_i^{\pm}$ ($i \in \wh{I}$) subject to the following defining relations: 
\begin{align*}
&k_i^+ k_i^- = k_i^- k_i^+=1, 
\quad 
[k_i^+, k_j^-]=0, 
\quad 
k_i^+ e_j k_i^- = q^{a_{ij}} e_j, 
\quad 
k_i^+ f_j k_i^- = q^{-a_{ij}} f_j, 
\\
&[e_i,f_j] = \d_{ij} \frac{k_i^+ - k_i^-}{q-q^{-1}}, 
\quad 
\sum_{s=0}^{1-a_{ij}} (-1)^s e_i^{(1-a_{ij} -s)} e_j e_i^{(s)} =0, 
\quad 
\sum_{s=0}^{1-a_{ij}} (-1)^s f_i^{(1-a_{ij} -s)} f_j f_i^{(s)} =0, 
\end{align*}
where we put $e_i^{(s)} = (e_i)^s / [s]!$ and $f_i^{(s)} = (f_i)^s / [s]!$ for $s \geq 0$. 
We note that $c := k_0^+ k_1^+ \dots k_{n-1}^+$ is the canonical central element of $U_q(\wh{\Fsl}_n)$. 

We also consider the quantum group $U_q(\Fgl_n)$ associated with the general linear Lie algebra $\Fgl_n$ 
which is an associative algebra over $\CC$ generated by 
$E_i$, $F_i$ ($i \in I$) and $T_j^{\pm}$ ($1\leq j \leq n$) subject to the following defining relations: 
\begin{align*}
&T_i^+ T_i^- = T_i^- T_i^+ =1, 
\quad 
[T_i^+, T_j^+]=0, 
\\
&T_i^+ E_j T_i^- = q^{\d_{i,j} - \d_{i,j+1}} E_j, 
\quad 
T_i^+ F_j T_i^- = q^{-(\d_{i,j} - \d_{i,j+1})} F_j, 
\\
& [E_i,F_j]= \d_{i,j} \frac{T_i^+ T_{i+1}^- - T_i^- T_{i+1}^+}{q-q^{-1}} , 
\\ 
&\sum_{s=0}^{1-a_{ij}} (-1)^s E_i^{(1-a_{ij}-s)} E_j E_i^{(s)}=0, 
\quad 
\sum_{s=0}^{1-a_{ij}} (-1)^s F_i^{(1-a_{ij}-s)} F_j F_i^{(s)}=0, 
\end{align*}
where we put $E_i^{(s)} = (E_i)^s/ [s]!$ and $F_i^{(s)} = (F_i)^s / [s]!$ for $s \geq 0$. 

For $\g \in \CC^{\times}$, 
we have the following evaluation homomorphism $\ev_{\g} : U_q(\wh{\Fsl}_n) \ra U_q (\Fgl_n)$. 

\begin{prop}[{\cite{J}}]
\label{Prop ev whsln gln}
For $\g \in \CC^{\times}$, 
there exists an algebra homomorphism 
$\ev_{\g} : U_q(\wh{\Fsl}_n) \ra U_q(\Fgl_n)$ such that 
\begin{align*}
&e_i \mapsto E_i, 
\quad 
f_i \mapsto F_i, 
\quad 
k_i^+ \mapsto T_i^+ T_{i+1}^-  \quad (i \in I),  
\quad 
k_0^+ \mapsto T_1^- T_n^+, 
\\
& e_0 \mapsto \g q^{-1} (T_1^+ T_n^+) [F_{n-1}, [F_{n-2}, \dots , [F_2, F_1]_{q^{-1}} \dots ]_{q^{-1}}]_{q^{-1}}, 
\\
& f_0 \mapsto (-1)^n \g^{-1} q^{n-1} (T_1^- T_n^-) [E_{n-1}, [E_{n-2}, \dots, [E_2, E_1]_{q^{-1}} \dots ]_{q^{-1}}]_{q^{-1}}.
\end{align*}
Moreover, 
the homomorphism $\ev_{\g}$ factors through the quotient algebra 
$U_q (\wh{\Fsl}_n)/ \lan c-1\ran$, 
where $\lan c-1 \ran$ is the two-sided ideal of $U_q(\wh{\Fsl}_n)$ generated by $c-1$. 
\end{prop}

It is known that the quotient algebra $U_q(\wh{\Fsl}_n)/ \lan c-1 \ran$ is isomorphic to 
the quantum loop algebra $U_q(L \Fsl_n) \cong \cU_{0,0}^{\lan \mathbf{0} \ran}$ as follows. 

\begin{prop}[{\cite{D}, \cite{B}}] 
\label{Prop iso quantum loop} 
There exists an algebra isomorphism 
$\Psi : U_q(\wh{\Fsl}_n)/ \lan c-1 \ran \ra \cU_{0,0}^{\lan \mathbf{0} \ran} \cong U_q (L \Fsl_n)$ such that 
\begin{align*}
&e_i \mapsto e_{i,0}, 
\quad 
f_i \mapsto f_{i,0}, 
\quad 
k_i^{+} \mapsto \psi_{i,0}^{+} 
\quad (i \in I), 
\quad 
k_0^{+} \mapsto \psi_{1,0}^- \psi_{2,0}^- \dots \psi_{n-1,0}^-, 
\\
& e_0 \mapsto [f_{n-1,0}, [f_{n-2,0}, \dots, [f_{2,0}, f_{1,1}]_{q^{-1}} \dots]_{q^{-1}}]_{q^{-1}} 
	(\psi_{1,0}^- \psi_{2,0}^- \dots \psi_{n-1,0}^-),  
\\
& f_0 \mapsto \mu (\psi_{1,0}^+ \psi_{2,0}^+ \dots \psi_{n-1,0}^+) 
	[e_{n-1,0}, [e_{n-2,0}, \dots, [e_{2,0}, e_{1,-1}]_{q^{-1}} \dots ]_{q^{-1}} ]_{q^{-1}}, 
\end{align*}
where 
$\mu \in \CC^{\times}$ is determined by the formula 
$[\Psi (e_0), \Psi (f_0)] = (q-q^{-1})^{-1} ( \Psi (k_0^+) - \Psi (k_0^-) )$. 
\end{prop}


\para 
Thanks to Proposition \ref{Prop ev whsln gln} and Proposition \ref{Prop iso quantum loop}, 
we have the algebra homomorphism $ \ev_{\g} \circ \Psi^{-1} : U_q (L \Fsl_n) \ra U_q(\Fgl_n)$, 
and we denote it by $\ev_{\g}$ again. 

Let $P = \bigoplus_{i=1}^n \ZZ \ve_i$ be the weight lattice of $\Fgl_n$, 
and put $\w_i = \ve_1+ \ve_2+ \dots + \ve_i$ for $i \in I$. 
Let $V(\w_i)$ be the simple highest weight $U_q(\Fgl_n)$-module of highest weight $\w_i$, 
and $v_0^{(i)} \in V( \w_i)$ be a highest weight vector. 
Then we have 
\begin{align*}
&E_j \cdot v_0^{(i)} =0 \text{ for all } j \in I, 
\quad 
F_j \cdot v_0^{(i)} =0 \text{ if } j \not=i, 
\quad 
F_i \cdot v_0^{(i)} \not=0, 
\quad 
F_i^2 \cdot v_0^{(i)} =0,
\\
&T_j^{\pm} \cdot v_0^{(i)} 
= \begin{cases} 
	q^{\pm 1} v_0^{(i)} & \text{ if } 1 \leq j \leq i, 
	\\
	v_0^{(i)} & \text{ if } i < j \leq n. 
	\end{cases}
\end{align*}
For each $\g \in \CC^{\times}$, 
we regard $U_q(\Fgl_n)$-module $V(\w_i)$ as  a $U_q (L \Fsl_n)$-module 
through the homomorphism $\ev_{\g} : U_q (L \Fsl_n) \ra U_q(\Fgl_n)$, 
and denote it by $V(\w_i)^{\ev_{\g}}$. 
The following proposition is obtained by the same argument with one in \cite[3.6]{CP94}. 

\begin{prop}[{cf. \cite[3.6]{CP94}}] 
\label{Prop ev loop}
For the $U_q (L \Fsl_n)$-module $V(\w_i)^{\ev_{\g}}$ ($i \in I$, $\g \in \CC^{\times}$), we have 
\begin{align*}
&e_{j,0} \cdot v_0^{(i)} =0 \text{ for all } j \in I, 
\quad 
f_{j,0} \cdot v_0^{(i)} =0 \text{ if } j \not=i, 
\quad 
f_{i,1} \cdot v_0^{(i)} = \g q^{-i+2} f_{i,0} \cdot v_0^{(i)}. 
\end{align*}
\end{prop}


\para 
\label{Def ev} 
Recall the injective homomorphism 
$\Theta^{\lan \mathbf{0} \ran} : U_q^{\lan \mathbf{0} \ran} \ra \cU_{0,0}^{\lan \mathbf{0} \ran} 
\cong U_q(L\Fsl_n)$ in Proposition \ref{Prop ThetaQ}. 
Then we have the algebra homomorphism 
$\ev_{\g} \circ \Theta^{\lan \mathbf{0} \ran} : U_q^{\lan \mathbf{0} \ran} \ra U_q(\Fgl_n)$, 
and we denote it by $\ev_{\g}^{\lan \mathbf{0} \ran}$. 
We cannot define the evaluation homomorphism 
$\ev_{0} : U_q(L \Fsl_n) \ra U_q (\Fgl_n)$ at $\g =0$. 
However, if we restrict $U_q(L\Fsl_n)$ to $U_q^{\lan \mathbf{0} \ran}$, 
we can also define the evaluation homomorphism 
$\ev_{0}^{\lan \mathbf{0} \ran} : U_q^{\lan \mathbf{0} \ran} \ra U_q(\Fgl_n)$ at $\g =0$ by 
\begin{align*}
X_{i,t}^+ \mapsto \d_{t,0} E_i, 
\quad 
X_{i,t}^- \mapsto \d_{t,0} F_i, 
\quad 
J_{i,t} \mapsto \d_{t,0}  \frac{1- (T_i^- T_{i+1}^+)^2}{q-q^{-1}}, 
\quad 
K_i^{+} \mapsto T_i^+ T_{i+1}^-. 
\end{align*} 
For each $\g \in \CC$, 
we regard  the 
$U_q(\Fgl_n)$-module $V(\w_i)$ as  a $U_q^{\lan \mathbf{0} \ran}$-module 
through the homomorphism $\ev_{\g}^{\lan \mathbf{0} \ran} : U_q^{\lan \mathbf{0} \ran} \ra U_q(\Fgl_n)$, 
and denote it by $V(\w_i)^{\ev_{\g}^{\lan \mathbf{0} \ran}}$. 
Then we have the following proposition. 


\begin{prop}
\label{Prop ev module}
For the $U_q^{\lan \mathbf{0} \ran}$-module $V(\w_i)^{\ev_{\g}^{\lan \mathbf{0} \ran}}$ 
($i \in I$, $\g \in \CC$), we have 
\begin{align*}
X_{j,t}^+ \cdot v_0^{(i)} =0, 
\quad 
J_{j,t} \cdot v_0^{(i)} = \begin{cases}
		 q^{-1} (\g q^{- i+2})^t  v_0^{(i)}  & \text{ if } j=i, 
		\\
		0 & \text{ if } j \not=i, 
	\end{cases} 
\quad 
K_j^+ \cdot v_0^{(i)} 
	= \begin{cases} 
		q v_0^{(i)} & \text{ if } j=i 
		\\ 
		v_0^{(i)}  & \text{ if } j \not=i 
		\end{cases} 
\end{align*}
for $j \in I$ and $t \geq 0$.  

\begin{proof}
>From the definitions, 
we have 
\begin{align}
\label{ev Kj+}
K_j^+ \cdot v_0^{(i)} 
	= \begin{cases} 
		q v_0^{(i)} & \text{ if } j=i 
		\\ 
		v_0^{(i)}  & \text{ if } j \not=i 
		\end{cases} 
\end{align} 
for $j \in I$ immediately. 
By Proposition \ref{Prop ev loop} in the case where $\g \not=0$ and direct calculation in the case where $\g =0$, 
we have 
\begin{align}
\label{ev Vwi}
\begin{split}
&X_{j,0}^+ \cdot v_0^{(i)} =0 \text{ for all } j \in I, 
\quad 
X_{j,0}^- \cdot v_0^{(i)} =0 \text{ if } j \not=i, 
\quad 
X_{i,1}^- \cdot v_0^{(i)} = \g q^{-i+2} X_{i,0}^- \cdot v_0^{(i)}. 
\end{split}
\end{align}
By the relation (Q1-1), we see that $J_{j,t}$ acts on $v_0^{(i)}$ as a scalar multiplication 
since the weight space of $V(\w_i)$ with the weight $\w_i$  is one-dimensional. 
Then, by the induction on $t$ using \eqref{ev Vwi} and \eqref{Xit+1pm}, we have 
\begin{align}
\label{ev Xjt+}
X_{j,t}^+ \cdot v_0^{(i)}=0 \text{ for all } j \in I \text{ and } t \geq 0. 
\end{align}
The equations $X_{j,0}^- \cdot v_0^{(i)}=0$ if $j \not=i$ in \eqref{ev Vwi} and \eqref{ev Xjt+} 
together with the relation (Q6) imply that 
$J_{j,t} \cdot v_0^{(i)} =0 $ for all $j \in I \setminus\{i\}$ and $t \geq 0$. 

For $t \geq 0$, 
applying $X_{i,t}^+$ to both sides of the equation 
$X_{i,1}^- \cdot v_0^{(i)} = \g q^{-i+2} X_{i,0}^- \cdot v_0^{(i)}$ in \eqref{ev Vwi}, 
we have $X_{i,t}^+ X_{i,1}^- \cdot v_0^{(i)} = \g q^{-i+2} X_{i,t}^+ X_{i,0}^- \cdot v_0^{(i)}$. 
By the relations (Q1-2), (Q6) and the equation \eqref{ev Xjt+}, 
the above equation implies 
$J_{i,t+1} \cdot v_0^{(i)} = \g q^{-i+2}  J_{i,t} \cdot v_0^{(i)}$. 
Thus we have 
\begin{align*}
J_{i,t} \cdot v_0^{(i)} 
= (\g q^{- i+2})^t J_{i,0} \cdot v_0^{(i)} 
= q^{-1} (\g q^{- i+2})^t  \cdot v_0^{(i)} 
\end{align*}
for $t \geq 0$, 
where the second equation follows from \eqref{ev Kj+},  (Q1-1) and  (Q1-2). 
\end{proof}
\end{prop}


\section{Highest weight $U_q^{\lan \bQ \ran}$-modules} 
\label{section hw}
In the rest of the paper, 
we assume that the parameter $q$ is not a root of unity. 

In this section, we give a notion of highest weight $U_q^{\lan \bQ \ran}$-modules 
with respect to the triangular decomposition 
\eqref{weak tri decom}. 
The argument is standard, so we give only notation and some statements. 


\para 
\textbf{Highest weight modules.} 
For a $U_q^{\lan \bQ \ran}$-module $M$, 
we say that $M$ is a highest weight module if there exists $v_0 \in M$ 
satisfying the following conditions: 
\begin{enumerate}
\item 
$M$ is generated by $v_0$ as a $U_q^{\lan \bQ \ran}$-module. 

\item 
$X_{i,t}^+ \cdot v_0 =0$ for all $(i,t) \in I \times \ZZ_{\geq 0}$. 

\item 
There exists 
$\bu= ((\la_i, (u_{i,t})_{t>0}))_{i \in I} \in ( \CC^{\times} \times \prod_{t>0} \CC)^I$ 
such that 
$K_i^+ \cdot v_0= \la_i v_0$ and $J_{i,t} \cdot v_0 = u_{i,t} v_0$ for each $i \in I$ and $t \in \ZZ_{>0}$. 
\end{enumerate}
In this case, 
we say that $\bu$ is the highest weight of $M$, 
and that $v_0$ is a highest weight vector of $M$. 
We remark that $J_{i,0} \cdot v_0 = (q-q^{-1})^{-1} (1 - \la_i^{-2})$ by the relation (Q1-2). 


\para 
\textbf{Verma modules.} 
For $\bu=((\la_i, (u_{i,t})_{t>0}))_{i \in I} \in ( \CC^{\times} \times \prod_{t>0} \CC)^I$, 
let $\fI(\bu)$ be the left ideal of $U_q^{\lan \bQ \ran}$ generated by 
$X_{i,t}^+$ ($(i,t) \in I \times \ZZ_{\geq 0}$), $K_i^+ - \la_i$ ($i \in I$) and 
$J_{i,t} - u_{i,t}$ ($(i,t) \in I \times \ZZ_{>0}$). 
Then, we define the Verma module as the  quotient module $M(\bu) = U_q^{\lan \bQ \ran} / \fI(\bu)$. 
By the standard argument, 
the Verma module $M(\bu)$ has the unique maximal proper submodule $\rad M(\bu)$, 
and we have the unique simple top $L(\bu) := M(\bu) / \rad M(\bu)$. 
We have the following proposition whose proof is also standard. 


\begin{prop}
\label{Prop fin simple hw}
For $\bu \in ( \CC^{\times} \times \prod_{t>0} \CC)^I$, 
a highest weight simple $U_q^{\lan \bQ \ran}$-module of highest weight $\bu$ is isomorphic to $L(\bu)$. 
Moreover, any finite dimensional simple $U_q^{\lan \bQ \ran}$-module 
is isomorphic to $L(\bu)$ for some $\bu \in ( \CC^{\times} \times \prod_{t>0} \CC)^I$. 
\end{prop}



\section{Some symmetric polynomials} 

In this section, we introduce some symmetric polynomials, and give some properties of them. 
These symmetric polynomials will be used to describe the highest weights of finite dimensional $U_q^{\lan \bQ \ran}$-modules.  

\para 
A partition is a non-increasing sequence $\la=(\la_1,\la_2, \dots)$ of non-negative integers 
with only finitely many non-zero terms. 
The size of a partition $\la$, denoted by $|\la|$, 
is $|\la|= \sum_{i \geq 1} \la_i$. 
We denote by $\la \vdash t$ if $\la$ is a partition of size $t$. 
The length of a partition $\la$ is the number of non-zero terms, and we denote it by $\ell(\la)$. 

Let $\CC[x_1,x_2,\dots, x_k]$ be the polynomial ring over $\CC$ with indeterminate variables $x_1,\dots, x_k$. 
For $t, k \in \ZZ_{>0}$, 
put 
\begin{align*}
& e_t (x_1,\dots, x_k) = \sum_{1 \leq i_1 < i_2 < \dots < i_t \leq k} x_{i_1} x_{i_2} \dots x_{i_t} \in \CC[x_1,\dots,x_k], 
\\
& p_t(x_1,\dots, x_k) = x_1^t + x_2^t + \dots + x_k^t \in \CC [x_1,\dots, x_k]
\end{align*}
and $e_0(x_1,\dots, x_k)=1$. 
Namely, these polynomials are the elementary symmetric polynomial and the power sum symmetric polynomial respectively. 
For a partition $\la =(\la_1,\dots, \la_k) \vdash t$ such that $\ell(\la) \leq k$, 
put 
\begin{align*}
m_{\la} (x_1,\dots, x_k) = \sum_{\mu \in \fS_k \cdot \la} x_1^{\mu_1} x_2^{\mu_2} \dots x_k^{\mu_k} 
\in \CC[x_1,\dots, x_k], 
\end{align*}
where $\fS_k \cdot \la = \{\mu =(\mu_1, \dots, \mu_k) \in \ZZ_{\geq 0}^k 
	\mid \mu_i = \la_{\s(i)} \,  (1\leq i \leq k) \text{ for some } \s \in \fS_k \}$. 
Namely, the polynomial $m_{\la}(x_1,\dots,x_k)$ is the monomial symmetric polynomial associated with $\la$. 

For $t,k \in \ZZ_{>0}$, we define a polynomial 
$p_t(q)(x_1,\dots, x_k) \in \CC[x_1,\dots, x_k]$ by 
\begin{align*}
p_t(q)(x_1,\dots, x_k) :=  \sum_{\la \vdash t \atop \ell(\la) \leq k} 
	q^{-\ell(\la)} (q-q^{-1})^{\ell(\la) -1} m_{\la} (x_1,\dots, x_k). 
\end{align*}
>From the definition, we see that $p_t(q)(x_1,\dots,x_k)$ is a symmetric polynomial. 


\begin{remark}
In the case where $q=1$, we have 
$p_t(1)(x_1,\dots, x_k) = p_t(x_1,\dots, x_k)$. 
Thus,  
the polynomial $p_t(q)(x_1,\dots, x_k)$ is a $q$-analogue of the power sum symmetric polynomial. 
\end{remark}

\begin{lem}
\label{Lemma ptq} 
For $t,k \in \ZZ_{>0}$, the polynomial $p_t(q)(x_1,\dots, x_k)$ satisfies the following equations: 
\begin{enumerate}
\item 
$p_t(q) (x_1,\dots, x_k) $
\\
$\dis =  p_t(q) (x_1,\dots, x_{k-1}) + q^{-1} x_k^t + q^{-1} (q-q^{-1}) \sum_{z=1}^{t-1} p_{z}(q) (x_1,\dots, x_{k-1}) x_k^{t-z}$. 

\item 
$p_t(q) (x_1,\dots, x_k) $
\\
$\dis = (-1)^{t-1} q^{-t} [t] e_t(x_1,\dots, x_k) + \sum_{z=1}^{t-1} (-1)^{t+z-1} p_z (q) (x_1,\dots, x_k) e_{t-z} (x_1,\dots, x_k)$. 

\item 
$\dis p_{k+t} (q) (x_1,\dots, x_k) = \sum_{z=0}^{k-1} (-1)^{k+z-1} p_{t +z} (q) (x_1,\dots, x_k) e_{k-z} (x_1,\dots, x_k)$.
\end{enumerate} 

\begin{proof}
(\roi). 
By the definition of the monomial symmetric polynomials,  
for $1 \leq l \leq k$, 
we see that 
\begin{align*}
&\sum_{\la \vdash t \atop \ell(\la) = l} m_{\la} (x_1,\dots, x_k) 
\\
&= 
\d_{(l \not=k)} \sum_{\la \vdash t \atop \ell(\la) =l} m_{\la} (x_1,\dots, x_{k-1}) 
	+ \d_{(l=1)} x_k^t 
	+ \d_{(l \not=1)} \sum_{z=1}^{t-l+1} \sum_{\la \vdash t-z \atop \ell(\la) = l-1} m_{\la}(x_1,\dots, x_{k-1}) x_k^z, 
\end{align*}
where $\d_{(*)}=1$ if the condition $*$ is true, and $\d_{(*)}=0$ if the condition $*$ is false. 
Thus, we have 
\begin{align*}
&p_t(q) (x_1,\dots, x_k) 
\\
&= \sum_{l=1}^{k-1} q^{-l} (q-q^{-1})^{l-1} \sum_{\la \vdash t \atop \ell(\la)=l} m_{\la}(x_1,\dots, x_{k-1}) 
	+ q^{-1} x_k^t 
	\\ & \quad 
	+ \sum_{l=2}^k q^{-l} (q-q^{-1})^{l-1} \sum_{z=1}^{t-l+1} \sum_{\la \vdash t-z \atop \ell(\la) =l-1} m_{\la} (x_1,\dots, x_{k-1}) x_k^z
\\
&= \sum_{\la \vdash t \atop \ell(\la) \leq k-1} q^{-\ell(\la)} (q-q^{-1})^{\ell(\la) -1}   m_{\la} (x_1,\dots, x_{k-1}) 
	+ q^{-1} x_k^t 
	\\ & \quad 
	+ \sum_{z=1}^{t-1} 
		\Big( \sum_{\la \vdash t-z \atop \ell(\la) \leq k-1} q^{-\ell(\la)-1} (q-q^{-1})^{\ell(\la)} m_{\la} (x_1,\dots, x_{k-1}) \Big) x_k^z
\\
&= p_t(q) (x_1,\dots, x_{k-1}) + q^{-1} x_k^t + q^{-1} (q-q^{-1}) \sum_{z=1}^{t-1}  p_{t-z}(q) (x_1,\dots, x_{k-1}) x_k^z.
\end{align*}

(\roii). 
Put 
\begin{align*}
&\wh{p}_t(q)(x_1,\dots, x_k) 
\\
&= (-1)^{t-1} q^{-t} [t] e_t(x_1,\dots, x_k) + \sum_{z=1}^{t-1} (-1)^{t+z-1} \wh{p}_z (q) (x_1,\dots, x_k) e_{t-z} (x_1,\dots, x_k).
\end{align*} 
We can prove 
the equation (\roi) replacing $p_t(q) (x_1,\dots, x_k)$ with $\wh{p}_t(q) (x_1,\dots, x_k)$
by the induction on the degree $t$. 
Then we can prove that 
$p_t(q)(x_1,\dots,x_k) = \wh{p}_t (q)(x_1,\dots, x_k)$ by the induction on the number $k$ of variables 
using the equation (\roi) for both polynomials. 
As a consequence, we obtain (\roii). 

(\roiii). Note that $e_{t'} (x_1,\dots,x_k) =0$ if $t' >k$, and the equation (\roiii) follows from the equation (\roii). 
\end{proof}
\end{lem}


\begin{cor}
\label{cor alg ind ptq}
For $k \in \ZZ_{>0}$, 
the set of polynomials 
\begin{align*}
\{ p_t(q) (x_1,\dots,x_k) \mid 1 \leq t \leq k\} 
\end{align*} 
is algebraically independent over $\CC$. 

\begin{proof}
Note that we assume that $q$ is not a root of unity. 
Then we can prove the corollary in the same way with the corresponding statement for power sum symmetric polynomials 
using the equation (\roii) in Lemma \ref{Lemma ptq}. 
\end{proof}
\end{cor}

\begin{cor}
\label{Cor et pla}
For $t=1,2,\dots,k$, 
there exist the unique $a_\la \in \CC$ ($\la \vdash t$) such that 
\begin{align*}
e_t (x_1,x_2,\dots, x_k) 
= \sum_{\la \vdash t} a_{\la} p_{\la}(q) (x_1,x_2, \dots, x_k), 
\end{align*}
where we put 
$p_{\la} (q)(x_1,\dots,x_k)= \prod_{i=1}^{\ell (\la)} p_{\la_i}(q) (x_1,\dots, x_k)$. 

\begin{proof}
We can prove the existence of the numbers $a_{\la}$ ($\la \vdash t$) by the induction on $t$ using Lemma \ref{Lemma ptq} (\roii). 
The uniqueness of $a_{\la}$ ($\la \vdash t$) follows from Corollary \ref{cor alg ind ptq}. 
\end{proof}
\end{cor}

\begin{prop}
\label{Prop gen series pt(q)}
For $k \in \ZZ_{>0}$, 
let $P^{(k)} (\w) = 1 + (q-q^{-1}) \sum_{t >0} p_t(q) (x_1,\dots, x_k) \w^t $ be 
the generating function. 
Then we have 
\begin{align*}
P^{(k)} (\w) 
=\frac{(1 - q^{-2} x_1 \w) (1- q^{-2} x_2 \w) \dots (1 - q^{-2} x_k \w)}{( 1 - x_1 \w) (1- x_2 \w) \dots (1- x_k \w)}. 
\end{align*}

\begin{proof}
In this proof, 
we denote $p_t(q) (x_1,\dots, x_k)$ (resp. $e_t(x_1,\dots, x_k)$) by 
$p_t(q)$ (resp. $e_t$) simply. 
We consider the generating function  
$E^{(k)}(\w) = \sum_{t \geq 0} (-1)^t e_t  \w^t$. 
Then, we have 
\begin{align*}
P^{(k)}(\w) E^{(k)}(\w) 
&= \Big( 1 + (q-q^{-1}) \sum_{t>0} p_t(q) \w^t \Big) \big( \sum_{t \geq 0} (-1)^t e_t \w^t \big) 
\\
&= \sum_{t \geq 0} (-1)^t e_t \w^t 
	+ \sum_{ t >0} \Big( (q-q^{-1}) \sum_{z=1}^t (-1)^{t-z}p_z(q) e_{t-z} \Big) \w^t. 
\end{align*} 
Applying Lemma \ref{Lemma ptq} (\roii), we have 
\begin{align*}
P^{(k)}(\w) E^{(k)}(\w) 
&= \sum_{t \geq 0} (-1)^t e_t \w^t 
	+\sum_{t>0}   (q-q^{-1}) (-1)^{t-1} q^{-t}[t] e_t \w^t
\\
&= \sum_{t \geq 0}(-1)^t  q^{-2 t}  e_t \w^t
\\
&= E^{(k)} (q^{-2} \w).
\end{align*}
On the other hand, 
we have 
\begin{align*}
E^{(k)} (\w) 
= \sum_{t=0}^k (-1)^t e_t (x_1, \dots, x_k) \w^t 
= ( 1 - x_1 \w) (1- x_2 \w) \dots (1- x_k \w)
\end{align*}
since $e_t(x_1,\dots, x_k)=0$ if $t>k$. 
As a consequence, we have 
\begin{align*}
P^{(k)} (\w) 
= \frac{E^{(k)} (q^{-2} \w)}{E^{(k)} (\w)} 
=\frac{(1 - q^{-2} x_1 \w) (1- q^{-2} x_2 \w) \dots (1 - q^{-2} x_k \w)}{( 1 - x_1 \w) (1- x_2 \w) \dots (1- x_k \w)}. 
\end{align*}
\end{proof}
\end{prop}

\remark 
The formula in Lemma \ref{Lemma ptq} (\roii) corresponds to the definition \eqref{def J0[t]} 
which is identified with \cite[Proposition 3.5 $(\mathrm{\roii})_r$]{CP91} under the injective algebra homomorphism $\Theta^{\lan 0 \ran}$ (see the paragraph \ref{Par J0t}). 
Thus, the formula in Lemma \ref{Lemma ptq} (\roii) is a $q$-analogue of Newton's formula relating the elementary symmetric polynomials and the power sums suggested in \cite[Remark 3.5]{CP91}. 
Under this correspondence, 
Proposition \ref{Prop gen series pt(q)} corresponds to \cite[Corollary 3.5]{CP91} 
(see also Corollary \ref{Cor Psi+ Q=0}). 

\para 
For $t,k \in \ZZ_{>0}$ and $Q, \b \in \CC^{\times}$, 
we define a polynomial $p_t^{\lan Q \ran} (q; \b) (x_1,\dots, x_k) \in \CC[x_1,\dots,x_k] $ by 
\begin{align}
\label{def ptQqb} 
\begin{split}
&p_t^{\lan Q \ran} (q;\b) (x_1,\dots, x_k)  
\\
&:=p_t(q) (x_1,\dots,x_k) + \wt{\b} Q^{-t} + (q-q^{-1}) \sum_{z=1}^{t-1} \wt{\b} Q^{-t+z} p_z (q) (x_1,\dots,x_k), 
\end{split} 
\end{align}
where we put $\wt{\b} = (q-q^{-1})^{-1} (1 - \b^{-2})$. 

By definition, the polynomial $p_t^{\lan Q \ran} (q;\b) (x_1,\dots,x_k)$ is a symmetric polynomial. 
In the case where $\b=\pm1$, 
we have $p_t^{\lan Q \ran} (q; \pm 1) (x_1,\dots,x_k)= p_t(q) (x_1,\dots,x_k)$. 

By the definition \eqref{def ptQqb} together with Corollary \ref{cor alg ind ptq}, we have the following lemma. 

\begin{lem}
\label{Lemma ptQ alg ind}
For $k \in \ZZ_{>0}$ and $Q, \b \in \CC^{\times}$, 
the set of polynomials 
\begin{align*}
\{ p_t^{\lan Q \ran} (q;\b)(x_1,\dots,x_k) \mid 1\leq t \leq k\} 
\end{align*} 
is algebraically independent over $\CC$. 
\end{lem}


\begin{lem} 
\label{Lemma ptQqb} 
Fix $k \in \ZZ_{>0}$, and 
we put $\bx =(x_1,\dots, x_k)$ for simplicity. 
For $t \in \ZZ_{> 0}$ and $Q, \b \in \CC^{\times}$, 
the polynomial $p_t^{\lan Q \ran} (q;\b) (\bx)$ satisfies the following equations: 
\begin{enumerate}
\item 
$p_t^{\lan Q \ran}(q;\b) (\bx) $
\\
$\dis = (-1)^{t-1} q^{-t} [t] e_t (\bx) + \wt{\b} Q^{-t}
+ \sum_{z=1}^{t-1} (-1)^{t-z+1} \big( p_{z}^{\lan Q \ran}(q;\b) (\bx) 
	- q^{-2 (t-z)} \wt{\b} Q^{-z}\big) e_{t-z} (\bx)$. 
	
\item 
$p_{k+t}^{\lan Q \ran} (q;\b) (\bx)$ 
\\
$\dis = Q^{-1} p_{k+t-1}^{\lan Q \ran} (q;\b) (\bx) + \sum_{z=0}^{k-1} (-1)^{k-z+1} 
	\big( p_{t+z}^{\lan Q \ran} (q;\b) (\bx) - Q^{-1} p_{t+z-1}^{\lan Q \ran}(q;\b)(\bx) \big) 
	e_{k-z} (\bx)$, 
where we put $\dis p_0^{\lan Q \ran} (q;\b)(\bx) = \frac{1 - (\b q^k)^{-2}}{q-q^{-1}}$. 
Note that the scalar $p_0^{\lan Q \ran} (q;\b)(\bx)$ appears only the case where $t=1$. 
\end{enumerate}

\begin{proof}
The equation (\roi) follows from the definition \eqref{def ptQqb} and Lemma \ref{Lemma ptq} (\roii). 
We prove (\roii). 
Note that $e_{t'}(\bx)=0$ if $t'>k$, 
then  the equation (\roi) implies 
\begin{align*}
&p_{k+t}^{\lan Q \ran}(q;\b)(\bx) - Q^{-1} p_{k+t-1}^{\lan Q \ran} (q;\b)(\bx) 
\\
&= \wt{\b} Q^{-(k+t)} + \sum_{z=t}^{k+t-1} (-1)^{k+t-z+1} 
	\big( p_z^{\lan Q \ran} (q;\b) (\bx) - q^{-2 (k+t-z)} \wt{\b} Q^{-z} \big) e_{k+t-z} (\bx)
	\\ & \quad 
	- Q^{-1} \big\{ \d_{t,1} (-1)^{k-1} q^{-k}[k] e_k(\bx) +   \wt{\b} Q^{-(k+t-1)} 
	\\ & \hspace{5em}
	+ \sum_{z= \max\{1,t-1\}}^{k+t-2} (-1)^{k+t-z} 
	\big( p_z^{\lan Q \ran} (q;\b) (\bx) - q^{-2 (k+t-z-1)} \wt{\b} Q^{-z} \big) e_{k+t-z-1} (\bx) \big\}
\\
&= \begin{cases}
	\dis \sum_{z=t}^{k+t-1} (-1)^{k+t-z+1} \big(  p_z^{\lan Q \ran} (q;\b) (\bx) - Q^{-1} p_{z-1}^{\lan Q \ran} (q;\b) (\bx) \big) 
		e_{k+t-z} (\bx) & \text{ if } t >1, 
	\\[1em]
	\dis (-1)^{k+1} \big( p_1^{\lan Q \ran} (q;\b) (\bx) - q^{-2 k}\wt{\b} Q^{-1}  - Q^{-1} q^{-k} [k] \big) e_{k} (\bx) 
		\\ \quad \dis 
		+ \sum_{z=2}^{k} (-1)^{k-z+2} \big( p_z^{\lan Q \ran} (q;\b) (\bx) - Q^{-1} p_{z-1}^{\lan \bQ \ran} (q;\b)(\bx) \big) 
			e_{k +1 -z} (\bx) 
		& \text{ if } t=1. 
	\end{cases}
\end{align*} 
Note that 
$q^{-2k} \wt{\b} + q^{-k}[k] = (q-q^{-1})^{-1} ( 1- (\b q^k)^{-2}) = p_0^{\lan Q \ran} (q;\b) (\bx)$, 
we have the equation (\roii) by replacing $z-t$ with $z$. 
\end{proof}
\end{lem}


\section{One-dimensional $U_q^{\lan \bQ \ran}$-modules}

In this section, we classify  one-dimensional 
$U_q^{\lan \bQ \ran}$-modules. 

\para 
Let $L  = \CC v$ be a  one-dimensional 
$U_q^{\lan \bQ \ran}$-module with a basis $ v$. 
Then $K_i^+$ ($i \in I$) acts on $v$ as a scalar multiplication. 
We denote the eigenvalue of the action of $K_i^+$  by $\b_i$. 
By the relation (Q1-2), we have $\b_i \not=0$. 

For $(j,t) \in I \times \ZZ_{\geq 0}$, 
the element $X_{j,t}^{\pm} \cdot v$ is an eigenvector of the eigenvalue $q^{\pm a_{ij}} \b_i$ 
for the action of $K_i^+$ ($i \in I$) if $X_{j,t}^{\pm} \cdot v \not=0$ 
by the relations (Q4-1) and (Q5-1). 
However, $L$ is one-dimensional, 
thus we have 
$X_{j,t}^{\pm} \cdot v=0$ for all $(j,t) \in I \times \ZZ_{\geq 0}$. 

For $(i,t) \in I \times \ZZ_{\geq 0}$, 
we have 
$(K_i^+ J_{i,t} - Q_i K_i^+ J_{i,t+1}) \cdot v = [X_{i,t}^+, X_{i,0}^-] \cdot v =0$ 
by the relation (Q6). 
This equation implies that 
$J_{i,t} \cdot v=0$ if $Q_i=0$, and $J_{i,t} \cdot v = Q_i J_{i,t+1} \cdot v$ if $Q_i \not=0$. 
Thus, we have $J_{i,t} \cdot v = Q_i^{-t} \cdot J_{i,0} \cdot v$ if $Q_i \not=0$. 
On the other hand, by the relation (Q1-2), 
we have $J_{i,0} \cdot v= (q-q^{-1})^{-1} (1 - \b_i^{-2}) v$. 
Then we have $\b_i = \pm 1$ if $Q_i=0$ since $J_{i,0} \cdot v =0$ in this case. 
As a consequence, we have 
\begin{align}
\label{act one dim rep}
X_{i,t}^{\pm} \cdot v =0, 
\quad 
K_i^{\pm} \cdot v = \b_i^{\pm 1} \cdot v, 
\quad 
J_{i,t} \cdot v = 
	\begin{cases} 
		0 & \text{ if } Q_i=0, 
		\\ 
		\dis \frac{1 - \b_i^{-2}}{q-q^{-1}} Q_i^{-t} v & \text{ if } Q_i \not=0 
	\end{cases}
\end{align}
for $i \in I$ and $t \in \ZZ_{\geq 0}$, where $\b_i = \pm 1$ if $Q_i=0$. 

\para 
For $Q \in \CC$, 
put $\BB^{\lan Q \ran} = \{\pm 1\}$ if $Q=0$ and $\BB^{\lan Q \ran} = \CC^{\times}$ if $Q \not=0$. 
Then we put $\BB^{\lan \bQ \ran} = \prod_{i \in I} \BB^{\lan Q_i \ran}$ 
for $\bQ=(Q_1,\dots, Q_{n-1}) \in \CC^{n-1}$. 
For $\Bb =(\b_i)_{i \in I} \in\BB^{\lan \bQ \ran}$, 
we define a one-dimensional 
$U_q^{\lan \bQ \ran}$-module $\cD_{\Bb}^{\lan \bQ \ran} = \CC v$ 
by \eqref{act one dim rep}. 
We can easily check that this action is well-defined. 
As a consequence of this section, we have the following proposition. 

\begin{prop}
\label{Prop one dim rep}
Any  one-dimensional 
$U_q^{\lan \bQ \ran}$-module is isomorphic to $\cD_{\Bb}^{\lan \bQ \ran}$ 
for some $\Bb \in \BB^{\lan \bQ \ran}$. 
\end{prop}

By \eqref{Def Psi Q=0}, \eqref{Def Psi Q not=0} and \eqref{act one dim rep}, 
we have the following corollary. 
\begin{cor}
For $\cD_{\Bb}^{\lan \bQ \ran} = \CC v $ ($\Bb \in \BB^{\lan \bQ \ran}$), 
we have 
\begin{align*}
\Psi_i^+ (\w) \cdot v = \begin{cases} 
		\b_i v &  \text{ if } Q_i=0, 
		\\
		(\b_i^{-1} - Q_i \b_i \w^{-1} ) v  & \text{ if } Q_i \not=0. 
	\end{cases}
\end{align*}
\end{cor}



\section{Finite dimensional simple modules of $U_q (\Fsl_2^{\lan 0 \ran}[x])$}

In this section, we classify the isomorphism classes of finite dimensional simple modules 
of the algebra $U_q^{\lan 0 \ran} = U_q (\Fsl_2^{\lan 0 \ran}[x])$ 
in the case of rank one and of $Q=0$. 

We recall that, in the case where  $Q=0$, 
the algebra $U_q^{\lan 0 \ran}$ is a subalgebra of 
the quantum loop algebra $U_q (L\Fsl_2)$ through the injective homomorphism 
$\Theta^{\lan 0 \ran}$ in Proposition \ref{Prop ThetaQ}. 
In this case, the argument to classify the finite dimensional simple
$U_q^{\lan 0 \ran}$-modules is essentially the same as  the argument for 
$U_q(L \Fsl_2)$ given in \cite{CP91}. 
However, we discuss the case where $Q=0$ in this section for completeness, 
and it is also useful in order to consider the case where $Q \not=0$ in the next section. 

In this and next sections, we consider only the case of rank one, namely $I=\{1\}$, 
so we omit the indices for $I$, 
e.g. we denote $X_{1,t}^{\pm}$ by $X_t^{\pm}$ simply, and so on. 


\para 
\label{Par J0t}
For $t,k \in \ZZ_{\geq 0}$, put 
\begin{align*}
X_t^{+(k)} = \frac{(X_t^+)^k}{[k]!}, 
\quad 
X_t^{-(k)} = \frac{(X_t^-)^k}{[k]!}. 
\end{align*}

For $t \in \ZZ_{\geq 0}$, 
we define the element $J_{[t]}^{\lan 0 \ran} \in U_q^{\lan 0 \ran}$ inductively by 
\begin{align}
\label{def J0[t]}
J_{[0]}^{\lan 0 \ran} = 1 
\text{ and }
J_{[t]}^{\lan 0 \ran} = q^t \frac{1}{[t]} \sum_{z=1}^t (-1)^{z-1}  J_z J_{[t-z]}^{\lan 0 \ran}
\text{ for } t >0.
\end{align}
For examples, we have 
\begin{align*}
&J_{[0]}^{\lan 0 \ran} =1, 
\quad 
J_{[1]}^{\lan 0 \ran} = q J_1, 
\quad 
J_{[2]}^{\lan 0 \ran} = \frac{1}{[2]} ( q^3 J_1^2 - q^2J_2), 
\\
&J_{[3]}^{\lan 0 \ran} =  \frac{1}{[3]!} \big(q^6 J_1^3 - (2 q^5 + q^3) J_1J_2 + (q^4 + q^2) J_3 \big).
\end{align*}

Compare the definition \eqref{def J0[t]} with \cite[Proposition 3.5 $(\textrm{\roii})_r$]{CP91} under 
the injective algebra homomorphism $\Theta^{\lan 0 \ran} : U_q^{\lan 0 \ran} \ra \cU^{\lan 0 \ran}_{0,0} \cong U_q (L \Fsl_2)$, 
then we see that 
$\Theta^{\lan 0 \ran} (J^{\lan 0 \ran}_{[t]}) = (-1)^t P_t$ 
for $t \in \ZZ_{\geq 0}$, 
where $P_t \in U_q(L \Fsl_2)$ is an element given in \cite[Proposition 3.5]{CP91}. 
The following lemma is a slight variation of \cite[Proposition 3.5 $(\textrm{\roiii})_r$]{CP91}. 

\begin{lem}
\label{Lemma X1+k X0-k+1}
For $k \in \ZZ_{>0}$, we have 
\begin{align*}
X_1^{+(k)} X_0^{-(k+1)} 
\equiv 
q^{-k(k+1)} \sum_{z=0}^k (-1)^z X_z^- (K^+)^k J_{[k-z]}^{\lan 0 \ran} 
\mod \fX_+,  
\end{align*} 
where $\fX_+$ is the left ideal of $U_q^{\lan 0 \ran}$ generated by $\{X_t^+ \mid t \geq 0\}$. 

\begin{proof}
See Appendix \ref{proof rel X1+k X0-k+1}. 
\end{proof}
\end{lem}

\para 
By Proposition \ref{Prop fin simple hw}, 
any finite dimensional simple $U_q^{\lan 0 \ran}$-module 
is isomorphic to the highest weight module $L( \bu)$ 
for some $\bu= (\la, (u_t)_{t>0}) \in \CC^{\times} \times \prod_{t>0} \CC$. 
We have the following necessary condition for $L(\bu)$ to be finite dimensional. 

\begin{prop}
\label{Prop fin simple sl20 1}
For $\bu= (\la, (u_t)_{t>0}) \in \CC^{\times} \times \prod_{t>0} \CC$, 
if the highest weight simple $U_q^{\lan 0 \ran}$-module $L(\bu)$ is finite dimensional, 
then there  exist 
$k \in \ZZ_{\geq 0}$ and $\g_1, \g_2, \dots, \g_k \in \CC$ such that 
\begin{align}
\label{weight sl20}
\begin{split}
\la = \pm q^k, 
\quad 
u_t = \begin{cases} 
	0 & \text{ if } k =0, 
	\\
	p_t(q) (\g_1, \g_2,\dots, \g_k) & \text{ if } k >0 
	\end{cases}
	\quad (t >0).
\end{split}
\end{align}
\begin{proof}
Let $v_0 \in L(\bu)$ be a highest weight vector. 
By the relation (Q5-1), 
we have 
$K^+ X_0^{-(k)} \cdot v_0 = q^{- 2k} \la X_0^{-(k)} \cdot v_0$. 
Namely $X_0^{-(k)} \cdot v_0$ is an eigenvector of the eigenvalue $q^{-2k}\la$ for the action of $K^+$  
if $X_0^{-(k)} \cdot v_0 \not=0$. 
Thus, there exists a non-negative integer $k$ such that 
$X_0^{-(k)} \cdot v_0 \not=0$ and $X_0^{-(k+1)} \cdot v_0 =0$ 
since $L(\bu)$ is finite dimensional. 

In the case where $k=0$, 
we can easily check that $L(\bu) $ is one-dimensional. 
Then we have \eqref{weight sl20} by Proposition \ref{Prop one dim rep}. 

Assume that $k >0$. 
By the induction on $c$ using the relation (Q1-2) and (Q6), 
we can show that $[X_0^+, X_0^{-(c)}] =  X_0^{-(c-1)} (q-q^{-1})^{-1} (q^{-c+1} K^+ - q^{c-1} K^-)$ 
for $c>0$.  
Then we have 
\begin{align*}
0 = X_0^+ X_0^{-(k+1)} \cdot v_0 
= \frac{ q^{-k } \la - q^{k} \la^{-1} }{q-q^{-1}} X_0^{-(k)} \cdot v_0.
\end{align*}
This implies that $\la = \pm q^k$ since $X_0^{-(k)} \cdot v_0 \not=0$. 

By Corollary \ref{cor alg ind ptq}, 
there  exist 
$\g_1, \g_2,\dots, \g_k \in \CC$ such that 
$u_t = p_t (q) (\g_1, \g_2,\dots,\g_k) $ for $t=1,2,\dots, k$. 

By the induction on $t$ using \eqref{def J0[t]} and Lemma \ref{Lemma ptq} (\roii), 
we see that 
\begin{align}
\label{J0[t]} 
J^{\lan 0 \ran}_{[t]} \cdot v_0 = e_t (\g_1,\g_2,\dots, \g_k) v_0
\end{align} 
for $t=1,2,\dots, k$. 

By Lemma \ref{Lemma X1+k X0-k+1} and the relation (Q6), 
for $t >0$, we have 
\begin{align*}
0 = X_t^+ X_1^{+(k)} X_0^{-(k+1)} \cdot v_0 
= q^{-k(k+1)} \la^{k+1} \sum_{z=0}^k (-1)^z J_{t+z} J_{[k-z]}^{\lan 0 \ran} \cdot v_0, 
\end{align*}
where we note that $X_t^+ (K^+)^k J_{[k-z]}^{\lan 0 \ran} \cdot v_0 =0$ since $v_0$ is a highest weight vector. 
Note that $J_{[0]}^{\lan 0 \ran}=1$, 
this equation implies that 
\begin{align}
\label{Jt+k v0}
 J_{t+k} \cdot v_0 
= \sum_{z=0}^{k-1} (-1)^{k-z+1} J_{t+z} J^{\lan 0 \ran}_{[k-z]} \cdot v_0. 
\end{align}
Then we can show that $u_{t+k} = p_{t+k} (q) (\g_1,\dots,\g_k)$ for $t >0$ by the induction on $t$ 
using \eqref{J0[t]}, \eqref{Jt+k v0} and Lemma \ref{Lemma ptq} (\roiii). 
\end{proof}
\end{prop}

\para 
In order to prove that 
the highest weight simple module $L(\bu)$ is finite dimensional 
if $\bu$ is given by \eqref{weight sl20}, 
we use evaluation modules through the following evaluation homomorphisms 
from $U_q^{\lan 0 \ran}$ to the quantum group $U_q(\Fsl_2)$. 
Let $e,f$ and $K^{\pm}$ be the usual Chevalley generators of $U_q(\Fsl_2)$. 
For $\g \in \CC$, 
we have the algebra homomorphism 
$\wt{\ev}_\g^{\lan 0 \ran} : U_q^{\lan 0 \ran} \ra U_q(\Fsl_2)$ such that 
\begin{align*}
&X_t^+ \mapsto \g^t q^{-t} (K^+)^t e, 
\quad 
X_t^- \mapsto \g^t q^{-t} f (K^+)^t, 
\quad 
K^{\pm} \mapsto K^{\pm}, 
\\
&J_t \mapsto \g^t q^{-t} (K^+)^t \frac{1 - (K^-)^2}{q-q^{-1}} - \g^t (q^t - q^{-t}) (K^+)^{t-1} fe. 
\end{align*}
We remark that, if $\g \not=0$, 
the homomorphism $\wt{\ev}_\g^{\lan 0 \ran}$ 
is the restriction of the evaluation homomorphism $\ev_{\g} : U_q(L \Fsl_2) \ra U_q (\Fsl_2)$ 
given in \cite[Proposition 4.1]{CP91} through the injection $\Theta^{\lan 0 \ran}: U_q^{\lan 0 \ran} \ra U_q(L\Fsl_2)$. 
In the case where $\g =0$, we can easily check the well-definedness of $\wt{\ev}_\g^{\lan 0 \ran}$ by direct calculations. 

Let $V_1= \CC v_0  \op \CC v_1$ be the two-dimensional simple $U_q(\Fsl_2)$-module of type $1$, 
namely the action of $U_q(\Fsl_2)$ is given by 
\begin{align*}
K^+ \cdot v_0 = q v_0, \quad e \cdot v_0 =0, \quad f \cdot v_0 = v_1, 
\quad 
K^+ \cdot v_1 = q^{-1} v_1, \quad e \cdot v_1 = v_0, \quad f \cdot v_1 =0. 
\end{align*}
For $\g \in \CC$, we regard $V_1$ as a $U_q^{\lan 0 \ran}$-module through the homomorphism 
$\wt{\ev}_{\g}^{\lan 0 \ran}$, 
and we denote it by $V_1^{\wt{\ev}_\g^{\lan 0 \ran}}$. 
By definition, we have 
\begin{align}
\label{act V1wtevg}
X_t^+ \cdot v_0 =0, 
\quad 
K^+ \cdot v_0 = q v_0, 
\quad 
J_t \cdot v_0 = q^{-1} \g^t v_0 
\quad (t \geq 0). 
\end{align}

\begin{remark} 
We can also discuss by using the evaluation homomorphisms 
$\ev_\g^{\lan 0 \ran} : U_q^{\lan 0 \ran} \ra U_q(\Fgl_2)$ 
given in \S \ref{section ev}. 
Both arguments are essentially the same although the eigenvalues for the action of $J_t$ are different. 
In this section, we use $\wt{\ev}_\g^{\lan 0 \ran}$  instead of $\ev_\g^{\lan 0 \ran}$ 
for a compatibility with the argument in \cite{CP91} 
(see Remark \ref{Remark Drinfeld poly}). 
\end{remark}

\begin{prop}
\label{Prop fin simple sl20 2} 
For $\bu= (\la, (u_t)_{t>0}) \in \CC^{\times} \times \prod_{t>0} \CC$, 
if there  exist  
$k \in \ZZ_{\geq 0}$ and $\g_1,\g_2,\dots, \g_k \in \CC$ 
such that 
\begin{align}
\label{weight sl20 2}
\begin{split}
\la = \pm q^k, 
\quad 
u_t = \begin{cases} 
	0 & \text{ if } k =0, 
	\\
	p_t(q) (\g_1, \g_2,\dots, \g_k) & \text{ if } k >0 
	\end{cases}
	\quad (t >0), 
\end{split}
\end{align}
then $L(\bu)$ is finite dimensional. 

\begin{proof}
Note that the coproduct $\D^{\lan 0 \ran}$ on $U_q^{\lan 0 \ran}$ 
is a restriction of the coproduct on $U_q(L\Fsl_2)$ through the injection $\Theta^{\lan 0 \ran}$. 
Then, by \cite[Proposition 4.4]{CP91}, 
we have 
\begin{align}
\label{D0 mod}
\begin{split} 
&\D^{\lan 0 \ran} (X_{t'}^+) 
	\equiv  X_{t'}^+ \otimes K^+ + 1 \otimes X_{t'}^+ + (q-q^{-1}) \sum_{z=1}^{t'} X_{t'-z}^+ \otimes K^+ J_z  
	\mod \fX_+^2 \otimes \fX_-, 
\\
& \D^{\lan 0 \ran} (J_t) 
	\equiv J_t \otimes 1 + 1 \otimes J_t + (q-q^{-1}) \sum_{z=1}^{t-1} J_z \otimes J_{t-z} 
	\mod \fX_+ \otimes \fX_-,  
\end{split}
\end{align}
for $t' \geq 0$ and $t >0$, 
where $\fX_+^2$ (resp. $\fX_+$, $\fX_-$) is the left ideal of $U_q^{\lan 0 \ran}$ 
generated by $\{X_s^+ X_{s'}^+ \mid s,s' \geq 0\}$ (resp. $\{X_s^+ \mid s \geq 0\}$, $\{X_s^- \mid s \geq 0 \}$). 

For each $\g_i$ ($1\leq i \leq k$), 
we consider the evaluation module $V_1^{\wt{\ev}_{\g_i}^{\lan 0 \ran}}$ at $\g_i$, and 
let $v_0^{(i)} \in V_1^{\wt{\ev}_{\g_i}^{\lan 0 \ran}}$ be a highest weight vector. 
We also consider the one-dimensional $U_q^{\lan 0 \ran}$-module $\cD_{\pm 1}^{\lan 0 \ran} = \CC v$ given by 
\eqref{act one dim rep}. 
Through the coproduct $\D^{\lan 0 \ran}$, 
we consider the $U_q^{\lan 0 \ran}$-module 
$\cD_{\pm 1}^{\lan 0 \ran} \otimes V_1^{\wt{\ev}_{\g_1}^{\lan 0 \ran}} \otimes V_1^{\wt{\ev}_{\g_2}^{\lan 0 \ran}} 
\otimes \dots \otimes V_1^{\wt{\ev}_{\g_k}^{\lan 0 \ran}}$. 
Let $V(\pm 1; \g_1,\dots, \g_k)$ be the $U_q^{\lan 0 \ran}$-submodule of 
$\cD_{\pm 1}^{\lan 0 \ran} \otimes V_1^{\wt{\ev}_{\g_1}^{\lan 0 \ran}} \otimes V_1^{\wt{\ev}_{\g_2}^{\lan 0 \ran}} 
\otimes \dots \otimes V_1^{\wt{\ev}_{\g_k}^{\lan 0 \ran}}$ 
generated by $v \otimes v_0^{(1)} \otimes v_0^{(2)} \otimes \dots \otimes v_0^{(k)}$. 
Then, by definition together with \eqref{D0 mod}, we see that 
\begin{align*}
&X_t^+ \cdot (v \otimes v_0^{(1)} \otimes \dots \otimes v_0^{(k)})=0 \quad (t \geq 0), 
\\
& K^+ \cdot (v \otimes v_0^{(1)} \otimes \dots \otimes v_0^{(k)}) = \pm q^k v \otimes v_0^{(1)} \otimes \dots \otimes v_0^{(k)}.
\end{align*}
For $t >0$, 
we have 
$J_t \cdot (v \otimes v_0^{(1)}) = q^{-1} \g_1^t v \otimes v_0^{(1)}= p_t(q) (\g_1) v \otimes v_0^{(1)}$ 
in $\cD_{\pm 1}^{\lan 0 \ran} \otimes V_1^{\wt{\ev}_{\g_1}^{\lan 0 \ran}}$
by \eqref{act one dim rep}, \eqref{act V1wtevg} and \eqref{D0 mod}. 
By the induction on $k$ using \eqref{act one dim rep}, \eqref{act V1wtevg}, \eqref{D0 mod} and Lemma \ref{Lemma ptq} (\roi), 
we can show that 
\begin{align*}
J_t \cdot (v \otimes v_0^{(1)} \otimes \dots \otimes v_0^{(k)}) 
= p_t(q) (\g_1,\g_2,\dots, \g_k)  v \otimes v_0^{(1)} \otimes \dots \otimes v_0^{(k)}. 
\end{align*}
As a consequence, 
the $U_q^{\lan 0 \ran}$-module 
$V(\pm1; \g_1, \dots, \g_k)$ is a highest weight module of the highest weight $\bu$ given by \eqref{weight sl20 2}, 
and $L(\bu)$ is a quotient of $V(\pm 1; \g_1,\dots, \g_n)$. 
Thus $L(\bu)$ is finite dimensional since $V(\pm1; \g_1, \dots, \g_k)$ is finite dimensional. 
\end{proof}
\end{prop}

\para 
Let $\CC[x]$ be the polynomial ring over $\CC$ with an indeterminate variable $x$. 
For $\vf \in \CC[x]$, we denote the leading coefficient of $\vf$ by $\b_\vf$. 
Put 
\begin{align*}
\CC[x]^{\lan 0 \ran} = \{ \vf \in \CC[x] \setminus \{0\} \mid \b_{\vf} = \pm 1\}.
\end{align*}
We define a map 
$\bu^{\lan 0 \ran} : \CC[x]^{\lan 0 \ran} \ra \CC^{\times} \times \prod_{t>0} \CC$ 
by 
\begin{align}
\label{def bu0}
\bu^{\lan 0 \ran} (\vf) = 
\begin{cases}
(\b_{\vf}, (0)_{t>0}) & \text{ if } \deg \vf =0,
\\
(\b_{\vf} q^{\deg \vf}, (p_t(q)(\g_1,\g_2,\dots, \g_k) )_{t>0}) & \text{ if } \deg \vf >0 
\end{cases}
\end{align}
for $\vf = \b_{\vf} (x-\g_1) (x-\g_2)\dots (x-\g_k) \in \CC[x]^{\lan 0 \ran}$. 

\begin{lem}
\label{Lemma inj u0} 
The map $\bu^{\lan 0 \ran} : \CC[x]^{\lan 0 \ran} \ra \CC^{\times} \times \prod_{t>0} \CC$  is injective. 

\begin{proof}
For $\vf, \vf'\in \CC[x]^{\lan 0 \ran}$, 
write 
$\vf=\b_{\vf} (x-\g_1)(x-\g_2) \dots (x-\g_k)$ and $\vf'= \b_{\vf'} (x-\g'_1) (x-\g'_2) \dots (x-\g'_l)$. 
If $\bu^{\lan 0 \ran} (\vf) = \bu^{\lan 0 \ran}(\vf')$, 
then we have 
\begin{align}
\label{ptq g} 
\b_{\vf} q^{k} = \b_{\vf'} q^{l}, 
\quad 
p_t(q) (\g_1,\g_2,\dots, \g_k) 
= 
p_t (q) (\g'_1,\g'_2, \dots, \g'_l) 
\quad (t >0). 
\end{align} 
The first equation implies that $\b_{\vf} = \b_{\vf'}$ and $k=l$ since $q$ is not a root of unity. 
Moreover,  
we have 
\begin{align*}
\vf 
&= \b_{\vf}x^k +  \b_{\vf} \sum_{z=1}^k (-1)^z e_z(\g_1,\g_2,\dots, \g_k) x^{k-z}
\\
&=\b_{\vf}x^k +  \b_{\vf} \sum_{z=1}^k (-1)^z \big( \sum_{\la \vdash z} a_{\la} p_{\la} (q)(\g_1, \g_2, \dots, \g_k) \big) x^{k-z}
\\
&=\b_{\vf'} x^k +  \b_{\vf'} \sum_{z=1}^k (-1)^z \big( \sum_{\la \vdash z} a_{\la} p_{\la} (q)(\g'_1,  \g'_2, \dots, \g'_k) \big) x^{k-z}
\\
&= \b_{\vf'}x^k +  \b_{\vf'} \sum_{z=1}^k (-1)^z e_z(\g'_1,\g'_2,\dots, \g'_k) x^{k-z}
\\
&= \vf'
\end{align*}
by Corollary \ref{Cor et pla} and \eqref{ptq g}. 
\end{proof}
\end{lem}


Proposition \ref{Prop fin simple sl20 1}, Proposition \ref{Prop fin simple sl20 2} and Lemma \ref{Lemma inj u0} imply 
the following theorem. 

\begin{thm}
\label{Thm simple sl20}
There exists the bijection between 
$\CC[x]^{\lan 0 \ran}$ and the isomorphism classes of finite dimensional simple $U_q (\Fsl_2^{\lan 0 \ran}[x])$-modules 
given by $\vf \mapsto L(\bu^{\lan 0 \ran}(\vf))$. 
\end{thm}

\begin{cor}
\label{Cor mult poly}
For $\vf, \psi \in \CC[x]^{\lan 0 \ran}$, 
let $v_0 \in L(\bu^{\lan 0 \ran}(\vf))$ 
(resp. $w_0 \in L (\bu^{\lan 0 \ran} (\psi))$) 
be a highest weight vector. 
Let $V(\vf, \psi)$ be a $U_q^{\lan 0 \ran}$-submodule of $L(\bu^{\lan 0 \ran}(\vf)) \otimes L (\bu^{\lan 0 \ran}(\psi))$ 
generated by $v_0 \otimes w_0$. 
Then $V(\vf, \psi)$ is a highest weight module of the highest weight $\bu^{\lan 0 \ran}(\vf \psi)$. 
In particular, we have $L(\bu^{\lan 0 \ran}(\vf)) \otimes L(\bu^{\lan 0 \ran}(\psi)) \cong L(\bu^{\lan 0 \ran} (\vf \psi))$ 
if $L(\bu^{\lan 0 \ran}(\vf)) \otimes L(\bu^{\lan 0 \ran}(\psi))$ is simple. 

\begin{proof}
For $\vf, \psi \in \CC[x]^{\lan 0 \ran}$, 
write $\vf$ and $\psi$ 
as   
$\vf= \ve (x-\g_1) (x-\g_2) \dots (x-\g_k)$ 
and 
$\psi = \ve' (x-\xi_1) (x-\xi_2) \dots (x-\xi_l)$ 
respectively. 
Let $v$ (resp. $w$) be a basis of one-dimensional 
$U_q^{\lan 0 \ran}$-module $\cD_{\ve}^{\lan 0 \ran}$ (resp. $\cD_{\ve'}^{\lan 0 \ran}$), 
and 
let $v_0^{(i)} \in V_1^{\wt{\ev}_{\g_i}^{\lan 0 \ran}} $ ($1 \leq i \leq k$) 
(resp. $w_0^{(i)} \in V_1^{\wt{\ev}_{\xi_i}^{\lan 0 \ran}}$ ($1 \leq i \leq l$))
be a highest weight vector. 
As in a proof of Proposition \ref{Prop fin simple sl20 2}, 
we have 
$L(\bu^{\lan 0 \ran} (\vf)) \cong \Top V(\ve; \g_1, \dots, \g_k)$ 
and 
$L(\bu^{\lan 0 \ran} (\psi)) \cong \Top V (\ve' ; \xi_1,\dots, \xi_l)$. 
By the definition of $\cD_{\pm 1}^{\lan 0 \ran}$ given in \eqref{act one dim rep}, 
we can easily check that 
$M \otimes \cD_{\pm 1}^{\lan 0 \ran} \cong \cD_{\pm 1}^{\lan 0 \ran} \otimes M$ 
as $U_q^{\lan 0 \ran}$-modules for any $U_q^{\lan 0 \ran}$-module $M$. 
As a consequence, 
we see that the highest weight of $V(\vf, \psi)$ is same as the highest weight of 
$V(\ve \ve' ; \g_1,\dots, \g_k, \xi_1, \dots, \xi_l)$ given by $\bu^{\lan 0 \ran}(\vf \psi)$. 
\end{proof} 
\end{cor}

\para 
We define a map $\flat : \CC[x] \ra \CC [\w]$ ($\vf \mapsto \vf^{\flat} (\w)$) by 
\begin{align*} 
\vf^{\flat}(\w) =  (1 - \g_1 \w) (1- \g_2 \w) \dots (1 - \g_k \w) 
\end{align*} 
if $\vf = \b_{\vf} (x - \g_1) (x-\g_2) \dots (x-\g_k)$.  
Then, 
Theorem \ref{Thm simple sl20} together with Proposition \ref{Prop gen series pt(q)} 
implies the following corollary. 
\begin{cor}
\label{Cor Psi+ Q=0}
For $\vf  \in \CC[x]^{\lan 0 \ran}$, let $v_0$ be a highest weight vector 
of $L(\bu^{\lan 0 \ran} (\vf))$. 
Then we have 
\begin{align*}
\Psi^+(\w) \cdot v_0 
= \b_{\vf} q^{\deg \vf} \frac{\vf^{\flat} (q^{-2} \w)}{\vf^{\flat} (\w)} v_0. 
\end{align*}
(Note that $p_t(q)(\g_1, \g_2, \dots, \g_{k-1}, 0) = p_t(q) (\g_1, \g_2,\dots, \g_{k-1})$ by Lemma \ref{Lemma ptq} (\roi). )
\end{cor}

\begin{remark}
\label{Remark Drinfeld poly}
Let $\CC[x]^D$ be the set of polynomials over $\CC$ with an indeterminate variable $x$ 
whose constant term is equal to $1$. 
By \cite[Theorem 3.4]{CP91}, 
there is a bijection between $\CC[x]^D$ and isomorphism classes of finite dimensional simple $U_q(L\Fsl_2)$-modules of type $1$. 
We call elements of $\CC[x]^D$ Drinfeld polynomials. 
For $\vf \in \CC[x]^{D}$, 
let $L^D(\vf)$ be the corresponding finite dimensional simple $U_q(L\Fsl_2)$-module, 
We regard $L^D(\vf)$ as a $U_q^{\lan 0 \ran}$-module through the injection $\Theta^{\lan 0 \ran}$. 
Then, 
we see that  $L^D (\vf)$ is still simple as a $U_q^{\lan 0 \ran}$-module 
(cf. \cite[Remark 3.2]{MTZ}). 

Let $\sharp : \CC[x]^D \ra \CC[x]^{\lan 0 \ran}$ ($\vf \mapsto \vf^{\sharp}$) be the injective map given by 
\begin{align*}
((1- \g_1 x) (1-\g_2 x) \dots ( 1 - \g_k x))^{\sharp} = (x- \g_1) (x-\g_2) \dots (x-\g_k), 
\end{align*}
where $\g_i \not=0$ ($1\leq i \leq k$). 
Then we see that 
$L^D (\vf) \cong L(\bu^{\lan 0 \ran} (\vf^{\sharp}))$ as $U_q^{\lan 0 \ran}$-modules 
by Corollary \ref{Cor Psi+ Q=0} and \cite[Theorem 3.4]{CP91}. 
\end{remark}



\section{Finite dimensional simple modules of $U_q(\Fsl_2^{\lan Q \ran}[x])$} 

In this section, we classify the isomorphism classes of finite dimensional simple modules of the algebra 
$U_q^{\lan Q \ran} = U_q (\Fsl_2^{\lan Q \ran}[x])$ in the case of rank one and of $Q\not=0$. 

\para 
For $k \geq 0$, 
put $J^{\lan Q \ran}_{[k;0]} = q^{- k(k+1)} Q^k$, 
For $t=1,2,\dots, k$, we define the element $J_{[k;t]}^{\lan Q \ran} \in U_q^{\lan Q \ran}$ inductively by 
\begin{align}
\label{def JktQ}
J_{[k;t]}^{\lan Q \ran} =q^t \frac{1}{[t]} \sum_{z=1}^t (-1)^{z-1} 
	\big( J_z - q^{2 (k-t+z)} Q^{-z} J_0  + q^{k- 2 (t-z)} [k] Q^{-z} \big) J_{[k;t-z]}^{\lan Q \ran}.
\end{align}
For examples, we have 
\begin{align*}
&J_{[1;0]}^{\lan Q \ran} = q^{-2} Q, 
	\quad 
	J_{[1;1]}^{\lan Q \ran} =  1 - q J_0  + q^{-1} Q J_1, 
\\
& J_{[2;0]}^{\lan Q \ran} = q^{-6} Q^2,  
	\quad 
	J_{[2;1]}^{\lan Q \ran} = q^{-3} [2] Q - q^{-1} Q J_0 + q^{-5} Q^2 J_1, 
\\
&J_{[2;2]}^{\lan Q \ran} 
=   1 - \frac{2q^2 + 1}{[2]} J_0 + q^{-2} [2] Q J_1 - \frac{q^{-4}}{[2]} Q^2 J_2  
	+ \frac{q^3}{[2]} J_0 J_0 -  Q J_0 J_1 + \frac{q^{-3}}{[2]} Q^2 J_1 J_1. 
\end{align*}

We have the following relations in $U_q^{\lan Q \ran}$. 

\begin{lem}
\label{Lemma X0+k X0-k+1}
For $k \in \ZZ_{>0}$, we have 
\begin{align*}
X_0^{+(k)} X_0^{-(k+1)} 
\equiv \sum_{z=0}^k (-1)^{k-z} X_z^- (K^+)^k J^{\lan Q \ran}_{[k; k-z]} 
\mod \fX_+, 
\end{align*}
where $\fX_+$  is the left ideal of $U_q^{\lan Q \ran}$ generated by $\{X_t^+ \mid t \geq 0\}$. 

\begin{proof}
See Appendix \ref{proof rel X0+k X0-k+1}.  
\end{proof}
\end{lem}

By using the above lemma, 
we have the following condition for $L(\bu)$ to be finite dimensional. 

\begin{prop}
\label{Prop finite cond UqQ sl2}
Assume that $Q \not=0$. 
For $\bu= (\la, (u_t)_{t>0}) \in \CC^{\times} \times \prod_{t>0} \CC$, 
the simple $U_q^{\lan Q \ran}$-module $L(\bu)$ is finite dimensional 
if and only if 
there exist  $k \in \ZZ_{\geq 0}$, $\b \in \CC^{\times}$ and $\g_1,\g_2,\dots, \g_k \in \CC$ 
such that 
\begin{align}
\label{hw Lu} 
\la = \b q^k,  
\quad 
u_t = \begin{cases} 
	\wt{\b} Q^{-t} & \text{ if } k=0, 
	\\
	p_t^{\lan Q \ran} (q;\b) (\g_1,  \g_2, \dots, \g_k) & \text{ if } k \not=0 
	\end{cases} 
\quad (t>0), 
\end{align}
where we put $\wt{\b} = (q-q^{-1})^{-1} (1- \b^{-2})$. 

\begin{proof}
We prove the only if part. 
Suppose that $L(\bu)$ is finite dimensional. 
Let $v_0 \in L(\bu)$ be a highest weight vector. 
By investigating the eigenvalues for the action of $K^+$, 
there exists $k \in \ZZ_{\geq 0}$ 
such that $X_0^{-(k)} \cdot v_0 \not=0$ and $X_0^{-(k+1)} \cdot v_0 =0$ since $L(\bu)$ is finite dimensional. 

In the case where $k=0$, we can easily check that $L(\bu)$ is one-dimensional. 
Then the condition \eqref{hw Lu} follows from Proposition \ref{Prop one dim rep}. 

Assume that $k >0$. 
Put $\b = \la q^{-k}$. 
By Lemma \ref{Lemma ptQ alg ind}, 
there exist  
$\g_1,\g_2,\dots, \g_k \in \CC$ such that 
$u_t = p_t^{\lan Q \ran} (q;\b)(\g_1,\g_2,\dots, \g_k)$ 
for $t =1,2, \dots, k$. 
By the induction on $t$ using \eqref{def JktQ} and Lemma \ref{Lemma ptQqb} (\roi), 
we can show that 
\begin{align}
\label{JktQ v0} 
J^{\lan Q \ran}_{[k;t]} \cdot v_0 
= q^{-k (k+1)} Q^k e_t (\g_1,\g_2,\dots, \g_k) v_0 
\quad (1 \leq t \leq k) 
\end{align}
where we note that $(q^{2k } J_0 - q^{k}[k]) \cdot v_0 = \wt{\b} v_0$. 

By Lemma \ref{Lemma X0+k X0-k+1} and the relation (Q6), 
for $t \geq 0$, we have 
\begin{align*}
&0 = X_t^+ X_0^{+(k)} X_0^{-(k+1)} \cdot v_0 
= \sum_{z=0}^k (-1)^{k-z} (J_{t+z} - Q J_{t+z+1}) (K^+)^{k+1} J^{\lan Q \ran}_{[k; k-z]} \cdot v_0. 
\end{align*}
Note that $(K^+)^{k+1} J^{\lan Q \ran}_{[k; k-z]} \cdot v_0 = \b^{k+1} Q^k e_{k-z}(\g_1,\dots, \g_k)$ 
by the choice of $\b$ and \eqref{JktQ v0}, 
then the above equation implies that 
\begin{align}
\label{Jk+t+1 v0} 
J_{k +t+1} \cdot v_0 
= Q^{-1} J_{k+t} \cdot v_0 + \sum_{z=0}^{k-1} (-1)^{k-z+1} (J_{t+1+z} - Q^{-1} J_{t+z}  ) e_{k-z} (\g_1,\dots, \g_k) \cdot v_0. 
\end{align}
Then we can show that 
$u_{k+t} = p_{k+t}^{\lan Q \ran} (q;\b) (\g_1,\dots, \g_k)$ for $t >0$ 
by the induction on $t$ using \eqref{Jk+t+1 v0} and Lemma \ref{Lemma ptQqb} (\roii). 

Next we prove the if part. 
Recall the algebra homomorphism $\D_r^{\lan Q \ran} : U_q^{\lan Q \ran} \ra U_q^{\lan Q \ran} \otimes U_q^{\lan 0 \ran}$ 
given in Theorem \ref{Thm DrQ}.
In a similar way  as in \cite[Proposition 4.4]{CP91}, 
we can show that 
\begin{align}
\label{DrQ mod}
\begin{split} 
&\D_r^{\lan Q \ran} (X_{t'}^+) 
	\equiv  X_{t'}^+ \otimes K^+ + 1 \otimes X_{t'}^+ + (q-q^{-1}) \sum_{z=1}^{t'} X_{t'-z}^+ \otimes K^+ J_z  
	\mod  \fX_{+2}^{\lan Q \ran} \otimes \fX_-^{\lan 0 \ran}, 
\\
& \D_r^{\lan Q \ran} (J_t) 
	\equiv J_t \otimes 1 + 1 \otimes J_t + (q-q^{-1}) \sum_{z=1}^{t-1} J_z \otimes J_{t-z} 
	\mod \fX_+^{\lan Q \ran} \otimes \fX_-^{\lan 0 \ran}, 
\end{split}
\end{align}
for $t' \geq 0$ and $t >0$, 
where 
$\fX_{+2}^{\lan Q \ran}$ (resp. $\fX_+^{\lan Q \ran}$, $\fX_-^{\lan Q \ran}$) 
is the left ideal of $U_q^{\lan Q \ran}$ 
generated by $\{X_s^+ X_{s'}^+ \mid s,s' \geq 0\}$ (resp. $\{X_s^+ \mid s \geq 0\}$, $\{X_s^- \mid s \geq 0 \}$).

For $\g_1, \g_2,\dots, \g_k \in \CC$, put $\vf=(x-\g_1)(x-\g_2) \dots (x-\g_k) \in \CC[x]^{\lan 0 \ran}$. 
Then $L(\bu^{\lan 0 \ran}(\vf))$ is finite dimensional simple $U_q^{\lan 0 \ran}$-module by Theorem \ref{Thm simple sl20}. 
Let $v_0 \in L(\bu^{\lan 0 \ran}(\vf))$ be a highest weight vector. 
For 
$\b \in \CC^{\times}$, 
we take the  one-dimensional $U_q^{\lan Q \ran}$-module $\cD_{\b}^{\lan Q \ran}= \CC v$. 
We consider the $U_q^{\lan Q \ran}$-module $\cD_{\b}^{\lan Q \ran} \otimes L(\bu^{\lan 0 \ran}(\vf))$ through the homomorphism 
$\D_r^{\lan Q \ran}$. 
Let $V(\b; \g_1,\dots, \g_k)$ be the $U_q^{\lan Q \ran}$-submodule of 
$\cD_{\b}^{\lan Q \ran} \otimes L(\bu^{\lan 0 \ran}(\vf))$ generated by $v \otimes v_0$. 
Then we have $X_t^+ \cdot (v \otimes v_0) =0$ for all $t \geq 0$  by \eqref{DrQ mod}. 
On the other hand, we have 
\begin{align*}
&J_t \cdot (v \otimes v_0)
\\
&=\big( \wt{\b} Q^{-t} + p_t(q) (\g_1,\dots, \g_k) + (q-q^{-1}) \sum_{z=1}^{t-1} \wt{\b} Q^{-z} p_{t-z}(q)(\g_1,\dots, \g_k) \big) 
	v \otimes v_0 
\end{align*}
for $t >0$ by \eqref{DrQ mod}, 
where $\wt{\b} = (q-q^{-1})^{-1} (1 - \b^{-2})$. 
Thus, we have 
$J_t \cdot (v \otimes v_0) = p_t^{\lan Q \ran} (q; \b)(\g_1,\dots, \g_k)$ by \eqref{def ptQqb} . 
We also see that $K^+ \cdot (v \otimes v_0) = \b q^k  v \otimes v_0$. 
As a consequence, 
the $U_q^{\lan Q \ran}$-module $V(\b;\g_1,\dots, \g_k)$ is a highest weight module of the highest weight $\bu$ 
given by \eqref{hw Lu}, 
and $L(\bu)$ is a quotient of $V(\b; \g_1,\dots, \g_k)$. 
Thus $L(\bu)$ is finite dimensional. 
\end{proof}
\end{prop} 

\para 
In order to give a correspondence between the elements of $\CC[x]$ and  finite dimensional simple $U_q^{\lan Q \ran}$-modules, 
We define a map $\bu^{\lan Q \ran} : \CC[x] \setminus \{0\} \ra \CC^{\times} \times \prod_{t>0} \CC$ by 
\begin{align*}
\bu^{\lan Q \ran}(\vf) 
= \begin{cases}
	\big( \b_{\vf},  ( \wt{\b}_{\vf}Q^{-t} )_{t>0} \big) & \text{ if } \deg \vf =0,
	\\
	\big( \b_{\vf} q^{\deg \vf}, (p_t^{\lan Q \ran}(q;\b_\vf) (\g_1,\g_2,\dots, \g_k) )_{t>0} \big) 
	& \text{ if } \deg \vf >0
	\end{cases} 
\end{align*}
for $\vf =\b_{\vf} (x-\g_1)(x-\g_2) \dots (x-\g_k) \in \CC[x] \setminus \{0\}$, 
where 
$\wt{\b}_{\vf} =(q-q^{-1})^{-1} (1- \b_{\vf}^{-2})$. 
Unfortunately, the map $\bu^{\lan Q \ran}$ is not injective. 
In order to obtain  an  index set of the isomorphism classes of finite dimensional simple $U_q^{\lan Q \ran}$-modules, 
we take a subset $\CC[x]^{\lan Q \ran}$ of $\CC[x]$ as 
\begin{align*}
\CC[x]^{\lan Q \ran} = \{ \vf \in \CC[x] \setminus \{0\} \mid \b_{\vf}^{-2} Q^{-1} \text{ is not a root of }\vf\}. 
\end{align*}
Then we have the following proposition. 

\begin{prop}
\label{Lemma uQ vf = UQ vf'}
\begin{enumerate} 
\item 
For $\vf, \vf' \in \CC[x] \setminus \{0\}$ such that $\deg \vf \geq \deg \vf'$, 
we have that 
$\bu^{\lan Q \ran}(\vf) = \bu^{\lan Q \ran}(\vf')$ 
if and only if 
\begin{align*}
\vf = q^{- (\deg \vf - \deg \vf')} \vf' \prod_{z=1}^{\deg \vf - \deg \vf'} ( x - q^{- 2 (z-1)} \b_{\vf}^{-2} Q^{-1}).
\end{align*}

\item 
The restriction of $\bu^{\lan \bQ \ran}$ to $\CC[x]^{\lan Q \ran}$ is injective. 
Moreover, 
for any $\vf\not=0 \in \CC[x] $, 
there exists  the unique 
$\vf' \in \CC[x]^{\lan Q \ran}$ such that $\bu^{\lan Q \ran} (\vf) = \bu^{\lan Q \ran} (\vf')$. 
\end{enumerate}

\begin{proof}
We prove the statement (\roi) by the induction on $\deg \vf - \deg \vf'$. 
For $\vf, \vf' \in \CC[x] \setminus \{0\}$, 
write 
$\vf=\b_{\vf} (x-\g_1) (x-\g_2) \dots (x-\g_k)$  
and 
$\vf' = \b_{\vf'} (x- \g'_1)(x-\g'_2) \dots (x-\g'_l)$. 

First, we consider the case where $k=l$. 
In this case, it is clear that $\bu^{\lan Q \ran} (\vf) = \bu^{\lan Q \ran}(\vf')$ if $\vf = \vf'$. 
Assume that $\bu^{\lan Q \ran}(\vf) = \bu^{\lan Q \ran}(\vf')$. 
If $k=l=0$, 
we can easily check that $\vf = \vf'$ from definitions. 
If $k=l >0$, 
the assumption $\bu^{\lan Q \ran} (\vf) = \bu^{\lan Q \ran}(\vf')$ implies that 
$\b_{\vf} = \b_{\vf'}$ and 
\begin{align*}
p_t(q) (\g_1,\dots, \g_k) 
&= p_t(q) (\g'_1, \dots, \g'_k)
\\
& \quad + (q-q^{-1}) \sum_{z=1}^{t-1} \wt{\b}_{\vf} Q^{-t+z} \big( (p_z(q) (\g'_1,\dots, \g'_k) - p_z(q) (\g_1,\dots, \g_k) \big)
\end{align*}
for $t >0$. 
Using this equation, 
we can prove that 
$p_t(q) (\g_1,\dots, \g_k) = p_t(q) (\g_1', \dots, \g_k')$ for $t>0$. 
Then we have $\vf = \vf'$ in the same way as the proof of Lemma \ref{Lemma inj u0}. 
As a consequence, 
we have $\vf = \vf'$ if $\bu^{\lan Q \ran}(\vf)= \bu^{\lan Q \ran}(\vf')$.

Next we consider the case where $k =l+1$. 
If $\vf = q^{-1} \vf' (x - \b_{\vf}^{-2} Q^{-1})$, 
we have 
$\b_{\vf} = q^{-1} \b_{\vf'}$ and 
$p_t^{\lan Q \ran}(q;\b_{\vf})( \g_1,\dots, \g_k) = p_t^{\lan Q \ran} (q; q^{-1} \b_{\vf'}) (\g'_1, \dots, \g'_l, \b_{\vf}^{-2} Q^{-1})$ 
for $t>0$. 
On the other hand, we have 
\begin{align*}
p_t^{\lan Q \ran} (q; q^{-1} \b_{\vf'}) (\g'_1, \dots, \g'_l, \b_{\vf}^{-2} Q^{-1})
&= p_t (q) (\g'_1, \dots, \g'_l, \b_{\vf}^{-2} Q^{-1}) + \frac{1- q^2 \b_{\vf'}^{-2}}{q-q^{-1}} Q^{-t} 
	\\ & \quad 
	+ \sum_{z=1}^{t-1} (1- q^2 \b_{\vf'}^{-2}) Q^{-t+z} p_z(q) (\g'_1, \dots, \g'_l, \b_{\vf}^{-2} Q^{-1}) 
\end{align*}
by the definition \eqref{def ptQqb}. 
Applying Lemma \ref{Lemma ptq} (\roi) to the right-hand side of the above equation, 
we have 
\begin{align*}
& p_t^{\lan Q \ran}(q;\b_{\vf})( \g_1,\dots, \g_k)
\\
&= p_t^{\lan Q \ran} (q; q^{-1} \b_{\vf'}) (\g'_1, \dots, \g'_l, \b_{\vf}^{-2} Q^{-1})
\\
&= p_t(q) (\g'_1, \dots, \g'_l)  + \frac{1- \b_{\vf'}^{-2}}{q-q^{-1}} Q^{-t}  
	+ \sum_{z=1}^{t-1} (1 - \b_{\vf'}^{-2} ) Q^{- t +z}  p_z(q) (\g'_1,\dots, \g'_l)  
\\
&= p_t^{\lan Q \ran} (q;  \b_{\vf'}) (\g'_1, \dots, \g'_l). 
\end{align*}
Thus, we have $\bu^{\lan Q \ran} (\vf) = \bu^{\lan Q \ran}(\vf')$ if $\vf = q^{-1} \vf'(x- \b_{\vf}^{-2} Q^{-1})$. 

Assume that $\bu^{\lan Q \ran} (\vf) = \bu^{\lan Q \ran} (\vf')$. 
Then we have 
$q \b_{\vf}  = \b_{\vf'}$ and  
\begin{align}
\label{ptQk ptQl}
p_t^{\lan Q \ran} (q; \b_{\vf}) (\g_1,\dots, \g_k) = p_t^{\lan Q \ran} (q; \b_{\vf'}) (\g'_1,\dots, \g'_l)
\end{align} 
for $t>0$. 
We note that $\wt{\b}_{\vf'} = (q-q^{-1})^{-1} (1 - q^{-2} \b_{\vf}^{-2}) = \wt{\b}_{\vf} + q^{-1} \b_{\vf}^{-2}$ 
since $\b_{\vf'}= q \b_{\vf}$. 
Then, by applying the definition \eqref{def ptQqb} to both sides of \eqref{ptQk ptQl}, we have 
\begin{align*}
&p_t (q) (\g_1, \dots, \g_k) 
\\
&= p_t (q) (\g'_1, \dots, \g'_l) + q^{-1} \b_{\vf}^{-2} Q^{-t} 
	+ (q-q^{-1}) \sum_{z=1}^{t-1} ( \wt{\b}_{\vf} + q^{-1} \b_{\vf}^{-2}) Q^{-t+z} p_z(q) (\g'_1,\dots, \g'_l) 
	\\& \quad 
	- (q-q^{-1}) \sum_{z=1}^{t-1} \wt{\b}_{\vf} Q^{-t+z} p_z(q) (\g_1,\dots, \g_k) 
\\
&= p_t(q) (\g'_1, \dots, \g'_l) + q^{-1} \b_{\vf}^{-2t} Q^{-t}  
	+ q^{-1} (q-q^{-1}) \sum_{z=1}^{t-1} \b_{\vf}^{-2(t-z)} Q^{-(t-z)} p_z(q) (\g'_1,\dots, \g'_l)
	\\ & \quad 
	+ q^{-1} (\b_{\vf}^{-2} - \b_{\vf}^{-2t}) Q^{-t}  
	+ (q-q^{-1}) \sum_{z=1}^{t-1} ( \wt{\b}_{\vf} + q^{-1} \b_{\vf}^{-2} - q^{-1} \b_{\vf}^{-2(t-z)} ) Q^{-t+z} p_z(q) (\g'_1,\dots, \g'_l)
	\\& \quad 
	- (q-q^{-1}) \sum_{z=1}^{t-1} \wt{\b}_{\vf} Q^{-t+z} p_z(q) (\g_1,\dots, \g_k). 
\end{align*}
Applying Lemma \ref{Lemma ptq} (\roi), we have
\begin{align}
\label{ptqg ptqg'}
\begin{split}
&p_t (q) (\g_1, \dots, \g_k) 
\\
&= p_t(q) (\g'_1, \dots, \g'_l, \b_{\vf}^{-2} Q^{-1}) 
	\\ & \quad 
	+ q^{-1} (\b_{\vf}^{-2} - \b_{\vf}^{-2t}) Q^{-t}  
	+ (q-q^{-1}) \sum_{z=1}^{t-1} ( \wt{\b}_{\vf} + q^{-1} \b_{\vf}^{-2} - q^{-1} \b_{\vf}^{-2(t-z)} ) Q^{-t+z} p_z(q) (\g'_1,\dots, \g'_l)
	\\& \quad 
	- (q-q^{-1}) \sum_{z=1}^{t-1} \wt{\b}_{\vf} Q^{-t+z} p_z(q) (\g_1,\dots, \g_k). 
\end{split}
\end{align}
Then, we can prove that 
$p_t(q) (\g_1, \dots, \g_k) = p_t(q) (\g'_1, \dots, \g'_l, \b_{\vf}^{-2} Q^{-1})$ for $t>0$ 
by the induction on $t$ using \eqref{ptqg ptqg'} with Lemma \ref{Lemma ptq} (\roi). 
This equations imply that 
$(x-\g_1) \dots (x-\g_k) = (x-\g'_1)\dots(x-\g'_l) (x- \b_{\vf}^{-2} Q^{-1})$ 
in a similar way as in the proof of Lemma \ref{Lemma inj u0}. 
As a consequence, 
we have $\vf= q^{-1} \vf' (x- \b_{\vf}^{-2} Q^{-1})$ if $\bu^{\lan Q \ran} (\vf) = \bu^{\lan Q \ran}(\vf')$. 

Finally, we consider the case where $k >l+1$ by the induvtion on $k-l$. 
Put $\vf'' = q^{-1} \vf'(x - q^2 \b_{\vf'}^{-2} Q^{-1})$, 
and we have $\b_{\vf''} = q^{-1} \b_{\vf'}$. 
Then, we have $\bu^{\lan Q \ran} (\vf'') = \bu^{\lan Q \ran} (\vf')$ by the above argument.  
On the other hand, 
we have 
\begin{align*}
\bu^{\lan Q \ran} (\vf) = \bu^{\lan Q \ran} (\vf'')
\LRa 
\vf 
&= q^{- (\deg \vf - \deg \vf'')} \vf'' \prod_{z=1}^{\deg \vf - \deg \vf''} ( x - q^{- 2 (z-1)} \b_{\vf}^{-2} Q^{-1}) 
\\
&= q^{- \deg \vf - \deg \vf'} \vf' \prod_{z=1}^{\deg \vf - \deg \vf'} ( x - q^{- 2 (z-1)} \b_{\vf}^{-2} Q^{-1}) 
\end{align*}
by the induction hypothesis. 
Thus, we have the statement (\roi). 

The statement (\roii) follows from the statement (\roi) and the definition of $\CC[x]^{\lan Q \ran}$. 
\end{proof}
\end{prop}


Proposition  \ref{Prop finite cond UqQ sl2}  and Proposition \ref{Lemma uQ vf = UQ vf'} 
imply the following theorem.  

\begin{thm}
\label{Theorem simple sl2Q}
Assume that $Q \not=0$. 
There exists the bijection between $\CC[x]^{\lan Q \ran}$ 
and the isomorphism classes of finite dimensional simple 
$U_q (\Fsl_2^{\lan Q \ran}[x])$-modules 
given by $\vf \mapsto L(\bu^{\lan Q \ran}(\vf))$. 
\end{thm}

\begin{cor}
\label{Cor Psi+ Q not=0}
For $\vf \in \CC[x]^{\lan Q \ran}$, 
let $v_0$ be a highest weight vector of $L(\bu^{\lan Q \ran}(\vf))$. 
Then we have 
\begin{align*}
\Psi^+(\w) \cdot v_0 
= q^{\deg \vf} \frac{\vf^{\flat} (q^{-2} \w)}{\vf^{\flat} (\w)} (\b_{\vf}^{-1} - Q \b_{\vf} \w^{-1}) v_0. 
\end{align*}
\begin{proof}
For $\vf \in \CC[x]^{\lan Q \ran}$, 
write $\vf = \b_{\vf} (x-\g_1) (x-\g_2) \dots (x- \g_k)$. 
By \eqref{Def Psi Q not=0}, 
for $t >0$, 
we have 
\begin{align*}
&\Psi_{t}^+ \cdot v_0 
\\
&= (q-q^{-1}) K^+ (J_{t} - Q J_{t+1}) \cdot v_0 
\\
&= (q-q^{-1}) \b_{\vf} q^{\deg \vf} \Big( p_t^{\lan Q \ran} (q; \b_{\vf}) (\g_1,\dots, \g_k) - Q p_{t+1}^{\lan Q \ran} (q; \b_{\vf}) (\g_1, \dots, \g_k) \Big) v_0. 
\end{align*}
Applying \eqref{def ptQqb}, we have 
\begin{align*}
\Psi_t^+ \cdot v_0 
= (q-q^{-1}) \b_{\vf} q^{\deg \vf} \big( \b_{\vf}^{-2} p_t(q) (\g_1,\dots, \g_k) 
	- Q p_{t+1} (q) (\g_1,\dots, \g_k) \big) v_0 
\end{align*} 
for $t >0$, 
where we note that $1 - (q-q^{-1}) \wt{\b}_{\vf} = \b_{\vf}^{-2}$. 
We also have 
\begin{align*}
\Psi_{-1}^+ \cdot v_0 = - Q \b_{\vf} q^{\deg \vf} v_0, 
\quad 
\Psi_{0}^+ \cdot v_0 = \b_{\vf} q ^{\deg \vf} \big( \b_{\vf}^{-2} - (q-q^{-1}) Q p_1(q)(\g_1,\dots, \g_k) \big) v_0. 
\end{align*}
Thus, we have 
\begin{align*}
&\Psi^+(\w) \cdot v_0 
\\
&= \Big\{ - Q \b_{\vf} q^{\deg \vf} \w^{-1} 
	+ \b_{\vf} q ^{\deg \vf} \big( \b_{\vf}^{-2} 
		- (q-q^{-1}) Q p_1(q)(\g_1,\dots, \g_k) \big)
	\\ & \qquad 
	+ \sum_{t >0}  (q-q^{-1}) \b_{\vf} q^{\deg \vf} \big( \b_{\vf}^{-2} p_t(q) (\g_1,\dots, \g_k) 
	- Q p_{t+1} (q) (\g_1,\dots, \g_k) \big) \w^t 
	\Big\} v_0 
\\
&= q^{ \deg \vf} \Big( 1 + (q-q^{-1}) \sum_{t>0} p_t(q) (\g_1,\dots, \g_k) \w^t \Big) \big( \b_{\vf}^{-1} - Q \b_{\vf} \w^{-1} \big) v_0 
\\
&= q^{\deg \vf} \frac{\vf^{\flat} (q^{-2} \w)}{\vf^{\flat} (\w)} (\b_{\vf}^{-1} - Q \b_{\vf} \w^{-1}) v_0 
\end{align*}
by Proposition \ref{Prop gen series pt(q)}. 
\end{proof}
\end{cor}



\section{Finite dimensional simple modules of $U_q(\Fsl_n^{\lan \bQ \ran}[x])$} 

In this section, we classify the isomorphism classes of finite dimensional simple $U_q (\Fsl_n^{\lan \bQ \ran}[x])$-modules.  
In this section, 
we denote by $U_q^{\lan \bQ \ran}$ the $(q, \bQ)$-current algebra $U_q (\Fsl_n^{\lan \bQ \ran}[x])$ of rank $n-1$ 
with a parameter $\bQ=(Q_1,\dots, Q_{n-1}) \in \CC^{n-1}$. 


\para 
For each $i \in I$, 
we can easily check that 
there exists the algebra homomorphism 
$\iota_i : U_q (\Fsl_2^{\lan Q_i \ran}[x]) \ra U_q^{\lan \bQ \ran}$ 
such that 
$\iota_i (X_t^{\pm}) = X_{i,t}^{\pm}$, 
$\iota_i (J_t) = (J_{i,t})$ 
and $\iota_i(K^{\pm}) = K_i^{\pm}$. 
For $\bu=((\la_i, (u_{i,t})_{t>0}))_{i \in I} \in ( \CC^{\times} \times \prod_{t>0} \CC)^I$, 
we regard the highest weight simple $U_q^{\lan \bQ \ran}$-module 
$L(\bu)$ as a $U_q (\Fsl_2^{\lan Q _i \ran}[x])$-module 
through the homomorphism $\iota_i$. 
Let $v_0 \in L(\bu)$ be a highest weight vector. 
Then we can easily check that 
the $U_q(\Fsl_2^{\lan Q_i \ran}[x])$-submodule of $L(\bu)$ generated by $v_0$ 
is a highest weight simple $U_q (\Fsl_2^{\lan Q_i \ran}[x])$-module of the highest weight $(\la_i, (u_{i,t})_{t>0})$. 
Thus, Theorem \ref{Thm simple sl20} and Theorem \ref{Theorem simple sl2Q} imply the following proposition. 

\begin{prop}
\label{Prop simple slnQ 1} 
For $\bu=((\la_i, (u_{i,t})_{t>0}))_{i \in I} \in ( \CC^{\times} \times \prod_{t>0} \CC)^I$, 
if the highest weight simple  $U_q^{\lan \bQ \ran}$-module 
$L(\bu)$ is finite dimensional, 
then there exists 
$(\vf_i)_{i \in I} \in \prod_{i \in I} \CC[x]^{\lan Q _i \ran}$ 
such that 
$((\la_i, (u_{i,t})_{t>0}))_{i\in I} = (\bu^{\lan Q_i \ran}(\vf_i))_{i \in I}$. 
\end{prop}


\para 
Recall the algebra homomorphism 
$\D_r^{\lan \bQ \ran} : U_q^{\lan \bQ \ran} \ra U_q^{\lan \bQ \ran} \otimes U_q^{\lan \mathbf{0} \ran}$ 
given in Theorem \ref{Thm DrQ}.
In a similar way  as in \cite[Proposition 4.4]{CP91}, 
we can show that 
\begin{align}
\label{DrQ mod slnQ}
\begin{split} 
&\D_r^{\lan \bQ \ran} (X_{i,t'}^+) 
	\equiv  X_{i,t'}^+ \otimes K_i^+ + 1 \otimes X_{i,t'}^+ + (q-q^{-1}) \sum_{z=1}^{t'} X_{i,t'-z}^+ \otimes K_i^+ J_{i,z}  
	\mod \fX_{+2}^{\lan \bQ \ran} \otimes \fX_-^{\lan \mathbf{0} \ran}, 
\\
& \D_r^{\lan \bQ \ran} (J_{i,t}) 
	\equiv J_{i,t} \otimes 1 + 1 \otimes J_{i,t} + (q-q^{-1}) \sum_{z=1}^{t-1} J_{i,z} \otimes J_{i,t-z} 
	\mod \fX_+^{\lan \bQ \ran} \otimes \fX_-^{\lan \mathbf{0} \ran}, 
\end{split}
\end{align}
for $i \in I$, $t' \geq 0$ and $t >0$, 
where 
$\fX_{+2}^{\lan \bQ \ran}$ (resp. $\fX_+^{\lan \bQ \ran}$, $\fX_-^{\lan \bQ \ran}$) is the left ideal of $U_q^{\lan \bQ \ran}$ 
generated by $\{X_{j,s}^+ X_{j',s'}^+ \mid (j,s), (j',s') \in I \times \ZZ_{\geq 0}\}$ 
	(resp. $\{X_{j,s}^+ \mid (j,s) \in I \times \ZZ_{\geq 0}\}$, $\{ X_{j,s}^- \mid  (j,s) \in I \times \ZZ_{\geq 0} \}$). 

For $ (\vf_i)_{i \in I} \in \prod_{i\in I}\CC[x]^{\lan Q_i \ran}$, 
write  
$\vf_i$ as $\vf_i =\b_{\vf_i} (x- \g_1^{(i)}) (x- \g_2^{(i)}) \dots (x- \g_{k_i}^{(i)})$ for each $i \in I$. 
Put $\g_{i,p} = q^{i-2} \g_p^{(i)}$ for $i \in I$ and $1\leq p \leq k_i$. 
We consider the evaluation module $V(\w_i)^{\ev_{\g_{i,p}}^{\lan \mathbf{0} \ran}}$ 
at $\g_{i,p}$ for $i \in I$ and $1\leq p \leq k_i$.  
Let $v_p^{(i)} \in V(\w_i)^{\ev_{\g_{i,p}}^{\lan \mathbf{0} \ran}}$ be a highest weight vector. 
We also consider the one-dimensional $U_q^{\lan \bQ \ran}$-module 
$\cD_{\Bb}^{\lan \bQ \ran} = \CC v$, where $\Bb =(\b_{\vf_i})_{i \in I} \in \prod_{i\in I} \BB^{\lan Q_i \ran}$. 
Then we have the $U_q^{\lan \bQ \ran}$-module 
$\cD_{\Bb}^{\lan \bQ \ran} \otimes \bigotimes_{i \in I} \bigotimes_{p=1}^{k_i} V(\w_i)^{\ev_{\g_{i,p}}^{\lan \mathbf{0} \ran}}$
through the algebra homomorphisms $\D_r^{\lan \bQ \ran}$ and $\D^{\lan \mathbf{0} \ran}$. 
Let $V((\vf_i)_{i \in I})$ be the $U_q^{\lan \bQ \ran}$-submodule of 
$\cD_{\Bb}^{\lan \bQ \ran} \otimes \bigotimes_{i \in I} \bigotimes_{p=1}^{k_i} V(\w_i)^{\ev_{\g_{i,p}}^{\lan \mathbf{0} \ran}}$ 
generated by $v \otimes (v_1^{(1)} \otimes \dots \otimes v_{k_1}^{(1)}) \otimes 
\dots \otimes (v_1^{(n-1)} \otimes \dots \otimes v_{k_{n-1}}^{(n-1)})$. 
By Proposition \ref{Prop ev module} and \eqref{DrQ mod slnQ}, 
we can show that 
$V((\vf_i)_{i \in I})$ is a highest weight $U_q^{\lan \bQ \ran}$-module of the highest weight 
$(\bu^{\lan Q_i \ran}(\vf_i))_{i \in I}$ 
in a similar way as in the proofs of Proposition \ref{Prop fin simple sl20 2} and of Proposition \ref{Prop finite cond UqQ sl2}. 
Thus, the highest weight simple $U_q^{\lan \bQ \ran}$-module 
$L((\bu^{\lan Q_i \ran}(\vf_i))_{i \in I})$ is a quotient of $V(( \vf_i)_{i \in I})$. 
In particular, it is finite dimensional. 
As a consequence, we have the following proposition. 


\begin{prop}
\label{Prop simple slnQ 2}
For $ (\vf_i)_{i \in I} \in \prod_{i\in I}\CC[x]^{\lan Q_i \ran}$, 
the highest weight simple $U_q^{\lan \bQ \ran}$-module 
$L((\bu^{\lan Q_i \ran}(\vf_i))_{i \in I})$  is finite dimensional. 
\end{prop}

We have the following theorem by Proposition \ref{Prop simple slnQ 1} and Proposition \ref{Prop simple slnQ 2}.

\begin{thm}
\label{Thm  slnQ}
There exists the bijection between $\prod_{i \in I} \CC[x]^{\lan Q_i \ran}$ 
and the isomorphism classes of finite dimensional simple 
$U_q( \Fsl_n^{\lan \bQ \ran}[x])$-modules 
given by $(\vf_i)_{i \in I} \mapsto L((\bu^{\lan Q_i \ran}(\vf_i))_{i \in I})$. 
\end{thm}

\appendix 

\section{A proof of Theorem \ref{Theorem PBW}}
\label{Proof PBW}

\para 
Let $\AA = \CC[v,v^{-1}]$ be the Laurent polynomial ring over $\CC$ with an indeterminate element $v$, 
and let $\KK= \CC(v)$ be the quotient field of $\AA$. 
We also consider the localization $\CC[v]_{(v=1)}$ of the polynomial ring $\CC[v]$ at $v=1$. 

For $\XX \in \{\AA, \KK, \CC[v]_{(v=1)}\}$, 
we define an associative algebra $\cA^{\XX}$ over $\XX$ 
by generators $x_{i,t}$ ($(i,t) \in I \times \ZZ_{\geq 0}$) with  defining relations 
\begin{align}
\label{def rel cA}
\begin{split}
&x_{i,t+1} x_{j,s} - v^{a_{ij}} x_{j,s} x_{i,t+1} = v^{a_{ij}} x_{i,t} x_{j,s+1} - x_{j,s+1} x_{i,t}, 
\\
&[x_{i,t}, x_{j,s}]=0 \text{ if } j \not= i, i \pm1, 
\\
&x_{i\pm 1,u}  (x_{i,s} x_{i,t} + x_{i,t} x_{i,s}) + (x_{i,s} x_{i,t} + x_{i,t} x_{i,s}) x_{i \pm 1,u} 
	= (v+v^{-1}) (x_{i,s} x_{i\pm 1,u} x_{i,t} + x_{i,t} x_{i \pm 1, u} x_{i,s} ).
\end{split}
\end{align}

For $q \in \CC^{\times}$, 
Let $\cA$ be the  scalar extension 
$\CC \otimes_{\AA} \cA^{\AA}$ of $\cA^{\AA}$ 
through the ring homomorphism $\cA \ra \CC$ ($v \mapsto q$). 
Clearly, the algebra $\cA$ is isomorphic to an associateive algebra over $\CC$ 
generated by $x_{i,t}$ ($(i,t) \in I \times \ZZ_{\geq 0}$) with defining relations \eqref{def rel cA}, 
where we replace $v$ with $q$. 
Then we have the surjective algebra homomorphisms 
\begin{align*}
\pi^+ : \cA \ra U_{q, \bQ}^+ \, (x_{i,t} \mapsto X_{i,t}^+), 
\quad 
\pi^- : \cA^{\opp} \ra U_{q, \bQ}^- \, (x_{i,t} \mapsto X_{i,t}^-), 
\end{align*}
where $\cA^{\opp}$ is the opposite algebra of $\cA$. 

\para 
Let $Q = \bigoplus_{i \in I} \ZZ \a_i$ be the root lattice of $\Fsl_n$, 
and we put $Q^+ = \sum_{i \in I} \ZZ_{\geq 0} \a_i$. 
>From the definition, we see that the algebra $\cA^{\XX}$ 
is a $Q$-graded algebra with $\deg_Q (x_{i,t}) = \a_i$, 
and $\cA^{\XX}$ is also a $\ZZ$-graded algebra with $\deg(x_{i,t}) =t$. 
Then the algebra $\cA^{\XX}$ decomposes into 
\begin{align*}
\cA^{\XX} = \bigoplus_{\g \in Q^+} \bigoplus_{s \geq 0} \cA_{\g,s}^{\XX}, 
\quad 
\cA_{\g,s}^{\XX} :=\{x \in \cA^{\XX} \mid \deg_Q (x) =\g, \, \deg (x)=s\} 
\end{align*}
as $\XX$-modules. 
It is clear that, for each $(\g,s) \in Q^+ \times \ZZ_{\geq 0}$, 
the $\XX$-module $\cA_{\g,s}^{\XX}$ is generated by 
$\{x_{i_1t_1} x_{i_2,t_2} \dots x_{i_k,t_k} \mid \a_{i_1}+ \dots + \a_{i_k}=\g, \, t_1+\dots + t_k =s \}$, 
and $\cA_{\g,s}^{\XX}$ is finitely generated over $\XX$. 

\para 
For $(\a_{i,j},t) \in \D^+ \times \ZZ_{\geq 0}$, put 
\begin{align*}
x_{\a_{i,j}} (t) := [[\dots [[ x_{j-1,0}^+, x_{j-2,0}^+]_v, x_{j-3,0}^+]_v, \dots, x_{i+1,0}^+]_v, x_{i,t}^+]_v 
\end{align*}
as an element of $\cA^{\XX}$. 
For $h \in H_{\geq 0}$, put 
\begin{align*}
x_h := \prod_{(\b,t) \in \D^+ \times \ZZ_{\geq 0}}^{\ra} x_{\b} (t)^{h(\b,t)}. 
\end{align*}
We also set 
\begin{align*}
&\cB^{\XX} :=\{ x_h \mid h \in H_{\geq 0}\}, 
\\
&\cB_{\g,s}^{\XX} :=\{ x_h \in \cB^{\XX} \mid \sum_{(\b,t) \in \D^+ \times \ZZ_{\geq 0}} h(\b,t) \cdot \b=\g, 
	\sum_{(\b,t) \in \D^+ \times \ZZ_{\geq 0}} h (\b,t) \cdot t =s \}. 
\end{align*}

\para 
Let $U_v (L\Fsl_n)$ be the quantum loop algebra over $\KK$ associated with $\Fsl_n$. 
Then we have an algebra homomorphism 
$\theta : \cA^{\KK} \ra U_v(L \Fsl_n)$ by $\theta(x_{i,t}) =e_{i,t}$. 
By \cite[Theorem 2.17]{T}, 
we see that the set $\{\theta (x_h) \mid h \in H_{\geq 0}\}$ is linearly independent, 
and we have 
\begin{align}
\label{dim cA B}
\dim_{\KK} \cA_{\g,s}^{\KK} \geq \sharp \cB_{\g,s}^{\KK} 
\end{align}
for $(\g,s) \in Q^+ \times \ZZ_{\geq 0}$. 

\para 
We note that the scalar extension $\CC \otimes_{\CC[v]_{(v=1)}} \cA^{\CC[v]_{(v=1)}}$ through the ring homomorphism 
$\CC[v]_{(v=1)} \ra \CC$ ($v \mapsto 1$) 
is isomorphic to the universal enveloping  algebra of the positive part of the polynomial current Lie algebra $\Fsl_n[x]$. 
Then, by the same argument using \eqref{dim cA B} as one of \cite[the proof of Proposition 1.13]{L}, 
we see that 
$\cB^{\KK}_{\g,s}$ gives  a 
$\KK$-basis of $\cA_{\g,s}^{\KK}$ for each $(\g,s) \in Q^+ \times \ZZ_{\geq 0}$.  
As a consequence, we have the following lemma. 
\begin{lem}
\label{Lemma basis cAKK} 
The set $\cB^{\KK}$ gives a $\KK$-basis of $\cA^{\KK}$, 
and the algebra homomorphism $\theta : \cA^{\KK} \ra U_v (L \Fsl_n) $ is injective. 
\end{lem}

\para 
For $(\b,t)\in \D^+ \times \ZZ_{\geq 0}$, 
put $\wt{x}_{\b} (t) := (v-v^{-1}) x_{\b} (t) \in \cA^{\KK}$, 
and we set $\wt{x}_h := \prod_{(\b,t) \in \D^+ \times \ZZ_{\geq 0}}^{\ra} \wt{x}_\b(t)^{h(\b,t)}$ for $h \in H_{\geq 0}$. 
Let $\wt{\cA}^{\AA}$ be the $\AA$-subalgebra of $\cA^{\KK}$ generated by 
$\{\wt{x}_\b (t) \mid (\b,t) \in \D^+ \times \ZZ_{\geq 0} \}$. 
By definitions, we have $\{\wt{x}_h \mid h \in H_{\geq 0}\} \subset \wt{\cA}^{\AA}$, 
and we see that the set $\{\wt{x}_h \mid h \in H_{\geq 0}\} \subset \wt{\cA}^{\AA}$ is linearly independent over $\AA$ 
thanks to Lemma \ref{Lemma basis cAKK}. 
On the other hand, 
for any $X \in \wt{\cA}^{\AA}$, 
we can write $X = \sum_{h \in H_{\geq 0}} c_h \wt{x}_h$ ($c_h \in \KK$) uniquely 
since $\wt{\cA}^{\AA} \subset \cA^{\KK}$ and $\{\wt{x}_h \mid h \in H_{\geq 0}\}$ is a $\KK$-basis of $\cA^{\KK}$ 
by Lemma  \ref{Lemma basis cAKK}. 
Then we have $\theta (X) = \sum_{h \in H_{\geq 0}} c_h \theta (\wt{x}_h)$, 
and we see that $c_h \in \AA$ by \cite[Theorem 2.19 (b)]{T}, 
where we note that the element $\theta (\wt{x}_h)$ coincides with the element $\wt{e}_h$ in \cite[Theorem 2.19 (b)]{T} by definitions.  
As a consequence, we have the following lemma. 

\begin{lem}
The set $\{ \wt{x}_h \mid h \in H_{\geq 0}\}$ gives a free $\AA$-basis of $\wt{\cA}^{\AA}$. 
\end{lem}

\para 
Recall that $\cA = \CC \otimes_{\AA} \cA^{\AA}$ through the ring homomorphism $\cA \ra \CC$ ($v \mapsto q$). 
Note that $\wt{\cA}^{\AA} \subset \cA^{\AA}$, this embedding induces the algebra homomorphism 
$\Phi : \CC \otimes_{\AA} \wt{\cA}^{\AA} \ra \CC \otimes_{\AA} \cA^{\AA} = \cA$. 
On the other hand, 
we can check that there exists an algebra homomorphism 
$\Psi : \cA \ra \CC \otimes_{\AA} \wt{\cA}^{\AA}$ such that $\Psi(x_{i,t}) = (q-q^{-1})^{-1} \otimes \wt{x}_{\a_{i,i+1}}(t)$ 
if $q \not= \pm 1$. 
Note that $\Psi (x_{\b}(t)) = (q-q^{-1})^{-1} \otimes \wt{x}_{\b}(t)$ for any $(\b,t) \in \D^+ \times \ZZ_{\geq 0}$, 
the homomorphism $\Psi $ is surjective. 
>From definitions, we have $\Phi \circ \Psi (x_{i,t}) = x_{i,t}$ for $(i,t) \in I \times \ZZ_{\geq 0}$. 
Thus, the homomorphism $\Phi \circ \Psi$ is the identity. 
In particular, $\Psi$ is an isomorphism. 
As a consequence, we have the following proposition. 


\begin{prop}
\label{Prop basis cA}
Assume  that 
$q \not= \pm 1$, the set $\{x_h \mid h \in H_{\geq 0}\}$ gives a $\CC$-basis of $\cA$. 
\end{prop}

\para 
Let $\cA^0$ be an associative algebra over $\CC$ generated by $\{ J_{i,t} , \, K_i^{\pm} \mid i \in I, \, t \in \ZZ_{\geq 0}\}$ 
subject to the defining relations (Q1-1) and (Q1-2). 
Then we have the surjective algebra homomorphism 
\begin{align*}
\pi^0 : \cA^0 \ra U_{q,\bQ}^0 \,\, ( J_{i,t} \mapsto J_{i,t}, \, K_i^{\pm} \mapsto K_i^{\pm}). 
\end{align*}
By definition, we see easily that $\{ K^{\bk} J_{h_0} \mid \bk \in \ZZ^{n-1}, h_0 \in H_0\}$ gives a $\CC$-basis of $\cA^0$, 
where we use the same notation with one in \S \ref{section qQCA}. 

By \eqref{weak tri decom}, we have the surjective linear map 
\begin{align*}
\pi : \cA^{\opp} \otimes_{\CC} \cA^0 \otimes_{\CC} \cA  
\xrightarrow{\pi^- \otimes \pi^0 \otimes \pi^+}
U_{q,\bQ}^- \otimes_{\CC} U_{q, \bQ}^0 \otimes_{\CC} U_{q,\bQ}^+ 
\xrightarrow{\text{multiplication}} 
U_q^{\lan \bQ \ran}. 
\end{align*}
Moreover, we see that the set 
\begin{align*}
\{ \Theta^{\lan \bQ \ran} \circ \pi (x_h \otimes K^\bk J_{h_0} \otimes x_{h'} )
\mid h, h' \in H_{\geq 0}, \, \bk \in \ZZ^{n-1}, \, h_0 \in H_0\}
\end{align*}
is linearly independent by \cite[Proposition 5.1]{FT} and \cite[Theorem 2.15]{T}. 
Thus, the set 
$\{ \pi (x_h \otimes K^\bk J_{h_0} \otimes x_{h'} )
\mid h, h' \in H_{\geq 0}, \, \bk \in \ZZ^{n-1}, \, h_0 \in H_0\}$ is linearly independent. 
Combining with Proposition \ref{Prop basis cA}, 
we see that $\pi$ (resp. $\pi^\pm$, $\pi^0$) is an isomorphism, 
and we obtain Theorem \ref{Theorem PBW}.


\section{A proof of Lemma \ref{Lemma X1+k X0-k+1}}
\label{proof rel X1+k X0-k+1}

In this appendix, we give a proof of Lemma \ref{Lemma X1+k X0-k+1}, 
so we consider some relations of $U_q (\Fsl_2^{\lan 0 \ran}[x]) $ in the case where rank one and $Q=0$. 

\para 
By the induction on $k \geq 0$, we can show that 
\begin{align}
\label{Xt+1+ Xt+k}
&X_{t+1}^+ X_t^{+(k)} = q^k \frac{1}{[2]} (J_1 X_t^{+(k+1)} - X_t^{+(k+1)} J_1), 
\\
\label{Xt-k Xt+1-} 
&X_t^{-(k)} X_{t+1}^-  = - q^k \frac{1}{[2]} (J_1 X_t^{-(k+1)} - X_t^{-(k+1)} J_1) 
\end{align}
for $t \geq 0$, 
where we note that $X_t^{+(k)} X_{t+1}^+ = q^{-2k} X_{t+1}^+ X_t^{+(k)}$ by the relation (Q2). 
We also remark that 
\eqref{Xt-k Xt+1-} follows from \eqref{Xt+1+ Xt+k} 
by applying the algebra anti-involution  
$\dag$ given in Lemma \ref{Lemma dag}. 
The relations \eqref{Xt+1+ Xt+k} and \eqref{Xt-k Xt+1-}  also hold in the case of $k=-1$ 
if we put $X_t^{\pm (-1)}=0$. 


\para 
By the induction on $k >0$, we can show that 
\begin{align}
\label{X1+ X0-k Q=0}
& X_1^+ X_0^{-(k)} = X_0^{-(k)} X_1^+ + q^{- k+1} X_0^{-(k-1)} K^+ J_1 - q^{- 2 (k-1)} X_0^{-(k-2)} X_1^- K^+, 
\\
\label{X1+k X0- Q=0} 
& X_1^{+(k)} X_0^- = X_0^- X_1^{+(k)} + q^{-k+1} K^+ J_1 X_1^{+(k-1)} - q^{- 2 (k-1)} K^+ X_2^+ X_1^{+(k-2)}, 
\end{align}
where we put $X_0^{-(-1)}= X_1^{+(-1)} =0$. 

For $k>0$, 
applying \eqref{X1+ X0-k Q=0} and \eqref{X1+k X0- Q=0} to the right-hand side of the equation 
\begin{align*}
X_1^{+(k)} X_0^{-(k+1)} 
&=q^{-k} ([k+1] - q^{-1} [k]) X_1^{+(k)} X_0^{-(k+1)} 
\\
&=q^{-k} X_1^{+(k)} X_0^- X_0^{-(k)} - q^{-k-1} X_1^{+(k-1)} X_1^+ X_0^{-(k+1)}, 
\end{align*}
we have 
\begin{align*}
& X_1^{+(k)} X_0^{-(k+1)} 
\\
&= q^{-k} \big\{ X_0^- X_1^{+(k)} + q^{-k+1} K^+ J_1 X_1^{+(k-1)} - q^{- 2 (k-1)} K^+ X_2^+ X_1^{+(k-2)} \big\} X_0^{-(k)} 
	\\ & \quad 
	- q^{-k-1} X_1^{+(k-1)} \big\{ X_0^{-(k+1)} X_1^+ + q^{- k} X_0^{-(k)} K^+ J_1 - q^{- 2 k} X_0^{-(k-1)} X_1^- K^+ \big\}. 
\end{align*}
Applying \eqref{Xt+1+ Xt+k} and \eqref{Xt-k Xt+1-} to the right-hand side of this equation, 
we have 
\begin{align}
\label{X1+k X0-k+1}
\begin{split}
& X_1^{+(k)} X_0^{-(k+1)} 
\\
&= q^{-k} \frac{1}{[k]} X_0^- X_1^+ X_1^{+(k-1)} X_0^{-(k)} 
	+ q^{-2 k} \frac{1}{[2]}  
	\big( J_1 X_1^{+(k-1)} X_0^{-(k)}  -  X_1^{+(k-1)}  X_0^{-(k)} J_1  \big) K^+ 
	\\ & \quad 
	- q^{-k-1} X_1^{+(k-1)}  X_0^{-(k+1)} X_1^+. 
\end{split}
\end{align}

\para 
We prove Lemma \ref{Lemma X1+k X0-k+1} by the induction on $k$. 
If $k=1$, the statement follows from \eqref{X1+ X0-k Q=0}. 
If $k>1$, 
we have 
\begin{align*}
X_1^{+(k)} X_0^{-(k+1)} 
&\equiv q^{-k} \frac{1}{[k]} X_0^- X_1^+ \big\{ q^{-(k-1)k} \sum_{z=0}^{k-1} (-1)^z X_z^- (K^+)^{k-1} J_{[k-1-z]}^{\lan 0 \ran} \big\} 
	\\ & \quad 
	+ q^{-2k} \frac{1}{[2]} J_1 \big\{ q^{-(k-1)k} \sum_{z=0}^{k-1} (-1)^z X_z^- (K^+)^{k-1} J_{[k-1-z]}^{\lan 0 \ran} \big\}  K^+ 
	\\ & \quad 
	- q^{-2k} \frac{1}{[2]} \big\{ q^{-(k-1)k} \sum_{z=0}^{k-1} (-1)^z X_z^- (K^+)^{k-1} J_{[k-1-z]}^{\lan 0 \ran} \big\}  J_1 K^+ 
	\mod \fX_+ 
\end{align*}
by applying the induction hypothesis to the right-hand side of the equation \eqref{X1+k X0-k+1}.
This equation together with the relations (Q6) and \eqref{Ji1 Xitpm} implies that 
\begin{align*}
X_1^{+(k)} X_0^{-(k+1)} 
&\equiv 
	q^{-k^2 }\frac{1}{[k]}  X_0^- \sum_{z=0}^{k-1} (-1)^z   (K^+)^k J_{z+1} J_{[k-1-z]}^{\lan 0 \ran} 
 	\\ & \quad 
	+ q^{-k^2 -k}  \frac{1}{[2]} \sum_{z=0}^{k-1} (-1)^z ( X_z^- J_1 -[2] X_{z+1}^- ) (K^+)^{k} J_{[k-1-z]}^{\lan 0 \ran}
	\\ & \quad 
	- q^{-k^2 -k}  \frac{1}{[2]}  \sum_{z=0}^{k-1} (-1)^z X_z^- (K^+)^{k} J_1 J_{[k-1-z]}^{\lan 0 \ran} 
	\mod \fX_+ 
\\
&=  q^{-k (k+1)}  X_0^- (K^+)^k  q^k \frac{1}{[k]} \sum_{z=0}^{k-1} (-1)^{(z+1) -1}   J_{z+1} J_{[k-(z+1)]}^{\lan 0 \ran}  
	\\ & \quad 
	+ q^{-k(k+1)}   \sum_{z=0}^{k-1} (-1)^{z+1}  X_{z+1}^-  (K^+)^{k} J_{[k-(z+1)]}^{\lan 0 \ran}. 
\end{align*}
Note the definition \eqref{def J0[t]}, 
and  this equation implies the statement of Lemma \ref{Lemma X1+k X0-k+1}.


\section{A proof of Lemma \ref{Lemma X0+k X0-k+1} }
\label{proof rel X0+k X0-k+1}

\begin{lem}
\label{Lemma JQkt} 
For $k \in \ZZ_{>0}$ and $t=1,2, \dots, k$, we have 
\begin{align}
\label{JQkt} 
&J^{\lan Q \ran}_{[k;t]} = 
	(K^-)^2 J^{\lan Q \ran}_{[k-1;t-1]} + q^{-2k} Q J^{\lan Q \ran}_{[k-1;t]} 
	\quad  \text{ if } t <k, 
\end{align}
\begin{align}
\label{JQkk}
&J^{\lan Q \ran}_{[k;k]} 
	=q^{-k} \frac{1}{[k]} \big\{ (Q J_1 - q^{2k} J_0 + q^k[k] ) J^{\lan Q \ran}_{[k-1;k-1]} 
	+ \sum_{z=1}^{k-1} (-1)^{z-1} (J_z - Q J_{z+1}) J^{\lan Q \ran}_{[k-1;k-z-1]} \big\}. 
\end{align}

\begin{proof}
We prove \eqref{JQkt} by the induction on $t$. 
In the case where $t=1$, we can check \eqref{JQkt}  by direct calculations using definitions. 
Suppose that $t>1$. 
Applying the induction hypothesis to the right-hand side of the definition \eqref{def JktQ}, we have 
\begin{align}
\label{JktQ-1}
\begin{split} 
J^{\lan Q \ran}_{[k;t]} 
&= q^t \frac{1}{[t]} \sum_{z=1}^{t-1} (-1)^{z-1} \big( J_z - q^{2(k-t+z)} Q^{-z} J_0 + q^{k- 2 (t-z)} [k] Q^{-z} \big) 
	\\ & \hspace{5em} \times  
	\big\{ (K^-)^2 J^{\lan Q \ran}_{[k-1;t-z-1]}  + q^{- 2k} Q J^{\lan Q \ran}_{[k-1; t-z]} \big\} 
	\\ & \quad 
	+ q^t \frac{1}{[t]} (-1)^{t-1} \big( J_t - q^{2k} Q^{-t} J_0 + q^k [k] Q^{-t} \big) J^{\lan Q \ran}_{[k;0]}
\\
&= (K^-)^2  q^t \frac{1}{[t]} \sum_{z=1}^{t-1} (-1)^{z-1} \big( J_z - q^{2 ( k  - t  +z)} Q^{-z} J_0 
		+ q^{k- 2 (t-z)}[k] Q^{-z} \big) J^{\lan Q \ran}_{[k-1; t-z-1]}
	\\ & \quad 
	+ q^{-2k} Q q^t \frac{1}{[t]} \sum_{z=1}^{t} (-1)^{z-1} \big( J_z - q^{2(k-t+z)} Q^{-z} J_0 + q^{k- 2 (t-z)} [k] Q^{-z} \big) 
		 J^{\lan Q \ran}_{[k-1;t-z]}, 
\end{split}
\end{align}
where we note that $ J^{\lan Q \ran}_{[k;0]} = q^{-2k} Q J^{\lan Q \ran}_{[k-1;0]}$ by definition. 
Note that  
\begin{align*}
&- q^{2(k-t+z)} Q^{-z} J_0 + q^{k- 2 (t-z)} [k] Q^{-z} 
\\
&= - \big\{ q^{2(k-1-t+z)} + (q-q^{-1}) q^{2(k-t+z)-1} \big\} Q^{-z} J_0  
	+ \big\{ q^{k-1-2 (t-z)}[k-1] + q^{ 2 (k-t+z)-1} \big\} Q^{-z} 
\\
&= - q^{2(k-1-t+z)} Q^{-z} J_0 + q^{k-1-2 (t-z)}[k-1] Q^{-z} + q^{2(k-t+z)-1} (K^-)^2 Q^{-z}
\end{align*}
by the relation (Q1-2), 
and the equation \eqref{JktQ-1} implies that  
\begin{align} 
\label{JktQ-2} 
\begin{split} 
J^{\lan Q \ran}_{[k;t]} 
&= (K^-)^2  q^t \frac{1}{[t]} \sum_{z=1}^{t-1} (-1)^{z-1} \big( J_z - q^{2 ( k  - t  +z)} Q^{-z} J_0 
		+ q^{k- 2 (t-z)}[k] Q^{-z} \big) J^{\lan Q \ran}_{[k-1; t-z-1]}
	\\ & \quad 
	+ q^{-2k} Q q^t \frac{1}{[t]} \sum_{z=1}^{t} (-1)^{z-1} \big( J_z - q^{2(k-1-t+z)} Q^{-z} J_0 + q^{k-1-2 (t-z)}[k-1] Q^{-z} \big) 
		 J^{\lan Q \ran}_{[k-1;t-z]} 
	\\ & \quad 
	+ q^{-2k} Q q^t \frac{1}{[t]} \sum_{z=1}^t (-1)^{z-1} q^{2(k-t+z)-1} (K^-)^2 Q^{-z} J^{\lan Q \ran}_{[k-1;t-z]}. 
\end{split} 
\end{align}
Note that 
\begin{align*}
&q^{-2k} Q q^t \frac{1}{[t]} \sum_{z=1}^t (-1)^{z-1} q^{2(k-t+z)-1} (K^-)^2 Q^{-z} J^{\lan Q \ran}_{[k-1;t-z]} 
\\
&= (K^-)^2 q^t \frac{1}{[t]} \sum_{z=1}^t  (-1)^{(z-1)-1} (- q^{- 2(t- (z-1)) +1} ) Q^{-(z-1)} J^{\lan Q \ran}_{[k-1; t-(z-1)-1]}, 
\end{align*} 
and 
$ q^{k- 2 (t-z)}[k] - q^{-2(t-z)+1} = q^{(k-1) - 2 (t-1-z)}[k-1]$, 
then the equation \eqref{JktQ-2} implies that 
\begin{align*}
J^{\lan Q \ran}_{[k;t]} 
&= (K^-)^2  q^t \frac{1}{[t]} \sum_{z=1}^{t-1} (-1)^{z-1} \big( J_z - q^{2 ( (k-1)  - (t-1)  +z)} Q^{-z} J_0 
		+ q^{(k-1) - 2 (t-1-z)}[k-1] Q^{-z} \big) J^{\lan Q \ran}_{[k-1; t-z-1]}
	\\ & \quad 
	+ q^{-2k} Q q^t \frac{1}{[t]} \sum_{z=1}^{t} (-1)^{z-1} \big( J_z - q^{2(k-1-t+z)} Q^{-z} J_0 + q^{k-1-2 (t-z)}[k-1] Q^{-z} \big) 
		 J^{\lan Q \ran}_{[k-1;t-z]} 
	\\ & \quad 
	+ (K^-)^2 q^t \frac{1}{[t]}  q^{- 2 t +1}  J^{\lan Q \ran}_{[k-1; t-1]}. 
\end{align*}
Applying the definition \eqref{def JktQ}, we have 
\begin{align*} 
J^{\lan Q \ran}_{[k;t]} 
&= (K^-)^2 q^t \frac{1}{[t]} q^{-t+1} [t-1] J^{\lan Q \ran}_{[k-1;t-1]}  
	+ q^{- 2 k} Q J^{\lan Q \ran}_{[k-1;t]}  
	 + (K^-)^2 q^t \frac{1}{[t]}  q^{- 2 t +1}  J^{\lan Q \ran}_{[k-1; t-1]}
\\
&= (K^-)^2  J^{\lan Q \ran}_{[k-1;t-1]}  + q^{- 2 k} Q J^{\lan Q \ran}_{[k-1;t]}. 
\end{align*}

Next we prove \eqref{JQkk}. 
Applying \eqref{JQkt} to the right-hand side of \eqref{def JktQ}, we have 
\begin{align*}
J_{[k;k]}^{\lan Q \ran} 
&=q^k \frac{1}{[k]} \sum_{z=1}^{k-1} (-1)^{z-1} 
	\big( J_z - q^{2  z} Q^{-z} J_0  + q^{ - k+ 2z} [k] Q^{-z} \big) 
	\\ & \hspace{5em} \times 
	\big\{ (K^-)^2 J^{\lan Q \ran}_{[k-1; k-z-1]} + q^{-2 k} Q J^{\lan Q \ran}_{[k-1; k-z]} \big\} 
	\\ & \quad 
	+ q^k \frac{1}{[k]} (-1)^{k-1} \big( J_k - q^{2k} Q^{-k} J_0 + q^k [k] Q^{-k} \big) J^{\lan Q \ran}_{[k;0]}. 
\end{align*}
This implies that 
\begin{align*}
J_{[k;k]}^{\lan Q \ran}
&= q^{-k} \frac{1}{[k]} \sum_{z=1}^{k-2} (-1)^{z-1} \big\{ 
		\big( q^{2k} J_z - q^{ 2k + 2 z} Q^{-z} J_0  + q^{  k+ 2z} [k] Q^{-z} \big)   (K^-)^2 
		\\ & \hspace{3em} 
		- \big( Q J_{z+1} - q^{2  z +2} Q^{-z} J_0  + q^{ - k+ 2z +2} [k] Q^{-z } \big) \big\}  J^{\lan Q \ran}_{[k-1;k-z-1]}
	\\ & \quad 
	+ q^{ -k} \frac{1}{[k]} (-1)^{k-1}\big\{ 
		 \big( Q J_k  - q^{2k} Q^{-k+1} J_0 + q^k [k] Q^{-k+1} \big) 
		 \\ & \hspace{3em} 
		 - \big(q^{2k} J_{k-1} - q^{4k -2} Q^{-k+1} J_0 + q^{3k-2}[k] Q^{-k+1} \big) (K^-)^2
	 	\big\} J^{\lan Q \ran}_{[k-1;0]} 
	\\ & \quad 
	+ q^{-k} \frac{1}{[k]} \big( Q J_1 - q^{2} J_0 + q^{-k+2}[k] \big) J^{\lan Q \ran}_{[k-1;k-1]} 
\\
&= q^{-k} \frac{1}{[k]} \big\{ (Q J_1 - q^{2k} J_0 + q^k[k] ) J^{\lan Q \ran}_{[k-1;k-1]} 
	+ \sum_{z=1}^{k-1} (-1)^{z-1} (J_z - Q J_{z+1}) J^{\lan Q \ran}_{[k-1;k-z-1]} \big\}
	\\ & \quad 
	+ q^{-k} \frac{1}{[k]} \sum_{z=1}^{k-2} (-1)^{z-1} \big\{ 
		\big( q^{2k} J_z - q^{ 2k + 2 z} Q^{-z} J_0  + q^{  k+ 2z} [k] Q^{-z} \big)   (K^-)^2 
		\\ & \hspace{3em} 
		- \big( Q J_{z+1} - q^{2  z +2} Q^{-z} J_0  + q^{ - k+ 2z +2} [k] Q^{-z } \big) 
		- \big( J_z - Q J_{z+1} \big) \big\}  J^{\lan Q \ran}_{[k-1;k-z-1]}
	\\ & \quad 
	+ q^{ -k} \frac{1}{[k]} (-1)^{k-1}\big\{ 
		 \big( Q J_k  - q^{2k} Q^{-k+1} J_0 + q^k [k] Q^{-k+1} \big) 
		 \\ & \hspace{3em} 
		 - \big(q^{2k} J_{k-1} - q^{4k -2} Q^{-k+1} J_0 + q^{3k-2}[k] Q^{-k+1} \big) (K^-)^2
		 + \big( J_{k-1} - Q J_k \big) 
	 	\big\} J^{\lan Q \ran}_{[k-1;0]} 
	\\ & \quad 
	+ q^{-k} \frac{1}{[k]} \big\{ \big( Q J_1 - q^{2} J_0 + q^{-k+2}[k] \big) 
		- \big( Q J_1 - q^{2k}J_0 +q^k [k] \big) 
	 \big\} J^{\lan Q \ran}_{[k-1;k-1]}
\\
&= q^{-k} \frac{1}{[k]} \big\{ (Q J_1 - q^{2k} J_0 + q^k[k] ) J^{\lan Q \ran}_{[k-1;k-1]} 
	+ \sum_{z=1}^{k-1} (-1)^{z-1} (J_z - Q J_{z+1}) J^{\lan Q \ran}_{[k-1;k-z-1]} \big\}
	\\ & \quad 
	+ q^{-k} \frac{1}{[k]} \sum_{z=1}^{k-1} (-1)^{z-1} \big\{ 
		q^{2k} \big( J_z - q^{  2 z} Q^{-z} J_0  + q^{  - k+ 2z} [k] Q^{-z} \big)   (K^-)^2 
		\\ & \hspace{3em} 
		- \big( J_{z} - q^{2  z +2} Q^{-z} J_0  + q^{ - k+ 2z +2} [k] Q^{-z } \big)  
		\big\}  J^{\lan Q \ran}_{[k-1;k-z-1]}
	\\ & \quad 
	- q^{-k} \frac{1}{[k]} \big\{  q^{k+1} (K^-)^2 - q^{-k+1}\big\} [k-1] J^{\lan Q \ran}_{[k-1;k-1]}, 
\end{align*}
where we use the relation (Q1-2) in the last term. 
Applying the definition \eqref{def JktQ} to the last term of the above equation, we have 
\begin{align*}
J_{[k;k]}^{\lan Q \ran}
&= q^{-k} \frac{1}{[k]} \big\{ (Q J_1 - q^{2k} J_0 + q^k[k] ) J^{\lan Q \ran}_{[k-1;k-1]} 
	+ \sum_{z=1}^{k-1} (-1)^{z-1} (J_z - Q J_{z+1}) J^{\lan Q \ran}_{[k-1;k-z-1]} \big\}
	\\ & \quad 
	+ q^{-k} \frac{1}{[k]} \sum_{z=1}^{k-1} (-1)^{z-1} \big\{ 
		q^{k + 2 z} \big(  [k]  - q [k-1]\big)  (K^-)^2 
		\\ & \hspace{3em} 
		- \big( - q^{2  z +1} (q-q^{-1})  J_0  + q^{ - k+ 2z +2} ([k] - q^{-1} [k-1]) \big)  
		\big\}  Q^{-z} J^{\lan Q \ran}_{[k-1;k-z-1]}
\\
&= q^{-k} \frac{1}{[k]} \big\{ (Q J_1 - q^{2k} J_0 + q^k[k] ) J^{\lan Q \ran}_{[k-1;k-1]} 
	+ \sum_{z=1}^{k-1} (-1)^{z-1} (J_z - Q J_{z+1}) J^{\lan Q \ran}_{[k-1;k-z-1]} \big\}, 
\end{align*}
where we also use the relation  (Q1-2). 
\end{proof}
\end{lem}


\para 
By the induction on $k>0$, we can show that 
\begin{align} 
\label{X0+ X0-k}
\begin{split}
&X_0^+ X_0^{-(k)} 
\\
&= X_0^{-(k)} X_0^+ + X_0^{-(k-1)} K^+ (q^{k-1} J_0 - q^{-k+1}Q J_1) 
	- (X_0^- - q^{-2} Q X_1^-) X_0^{-(k-2)} K^+,
\end{split}
\\
\label{X0+k X0-}
\begin{split}
&X_0^{+(k)} X_0^- 
\\
&= X_0^- X_0^{+(k)} + K^+ (q^{k-1} J_0 - q^{-k+1} Q J_1) X_0^{+(k-1)} 
	- K^+ X_0^{+(k-2)} (X_0^+ - q^{-2} Q X_1^+), 
\end{split}
\end{align}
where we put $X_0^{+(-1)} = X_0^{-(-1)} =0$.  
We remark that the relation \eqref{X0+k X0-} follows from \eqref{X0+ X0-k} 
by applying the algebra anti-involution $\dag$ given in Lemma \ref{Lemma dag}. 

\para 
For $k >0$, applying \eqref{X0+ X0-k} and \eqref{X0+k X0-} to the right-hand side of the equation 
\begin{align*}
X_0^{+(k)} X_0^{-(k+1)} 
&= q^{-k} ( [k+1] - q^{-1} [k]) X_0^{+(k)} X_0^{-(k+1)} 
\\
&= q^{-k} X_0^{+(k)} X_0^- X_0^{-(k)}  - q^{-k-1} X_0^{+(k-1)} X_0^+ X_0^{-(k+1)}, 
\end{align*}
we have 
\begin{align*}
X_0^{+(k)} X_0^{-(k+1)}
&=q^{-k} \big\{ X_0^- X_0^{+(k)} + K^+ (q^{k-1} J_0 - q^{-k+1} Q J_1) X_0^{+(k-1)} 
		\\ & \hspace{4em} 
				- K^+ X_0^{+(k-2)} (X_0^+ - q^{-2} Q X_1^+) \big\} X_0^{-(k)} 
	\\ & \quad 
	- q^{-k-1} X_0^{+(k-1)} 
	\big\{ X_0^{-(k+1)} X_0^+ + X_0^{-(k)} K^+ (q^k J_0 - q^{-k} Q J_1) 
		\\ & \hspace{9em} 
		- (X_0^- - q^{-2} Q X_1^-) X_0^{-(k-1)} K^+ \big\}. 
\end{align*}
We note that the relations \eqref{Xt+1+ Xt+k} and \eqref{Xt-k Xt+1-} also hold in the case where $Q \not=0$. 
Applying the relations (Q2), (Q3), \eqref{Xt+1+ Xt+k} and \eqref{Xt-k Xt+1-} to the right-hand side of the above equation, 
we have 
\begin{align}
\label{X0+k X0-k+1}
\begin{split}
&X_0^{+(k)} X_0^{-(k+1)} 
\\
&=  q^{-k} \frac{1}{[k]}  X_0^- X_0^+ X_0^{+(k-1)} X_0^{-(k)} 
	+ \big( q^{-3} J_0 - q^{-2k} \frac{1}{[2]} Q   J_1 \big) X_0^{+(k-1)} X_0^{-(k)} K^+ 
	\\ & \quad 
	-  X_0^{+(k-1)}  X_0^{-(k)} \big( q^{-1} J_0  -  q^{-2k} \frac{1}{[2]} Q J_1 -  q^{-2} \big) K^+
	- q^{-k-1} X_0^{+(k-1)} X_0^{-(k+1)} X_0^+. 
\end{split}
\end{align}

\para 
We prove Lemma \ref{Lemma X0+k X0-k+1} by the induction on $k$. 
In the case where $k=1$, 
the statement follows from \eqref{X0+ X0-k}. 
Suppose that $k >1$. 
By \eqref{X0+k X0-k+1}, we have 
\begin{align*}
X_0^{+(k)} X_0^{-(k+1)} 
&\equiv 
q^{-k} \frac{1}{[k]}  X_0^- X_0^+ X_0^{+(k-1)} X_0^{-(k)} 
	+ \big( q^{-3} J_0 - q^{-2k} \frac{1}{[2]} Q   J_1 \big) X_0^{+(k-1)} X_0^{-(k)} K^+ 
	\\ & \quad 
	-  X_0^{+(k-1)}  X_0^{-(k)} \big( q^{-1} J_0  -  q^{-2k} \frac{1}{[2]} Q J_1 -  q^{-2} \big) K^+
	\quad \mod \fX_+. 
\end{align*}
Applying the induction hypothesis, 
we have 
\begin{align*}
&X_0^{+(k)} X_0^{-(k+1)} 
\\
&\equiv q^{-k} \frac{1}{[k]} X_0^- X_0^+ 
		\big\{ \sum_{z=0}^{k-1} (-1)^{k-z-1} X_z^- (K^+)^{k-1} J^{\lan Q \ran}_{[k-1; k-z-1]} \big\} 
		\\
		& \quad 
		+\big( q^{-3} J_0 - q^{- 2k}\frac{1}{[2]} Q J_1 \big) 
			\big\{ \sum_{z=0}^{k-1} (-1)^{k-z-1} X_z^- (K^+)^{k-1} J^{\lan Q \ran}_{[k-1; k-z-1]} \big\}  K^+ 
		\\
		&\quad 
		- \big\{ \sum_{z=0}^{k-1} (-1)^{k-z-1} X_z^- (K^+)^{k-1} J^{\lan Q \ran}_{[k-1; k-z-1]} \big\}  
			\big( q^{-1} J_0  -  q^{-2k} \frac{1}{[2]} Q J_1 -  q^{-2} \big) K^+ 
		\mod \fX^+
\\
&\equiv 
	q^{-k} \frac{1}{[k]} X_0^- \sum_{z=0}^{k-1} (-1)^{k-z-1} \big( J_z - Q J_{z+1} \big) 
		(K^+)^k J^{\lan Q \ran}_{[k-1; k-z-1]}
	\\
	&\quad 
	+ \sum_{z=0}^{k-1} (-1)^{k-z-1} 
		\big\{ q^{-3}  ( q^4 X_z^- J_0 - q^2 [2] X_z^- ) - q^{-2 k} \frac{1}{[2]} Q ( X_z^- J_1 - [2] X_{z+1}^-) \big\} 
		(K^+)^k J^{\lan Q \ran}_{[k-1;k-z-1]} 
	\\
	& \quad 
	- \sum_{z=0}^{k-1} (-1)^{k-z-1} X_z^- (K^+)^{k} J^{\lan Q \ran}_{[k-1; k-z-1]} 
			\big( q^{-1} J_0  -  q^{-2k} \frac{1}{[2]} Q J_1 -  q^{-2} \big)
	\mod \fX^+
\end{align*}
where we use the relations (Q1-1), (Q1-2), (Q6), \eqref{Ji1 Xitpm} 
and the fact $X_s^+ U_{q,Q}^0 \subset \fX_+$ for all $s \geq 0$ which follows from defining relations immediately. 
This equation implies 
\begin{align*}
&X_0^{+(k)} X_0^{-(k+1)} 
\\
&\equiv 
	(-1)^k X_0^- (K^+)^k \big\{ 
		q^{-k} \frac{1}{[k]} \sum_{z=0}^{k-1} (-1)^{-z-1} (J_z - Q J_{z+1} ) J^{\lan Q \ran}_{[k-1;k-z-1]}
		+ \big(  1 - (q-q^{-1}) J_0 \big) J^{\lan Q \ran}_{[k-1;k-1]} \big\} 
	\\ & \quad 
	+ \sum_{z=1}^{k-1} (-1)^{k-z} X_z^- (K^+)^k \big\{ 1 - (q-q^{-1}) J_0  \big\} J^{\lan Q \ran}_{[k-1; k-z-1]} 
	\\ & \quad 
	+ \sum_{z=0}^{k-1} (-1)^{k-(z+1)} X_{z+1}^- (K^+)^k \big\{ q^{-2k} Q \big\} J_{[k-1; k- (z+1)]}^{\lan Q \ran} 
	\mod \fX_+
\\
&= (-1)^k X_0^- (K^+)^k 
		\\ & \qquad \times 
		q^{-k} \frac{1}{[k]} \big\{  \big(  Q J_1 - q^{2 k}   J_0  + q^k[k]  \big) J^{\lan Q \ran}_{[k-1;k-1]} 
			+ \sum_{z=1}^{k-1} (-1)^{-z-1} (J_z - Q J_{z+1} ) J^{\lan Q \ran}_{[k-1;k-z-1]} \big\} 
	\\ & \quad 
	+ \sum_{z=1}^{k-1} (-1)^{k-z} X_z^- (K^+)^k \big\{ (K^-)^2 J^{\lan Q \ran}_{[k-1; k-z-1]} 
		+ q^{-2k} Q J^{\lan Q \ran}_{[k-1; k-z]} \big\} 
	\\ & \quad 
	+ X_{k}^- (K^+)^k (q^{-2k} Q ) J_{[k-1; 0]}^{\lan Q \ran}, 
\end{align*}
where we note that 
$1- (q-q^{-1}) J_0 =(K^-)^2$ by (Q1-2). 
Applying Lemma \ref{Lemma JQkt}, we have 
\begin{align*}
X_0^{+(k)} X_0^{-(k+1)} 
\equiv \sum_{z=0}^k (-1)^{k-z} X_z^- (K^+)^k J_{[k;k-z]}^{\lan Q \ran} 
\mod \fX_+,
\end{align*}
where we note that $ J^{\lan Q \ran}_{[k;0]} = q^{-2k} Q J^{\lan Q \ran}_{[k-1;0]}$ by definition. 


\section{The $(q, \bQ)$-current algebra $U_q(\Fgl_n^{\lan \bQ \ran}[x])$  and cyclotomic $q$-Schur algebras}
\label{section gl}

In this appendix, 
we consider the $(q, \bQ)$-current algebra $U_q(\Fgl_n^{\lan \bQ \ran}[x])$ associated with 
the general linear  Lie algebra $\Fgl_n$. 
We show that the algebra $U_q (\Fgl_n^{\lan \bQ \ran}[x])$ with  special parameters 
is isomorphic to the algebra $\cU_{q, \wh{\bQ}}(\bn)$ introduced in \cite{W16} (see Proposition \ref{Prop iso gl g}).  
We also give some connection with cyclotomic $q$-Schur algebras according to \cite{W16}. 


\para 
Recall that 
$A=(a_{ij})_{1\leq i,j \leq n-1}$ is the Cartan matrix of type $A_{n-1}$. 
We also 
put $\wt{a}_{ii}=1$, $\wt{a}_{i+1,i}=-1$ and $\wt{a}_{ij}=0$ if $i\not=j,j+1$ 
for $1\leq i,j \leq n$. 
We define the $(q, \bQ)$-current algebra $U_q (\Fgl_n^{\lan \bQ \ran}[x])$ 
associated with the general linear Lie algebra $\Fgl_n$ as follows. 

\begin{definition}
For $q \in \CC^{\times}$ and $\bQ =(Q_1,Q_2,\dots, Q_{n-1}) \in \CC^{n-1}$, 
we define the associative algebra $U_q (\Fgl_n^{\lan Q \ran}[x])$ over $\CC$ 
by the following generators and defining relations: 

\begin{description}
\item[Generators] 
$X_{i,t}^{\pm}$, $I_{j,t}^{\pm}$, $\wt{K}^{\pm}_j$ ($1\leq i \leq n-1$, $1\leq j \leq n$, $t \in \ZZ_{\geq 0}$), 
\item[Defining relations] 
\begin{align*}
&\tag{Q'1-1} 
[\wt{K}^{+}_i, \wt{K}^{+}_j] = [\wt{K}^{+}_i, I^{\s}_{j,t}] = [I^{\s}_{i,s}, I^{\s'}_{j,t}] =0 
	\quad (\s, \s' \in \{+,-\}), 
\\
& \tag{Q'1-2}
\wt{K}_j^+ \wt{K}_j^- = 1 = \wt{K}_j^- \wt{K}_j^+, 
\quad 
(\wt{K}_j^{\pm})^2 = 1 \pm (q-q^{-1}) I^{\mp}_{j,0}, 
\\
\tag{Q'2} 
\begin{split} 
&X_{i,t+1}^+ X_{i,s}^+ - q^2 X_{i,s}^+ X_{i,t+1}^+ = q^2 X_{i,t}^+ X_{i,s+1}^+ - X_{i,s+1}^+ X_{i,t}^+, 
\\
&  X_{i,t+1}^+ X_{i+1,s}^+ - q^{-1} X_{i+1,s}^+ X_{i,t+1}^+ = X_{i,t}^+ X_{i+1,s+1}^+ - q X_{i+1,s+1}^+ X_{i,t}^+, 
\end{split}
\\
\tag{Q'3} 
\begin{split} 
&X_{i,t+1}^- X_{i,s}^- - q^{-2} X_{i,s}^- X_{i,t+1}^- = q^{-2} X_{i,t}^- X_{i,s+1}^- - X_{i,s+1}^- X_{i,t}^-, 
\\
&  X_{i+1,s}^- X_{i,t+1}^- - q^{-1} X_{i,t+1}^- X_{i+1,s}^-  = X_{i+1,s+1}^- X_{i,t}^- - q X_{i,t}^- X_{i+1,s+1}^-, 
\end{split}\\
&\tag{Q'4-1} 
\wt{K}_i^+ X_{j,t}^+ \wt{K}_i^- = q^{\wt{a}_{ij}} X_{j,t}^+, 
\\
& \tag{Q'4-2} 
q^{\pm \wt{a}_{ij}} I^{\pm}_{i,0} X_{j,t}^+ - q^{\mp \wt{a}_{ij}} X_{j,t}^+ I^{\pm}_{i,0} 
	= \wt{a}_{ij} X_{j,t}^+, 
\\
& \tag{Q'4-3} 
[I^{\pm}_{i,s+1}, X_{j,t}^+] = q^{\pm \wt{a}_{ij}} I^{\pm}_{i,s} X_{j,t+1}^+ - q^{\mp \wt{a}_{ij}} X_{j,t+1}^+ I^{\pm}_{i,s}, 
\\
& \tag{Q'5-1} 
\wt{K}_i^+ X_{j,t}^- \wt{K}_i^- = q^{- \wt{a}_{ij}} X_{j,t}^{\pm},  
\\
&\tag{Q'5-2} 
q^{\mp \wt{a}_{ij}} I_{i,0}^{\pm} X_{j,t}^- - q^{\pm \wt{a}_{ij}} X_{j,t}^- I_{i,0}^{\pm} 
	= - \wt{a}_{ij} X_{j,t}^-, 
\\
&\tag{Q'5-3} 
[I^{\pm}_{i,s+1}, X_{j,t}^-] = q^{\mp \wt{a}_{ij}} I^{\pm}_{i,s} X_{j,t+1}^- - q^{\pm \wt{a}_{ij}} X_{j,t+1}^- I^{\pm}_{i,s}, 
\\
&\tag{Q'6}
	[X_{i,t}^+, X_{j,s}^-] = \d_{i,j} K_i^+ (  J_{i,s+t} - Q_i  J_{i,s+t+1}), 
\\
&\tag{Q'7} 
	[X_{i,t}^+, X_{j,s}^+] =0  \text{ if } j \not=i,i \pm 1, 
	\\& 
	X_{i\pm 1,u}^+ (X_{i,s}^+ X_{i,t}^+ + X_{i,t}^+ X_{i,s}^+) + (X_{i,s}^+ X_{i,t}^+ + X_{i,t}^+ X_{i,s}^+) X_{i \pm 1, u}^+  
	\\ & = (q+q^{-1}) (X_{i,s}^+ X_{i \pm 1,u}^+ X_{i,t}^+ + X_{i,t}^+ X_{i \pm 1, u}^+ X_{i,s}^+), 
\\
&\tag{Q'8} 
	[X_{i,t}^-, X_{j,s}^-] =0  \text{ if } j \not=i,i \pm 1, 
	\\& 
	X_{i\pm 1,u}^- (X_{i,s}^- X_{i,t}^- + X_{i,t}^- X_{i,s}^-) + (X_{i,s}^- X_{i,t}^- + X_{i,t}^- X_{i,s}^-) X_{i \pm 1, u}^-  
	\\ & = (q+q^{-1}) (X_{i,s}^- X_{i \pm 1,u}^- X_{i,t}^- + X_{i,t}^- X_{i \pm 1, u}^- X_{i,s}^-), 
\end{align*}
where we put 
$K_i^+ = \wt{K}_i^+ \wt{K}_{i+1}^-$, $K_i^- = \wt{K}_i^- \wt{K}_{i+1}^+$, 
\begin{align*}
J_{i,t} = \begin{cases} 
	I^+_{i,0} - I^-_{i+1,0} + (q-q^{-1}) I_{i,0}^+ I_{i+1,0}^- & \text{ if } t=0, 
	\\ \dis 
	q^{-t} I^+_{i,t} - q^t I^-_{i+1,t} - (q-q^{-1}) \sum_{z=1}^{t-1} q^{-t + 2z} I^+_{i,t-z} I^-_{i+1,z} 
	& \text{ if } t>0.
	\end{cases} 
\end{align*}
\end{description}
\end{definition} 

\begin{remark}
In the case where $q=1$, 
let $\fI$ be the two-sided ideal of $U_1(\Fgl_n^{\lan \bQ \ran}[x])$ generated by 
$\{ \wt{K}_j^+ -1 , \, I_{j,t}^+ - I_{j,t}^- \mid 1\leq j \leq n , \, t \in \ZZ_{\geq 0}\}$. 
Then we see easily that 
$U_1(\Fgl_n^{\lan \bQ \ran}[x])/ \fI$ 
is isomorphic to the universal envelope algebra of the deformed current Lie algebra $\Fgl_n^{\lan \bQ \ran}[x]$ 
given in \cite[Definition 1.1]{W18}.
\end{remark}
\\

>From the defining relations, 
we can easily check the following lemma. 

\begin{lem}
\label{Lemma dag gl}
There exists the algebra anti-involution 
$\dag : U_q (\Fgl_n^{\lan \bQ \ran}[x]) \ra U_q (\Fgl_n^{\lan\bQ \ran}[x])$  
such that 
$\dag(X_{i,t}^{\pm}) = X_{i,t}^{\mp}$, 
$\dag(I_{j,t}^{\pm}) = I_{j,t}^{\pm}$ 
and $\dag(\wt{K}_j^{\pm})  = \wt{K}_j^{\pm}$ 
for $1 \leq i \leq n-1$, $1 \leq j \leq n$ and $t \in \ZZ_{\geq 0}$. 
\end{lem}

\begin{lem}
\label{Lemma some rel UqglQ}
We have the following relations in $U_q(\Fgl_n^{\lan \bQ \ran}[x])$.  

\begin{enumerate}
\item 
$K_i^+ X_{j,t}^{\pm} K_i^- = q^{\pm a_{ij}} X_{j,t}^{\pm}$.  

\item 
$(K_i^-)^2 = 1 - (q-q^{-1}) J_{i,0}$. 

\item 
$ q^{\pm a_{ij}} J_{i,0} X_{j,t}^{\pm} - q^{\mp a_{ij}} X_{j,t}^{\pm} J_{i,0} = [ \pm a_{ij}] X_{j,t}^{\pm}$. 

\item 
$[J_{i,s+1}, X_{j,t}^+] = q^{2 \wt{a}_{i,j}}J_{i,s} X_{j,t+1}^+ - q^{2 \wt{a}_{i+1,j}} X_{j,t+1}^+ J_{i,s}$.  

\item 
$[J_{i,s+1}, X_{j,t}^-] = q^{2 \wt{a}_{i+1,j}} J_{i,s} X_{j,t+1}^- - q^{2 \wt{a}_{i,j}} X_{j,t+1}^- J_{i,s}$. 
\end{enumerate}

\begin{proof}
Note that $K^{\pm}_i = \wt{K}_i^{\pm} \wt{K}_{i+1}^{\mp}$, 
then the relation (\roi) follows from the relations (Q'1-1), (Q'4-1) and (Q'5-1) immediately. 
We also have the relation (\roii) by direct calculation using the relations (Q'1-1) and (Q'1-2).  
By the relations (\roi) and (\roii), we have 
\begin{align*}
J_{i,0} X_{j,t}^{\pm} 
= \frac{ 1 - (K_i^-)^2}{q-q^{-1}} X_{j,t}^{\pm} 
=X_{j,t}^{\pm} \frac{1 - q^{\mp 2 a_{ij}} (K_i^-)^2}{q-q^{-1}} 
= \frac{1 - q^{\mp 2 a_{ij}}}{q-q^{-1}} X_{j,t}^{\pm} + q^{\mp 2 a_{ij}} X_{j,t} J_{i,0}.  
\end{align*}
This implies the relation (\roiii). 
We prove (\roiv). 
By the definition of $J_{i,s+1}$, we have 
\begin{align*}
J_{i,s+1} X_{j,t}^+ 
= \big( q^{-(s+1)} I_{i,s+1}^+ - q^{s+1} I^-_{i+1,s+1} - (q-q^{-1}) \sum_{z=1}^{s} q^{-(s+1) + 2 z} I^+_{i,s-z+1} I^-_{i+1,z} \big) 
	X_{j,t}^+. 
\end{align*}
Applying the relation (Q'4-3), we have 
\begin{align*}
&J_{i,s+1} X_{j,t}^+ 
\\
&= X_{j,t}^+ 
	\big( q^{-(s+1)} I_{i,s+1}^+ - q^{s+1} I^-_{i+1,s+1} - (q-q^{-1}) \sum_{z=1}^{s} q^{-(s+1) + 2 z} I^+_{i,s-z+1} I^-_{i+1,z} \big)
	\\ & \quad 
	+ q^{- (s+1)} \big( q^{\wt{a}_{ij}} I^+_{i,s} X_{j,t+1}^+ - q^{- \wt{a}_{ij}} X_{j,t+1}^+ I^+_{i,s} \big) 
	\\ & \quad 
	- q^{s+1} \big( q^{-\wt{a}_{i+1,j}} I^-_{i+1,s} X^+_{j,t+1} - q^{\wt{a}_{i+1,j}} X_{j,t+1}^+ I^-_{i+1,s} \big)  
	\\ & \quad 
	- (q-q^{-1}) \sum_{z=1}^s q^{-(s+1) + 2 z} 
	\big( q^{\wt{a}_{ij}} I^+_{i,s-z} X_{j,t+1}^+ I^-_{i+1,z} - q^{- \wt{a}_{ij}} X^+_{j,t+1} I^+_{i,s-z} I^-_{i+1,z} 
		\\ & \hspace{5em} 
		+ q^{- \wt{a}_{i+1,j}} I^+_{i,s-z+1} I^-_{i+1,z-1} X^+_{j,t+1} 
		- q^{\wt{a}_{i+1,j}} I^+_{i,s-z+1} X^+_{j,t+1} I^-_{i+1,z-1} \big)
\\
&= X_{j,t}^+ J_{i,s+1} 
	+ \big\{ q^{\wt{a}_{ij} -s -1} I_{i,s}^+ - q^{- \wt{a}_{i+1,j} + s+1} I^-_{i+1,s} 
		\\ & \hspace{7em} 
		- (q-q^{-1}) \sum_{z=1}^s q^{- \wt{a}_{i+1,j}-s + 2 z -1}  I^+_{i,s-z+1} I^-_{i+1,z-1}  \big\} X_{j,t+1}^+ 
	\\ & \quad 
	- X_{j,t+1}^+ \big\{ q^{-\wt{a}_{ij} - s-1} I_{i,s}^+ - q^{\wt{a}_{i+1,j} +s+1} I^-_{i+1,s} 
		- (q-q^{-1}) \sum_{z=1}^s q^{-\wt{a}_{ij} -s + 2 z -1} I^+_{i,s-z} I^-_{i+1,z} \big\} 
	\\ & \quad 
	- (q-q^{-1}) \sum_{z=1}^s q^{- (s+1) + 2 z}
	 \big( q^{\wt{a}_{ij}} I^+_{i,s-z} X_{j,t+1}^+ I^-_{i+1,z} - q^{\wt{a}_{i+1,j}} I^+_{i,s-z+1} X_{j,t+1}^+ I^-_{i+1,z-1} \big), 
\end{align*}
where we note that 
\begin{align*}
&\sum_{z=1}^s q^{- (s+1) + 2 z}
	 \big( q^{\wt{a}_{ij}} I^+_{i,s-z} X_{j,t+1}^+ I^-_{i+1,z} - q^{\wt{a}_{i+1,j}} I^+_{i,s-z+1} X_{j,t+1}^+ I^-_{i+1,z-1} \big)
\\
&= \sum_{z=1}^{s-1} q^{-(s+1)+ 2z} ( q^{\wt{a}_{ij}} - q^{\wt{a}_{i+1,j}+2}) I^+_{i,s-z} X_{j,t+1}^+ I^-_{i+1,z}  
	\\ & \quad 
	+ q^{\wt{a}_{ij} + s-1} I^+_{i,0} X_{j,t+1}^+ I^-_{i+1,s} 
	- q^{\wt{a}_{i+1,j} -s +1} I^+_{i,s} X_{j,t+1}^+ I^-_{i+1,0}. 
\end{align*}
Thus, we have 
\begin{align*}
&[J_{i,s+1}, X_{j,t}^+ ]
\\
&= \big\{ q^{\wt{a}_{ij} -s -1} I_{i,s}^+ - q^{- \wt{a}_{i+1,j} + s+1} I^-_{i+1,s} 
		- (q-q^{-1})  \sum_{z=1}^{s-1} q^{- \wt{a}_{i+1,j}-s + 2 z +1}  I^+_{i,s-z} I^-_{i+1,z}  \big\} X_{j,t+1}^+ 
	\\ & \quad 
	- X_{j,t+1}^+ \big\{ q^{-\wt{a}_{ij} - s-1} I_{i,s}^+ - q^{\wt{a}_{i+1,j} +s+1} I^-_{i+1,s} 
		- (q-q^{-1}) \sum_{z=1}^{s-1} q^{-\wt{a}_{ij} -s + 2 z -1} I^+_{i,s-z} I^-_{i+1,z} \big\} 
	\\ & \quad 
	- (q-q^{-1}) \sum_{z=1}^{s-1} q^{-(s+1)+ 2z} ( q^{\wt{a}_{ij}} - q^{\wt{a}_{i+1,j}+2}) I^+_{i,s-z} X_{j,t+1}^+ I^-_{i+1,z}  
	\\ & \quad 
	- q^{s-1} (q-q^{-1})   \big( q^{\wt{a}_{ij} } I^+_{i,0} X_{j,t+1}^+ - q^{-\wt{a}_{ij} } X_{j,t+1}^+ I^+_{i,0} \big)  I^-_{i+1,s} 
	\\ & \quad 
	+ q^{-s+1} (q-q^{-1})  I^+_{i,s} \big( q^{\wt{a}_{i+1,j}} X_{j,t+1}^+ I^-_{i+1,0} - q^{- \wt{a}_{i+1,j}}  I^-_{i+1,0} X^+_{j,t+1} \big). 
\end{align*}
Applying $J_{i,s} = q^{-s} I^+_{i,s} - q^s I^-_{i+1,s} - (q-q^{-1}) \sum_{z=1}^{s-1} q^{-s + 2 z} I^+_{i,s-z} I^-_{i+1,z}$ 
and the relation (Q'4-2),  
we have 
\begin{align}
\label{Jis+1 Xjt+} 
\begin{split}
[J_{i,s+1}, X_{j,t}^+ ]
&= q^{-\wt{a}_{i+1,j} +1} J_{i,s} X_{j,t+1}^+ 
	+ q^{-s} ( q^{\wt{a}_{ij} -1} - q^{-\wt{a}_{i+1,j} +1 } ) I^+_{i,s} X_{j,t+1}^+  
	\\ & \quad 
	- q^{-\wt{a}_{ij} -1} X_{j,t+1}^+ J_{i,s} + q^s (q^{\wt{a}_{i+1,j}+1} - q^{- \wt{a}_{ij} -1} ) X_{j,t+1}^+ I^-_{i+1,s}
	\\ & \quad 
	- (q-q^{-1}) \sum_{z=1}^{s-1} q^{- s + 2z} ( q^{\wt{a}_{ij} -1} - q^{\wt{a}_{i+1,j}+1}) I^+_{i,s-z} X_{j,t+1}^+ I^-_{i+1,z}  
	\\ & \quad 
	-  \wt{a}_{ij} q^{s-1} (q-q^{-1})   X_{j,t+1}^+  I^-_{i+1,s} 
	- \wt{a}_{i+1,j} q^{-s+1} (q-q^{-1})  I^+_{i,s}  X_{j,t+1}^+. 
\end{split}
\end{align}
On the other hand, 
by (Q'4-2) and (Q'4-3), we have 
\begin{align}
\label{I=iu X+jt+1}
\begin{split}
& I^+_{i,u} X_{j,t+1}^+ = X_{j,t+1}^+ I^+_{i,u}  \text{ if } i \not= j, j+1, 
\\
&X_{j,t+1}^+ I_{i+1,u}^- = I^-_{i+1,u} X_{j,t+1}^+  \text{ if } i \not= j -1, j
\end{split}
\end{align}
for $u \geq 0$. 
Then, the equations \eqref{Jis+1 Xjt+} and \eqref{I=iu X+jt+1} imply the relation (\roiv). 
The relation (\rov) follows from (\roiv) by applying the algbera anti-involution $\dag$. 
\end{proof}
\end{lem}

\begin{prop}
\label{Prop sl gl}
Put $\bQ_{[q]} = (q^{-1} Q_1, q^{-2} Q_2, \dots, q^{-(n-1)} Q_{n-1})$.  
Then, there exists the algebra homomorphism 
$\Up^{\lan \bQ \ran} : U_q(\Fsl_n^{\lan \bQ_{[q]} \ran}[x]) \ra U_q(\Fgl_n^{\lan \bQ \ran}[x])$ 
such that 
\begin{align*}
\Up^{\lan \bQ \ran} (X_{i,t}^{\pm}) = q^{ t i } X_{i,t}^{\pm},
\quad 
\Up^{\lan \bQ \ran} (J_{i,t}) = q^{ t i }J_{i,t}, 
\quad 
\Up^{\lan \bQ \ran} (K_i^{\pm}) = K_i^{\pm}.
\end{align*}
\begin{proof}
We can check the well-defindness of the homomorphism $\Up^{\lan \bQ \ran}$ by direct calculations 
using the defining relations of $U_q(\Fgl_n^{\lan \bQ \ran}[x])$ and Lemma \ref{Lemma some rel UqglQ}. 
\end{proof}
\end{prop}


\para 
For $q \in \CC^{\times}$ and $\wh{\bQ}=(\wh{Q}_0, \wh{Q}_1, \dots, \wh{Q}_{r-1}) \in \CC^{r}$, 
let $\sH_{m,r}$ be the Ariki-Koike algebra associated to the complex reflection group $\fS_m \ltimes (\ZZ/ r \ZZ)^m$ 
of type $G(r,1,m)$ with parameters $q$ and $\wh{\bQ}$. 
Namely, 
$\sH_{m,r}$ is the associative algebra over $\CC$ generated by $T_0$, $T_1$, \dots, $T_{m-1}$ 
subject to the defining relations 
\begin{align*}
&(T_0 - \wh{Q}_0) (T_0 - \wh{Q}_1) \dots (T_0- \wh{Q}_{r-1}) =0, 
\quad 
(T_i -q)(T_i+q^{-1})=0 \quad (1\leq i \leq m-1), 
\\
&T_0 T_1 T_0 T_1 = T_1 T_0 T_1 T_0, 
\quad 
T_i T_{i+1} T_i = T_{i+1} T_i T_{i+1} \quad (1\leq i \leq m-2), 
\\
&T_i T_j = T_j T_i \quad (|i-j| >1). 
\end{align*}

For $\bn=(n_1, n_2,\dots, n_r) \in \ZZ_{>0}^r$, 
let $\sS_{m,r}(\bn)$ be the cyclotomic $q$-Schur algebra associated to the Ariki-Koike algebra $\sH_{m,r}$ 
with respect to $\bn$ defined in \cite{DJM} (see also \cite[\S6]{W16} for definitions). 

An $r$-tuple of partitions $\la=(\la^{(1)}, \dots, \la^{(r)})$ is called an $r$-partition. 
For an $r$-partition $\la =(\la^{(1)}, \dots, \la^{(r)})$, 
we denote $\sum_{k=1}^r |\la^{(k)} |$ by $|\la|$, and we call it the size of $\la$. 
Set $\vL^+_{m,r} =\{ \la=(\la^{(1)}, \dots, \la^{(r)}) \text{ : $r$-partition} \mid |\la|=m, \, \ell(\la^{(r)}) \leq n_r\}$. 
For $\la \in \vL^+_{m,r}$, 
let $\D(\la)$ be the Weyl (cell) module corresponding to $\la$ constructed in \cite{DJM}. 
It is known that 
$\sS_{m,r}(\bn)$ is a quasi-hereditary cellular algebra with the set of standard modules $\{\D(\la) \mid \la \in \vL_{m,r}^+\}$
if $n_k \geq m$ for all $k$ by \cite{DJM}. 

\para 
For $\bn=(n_1,\dots, n_r) \in \ZZ^r_{>0}$, 
set $n = n_1+\dots+n_r$, 
$\vG(\bn) = \{(i,k) \mid 1 \leq i \leq n_k, \, 1 \leq k \leq r\}$ and $\vG'(\bn)= \vG(\bn) \setminus \{(n_r,r)\}$. 
We identify $\vG(\bn)$ with $\{1,2, \dots, n\}$ by the bijection 
\begin{align*}
\xi : \vG(\bn) \ra \{1,2,\dots,n\}  \text{ such that } \xi (i,k) = \sum_{j=1}^{k-1} n_j +i. 
\end{align*} 
Namely, $\vG(\bn)$ gives the separation of the set $\{1,2,\dots, n\}$ to $r$-parts with respect to $\bn$. 
Under this identification, we regard $(n_k+1, k)$ (resp. $(0,k)$) as $(1,k+1)$ (resp. $(n_{k-1},k-1)$). 
For $(i,k),(j,l) \in \vG(\bn)$,  
set  $\wt{a}_{(i,k)(j,l)} = \wt{a}_{\xi(i,k), \xi(j,l)}$. 
By \cite{W16}, 
the cyclotomic $q$-Schur algebra $\sS_{m,r}(\bn)$ is realized as a quotient of the algebra $\cU_{q, \wh{\bQ}}(\bn)$ 
defined as follows. 


\begin{definition}[{\cite[Definition 4.2]{W16}}] 
We define the associative algebra $\cU_{q, \wh{\bQ}}(\bn)$ over $\CC$ by the following generators and defining relations: 

\begin{description}
\item[Generators] 
$\cX_{(i,k),t}^{\pm}$, $\cI_{(j,l),t}^{\pm}$, $\wt{\cK}^{\pm}_{(j,l)}$ 
	($(i,k) \in \vG'(\bn)$, $(j,l) \in \vG(\bn)$, $t \in \ZZ_{\geq 0}$), 
	
\item[Defining relations] 
\begin{align*}
\tag{R1} 
&\wt{\cK}^+_{(j,l)} \wt{\cK}^-_{(j,l)} = 1 = \wt{\cK}^-_{(j,l)} \wt{\cK}^+_{(j,l)}, 
	\quad (\wt{\cK}^{\pm}_{(j,l)})^2 
	= 1 \pm (q-q^{-1}) \cI^{\mp}_{(j,l),0}, 
\\
\tag{R2} 
&[\wt{\cK}^+_{(i,k)}, \wt{\cK}^+_{(j,l)}] 
	= [\wt{\cK}^+_{(i,k)}, \cI^{\s}_{(j,l),t}] 
	= [\cI^\s_{(i,k),s}, \cI^{\s'}_{(j,l),t}]=0 
	\quad (\s, \s' \in \{+, - \}), 
\\
\tag{R3} 
&\wt{\cK}^+_{(i,k)} \cX^{\pm}_{(j,l),t} \wt{\cK}^-_{(i,k)} = q^{\pm \wt{a}_{(i,k) (j,l)}} \cX^{\pm}_{(j,l),t}, 
\\
\tag{R4}
\begin{split}
& q^{\pm \wt{a}_{(i,k)(j,l)}} \cI^{\pm}_{(i,k),0} \cX^+_{(j,l),t} - q^{\mp \wt{a}_{(i,k)(j,l)}} \cX_{(j,l),t}^+ \cI^{\pm}_{(i,k),0} 
	= \wt{a}_{(i,k)(j,l)} \cX_{(j,l),t}^+, 
\\
& q^{\mp \wt{a}_{(i,k)(j,l)}} \cI^{\pm}_{(i,k),0} \cX^-_{(j,l),t} - q^{\pm \wt{a}_{(i,k) (j,l)}} \cX_{(j,l),t}^- \cI^{\pm}_{(i,k),0} 
	= - \wt{a}_{(i,k)(j,l)} \cX^-_{(j,l),t}, 
\end{split}
\\
\tag{R5}
\begin{split}
& [\cI^{\pm}_{(i,k),s+1}, \cX^+_{(j,l),t}] 
	= q^{\pm \wt{a}_{(i,k)(j,l)}} \cI^{\pm}_{(i,k),s} \cX_{(j,l),t+1}^+ 
		- q^{\mp \wt{a}_{(i,k)(j,l)}} \cX_{(j,l),t+1}^+ \cI^{\pm}_{(i,k),s}, 
\\
& [\cI^{\pm}_{(i,k),s+1}, \cX^-_{(j,l),t}] 
	= q^{\mp \wt{a}_{(i,k)(j,l)}} \cI^{\pm}_{(i,k),s} \cX^-_{(j,l),t+1} - q^{\pm \wt{a}_{(i,k)(j,l)}} \cX^-_{(j,l),t+1} \cI^{\pm}_{(i,k),s}, 
\end{split} 
\end{align*}
\begin{align*}
\tag{R6} 
\\[-4em]
&[\cX^+_{(i,k),t}, \cX^-_{(j,l),s}] 
= \d_{(i,k)(j,l)} \begin{cases}
		\cK^+_{(i,k)} \cJ_{(i,k),s+t} & \text{ if } i \not=n_k, 
		\\
		- \wh{Q}_k \cK^+_{(n_k,k)} \cJ_{(n_k,k),s+t} + \cK^+_{(n_k,k)} \cJ_{(n_k,k), s+t+1} 
		& \text{ if } i =n_k, 
	\end{cases}
\\
\tag{R7} 
\begin{split}
& [\cX_{(i,k),t}^{\pm}, \cX_{(j,l),s}^{\pm}]=0 
	\quad \text{ if } (j,l) \not=(i,k), (i \pm1,k), 
\\
& \cX_{(i,k),t+1}^{\pm} \cX_{(i,k),s}^{\pm} - q^{\pm 2} \cX_{(i,k),s}^{\pm} \cX_{(i,k),t+1}^{\pm} 
	= q^{\pm 2} \cX_{(i,k),t}^{\pm} \cX_{(i,k),s+1}^{\pm} - \cX_{(i,k),s+1}^{\pm} \cX_{(i,k),t}^{\pm},  
\\
& \cX_{(i,k),t+1}^+ \cX_{(i+1,k),s}^+ - q^{-1} \cX_{(i+1,k),s}^+ \cX_{(i,k),t+1}^+ 
	= \cX_{(i,k),t}^+ \cX_{(i+1,k),s+1}^+ - q \cX_{(i+1,k),s+1}^+ \cX_{(i,k),t}^+, 
\\
& \cX_{(i+1,k),s}^- \cX_{(i,k),t+1}^- - q^{-1} \cX_{(i,k),t+1}^- \cX_{(i+1,k),s}^- 
	= \cX_{(i+1,k),s+1}^- \cX_{(i,k),t}^- - q \cX_{(i,k),t}^- \cX_{(i+1,k),s+1}^-, 
\end{split}
\\
\tag{R8}
\begin{split}
&\cX^+_{(i \pm 1,k),u} (\cX^+_{(i,k),s} \cX^+_{(i,k),t} + \cX^+_{(i,k),t} \cX^+_{(i,k),s} ) 
	+ (\cX^+_{(i,k),s} \cX^+_{(i,k),t} + \cX^+_{(i,k),t} \cX^+_{(i,k),s} ) \cX^+_{(i \pm 1,k),u} 
	\\
	&= (q+q^{-1}) (\cX^+_{(i,k),s} \cX^+_{(i \pm 1,k),u} \cX^+_{(i,k),t}  + \cX^+_{(i,k),t} \cX^+_{(i\pm1, k),u} \cX^+_{(i,k),s}), 
\\
&\cX^-_{(i \pm 1,k),u} (\cX^-_{(i,k),s} \cX^-_{(i,k),t} + \cX^-_{(i,k),t} \cX^-_{(i,k),s} ) 
	+ (\cX^-_{(i,k),s} \cX^-_{(i,k),t} + \cX^-_{(i,k),t} \cX^-_{(i,k),s} ) \cX^-_{(i \pm 1,k),u} 
	\\
	&= (q+q^{-1}) (\cX^-_{(i,k),s} \cX^-_{(i \pm 1,k),u} \cX^-_{(i,k),t}  + \cX^-_{(i,k),t} \cX^-_{(i\pm1, k),u} \cX^-_{(i,k),s}), 
\end{split}
\end{align*}
where we put $\cK^+_{(i,k)} = \wt{\cK}^+_{(i,k)} \wt{\cK}^-_{(i+1,k)}$, $\cK^-_{(i,k)} = \wt{\cK}^-_{(i,k)} \wt{\cK}^+_{(i+1,k)}$ 
and 
\begin{align*}
\cJ_{(i,k),t} 
= \begin{cases} 
	\cI^+_{(i,k),0} - \cI^-_{(i+1,k),0} + (q-q^{-1}) \cI^+_{(i,k),0} \cI^-_{(i+1,k),0} & \text{ if } t=0, 
	\\ \dis 
	q^{-t} \cI^+_{(i,k),t} - q^t \cI^-_{(i+1,k),t} - (q-q^{-1}) \sum_{z=1}^{t-1} q^{-t + 2 z} \cI^+_{(i,k),t-z} \cI^-_{(i+1,k),z} 
	& \text{ if } t>0.
	\end{cases} 
\end{align*}
\end{description}
\end{definition}

We remark that the parameter $\wh{Q}_0$ does not appear, 
and we do not need it, in the definition of the algebra $\cU_{q, \wh{\bQ}}(\bn)$. 
The parameter $\wh{Q}_0$ appears in the algebra homomorphism from $\cU_{q, \wh{\bQ}}(\bn)$ to the cyclotomic $q$-Schur algebra 
$\sS_{m,r}(\bn)$ given in \cite[Theorem 8.1]{W16}. 

We can easily prove the following proposition by checking defining relations.

\begin{prop}
\label{Prop iso gl g}
Assume that $\wh{Q}_i \not=0$ for all $1 \leq i\leq  r-1$. 
Set $\bQ'=(Q'_1, Q'_2, \dots, Q'_{n-1}) \in \CC^{n-1}$ as  
\begin{align}
\label{set Qi}
Q'_i= \begin{cases} 
	\wh{Q}_k^{-1}  & \text{ if } \xi^{-1}(i)=(n_k,k)  \text{ for some } k, 
	\\
	 0 & \text{ otherwise.} 
	 \end{cases} 
\end{align}
Then, there exists the algebra isomorphism 
$\Om^{\lan \wh{\bQ} \ran}_{\bn} : U_q(\Fgl_n^{\lan \bQ' \ran}[x]) \ra \cU_{q, \wh{\bQ}}(\bn)$
such that 
\begin{align*}
&\Om^{\lan \wh{\bQ} \ran}_{\bn} (X_{i,t}^+) = 
	\begin{cases}
		\cX^+_{\xi^{-1}(i),t} & \text{ if } \xi^{-1}(i) \not= (n_k,k) \text{ for all } k, 
		\\
		- \wh{Q}_k^{-1} \cX^+_{\xi^{-1}(i),t} & \text{ if } \xi^{-1}(i) =(n_k,k) \text{ for some } k, 
	\end{cases}
\\
& \Om^{\lan \wh{\bQ} \ran}_{\bn} (X_{i,t}^-) = \cX^-_{\xi^{-1}(i), t}, 
\quad 
\Om^{\lan \wh{\bQ} \ran}_{\bn} ( I^{\pm}_{j,t}) = \cI^{\pm}_{\xi^{-1}(j), t}, 
\quad \Om^{\lan \wh{\bQ} \ran}_{\bn} (\wt{K}^{\pm}_j) = \wt{\cK}^{\pm}_{\xi^{-1}(j)}.
\end{align*}
\end{prop}

\para 
Let 
$\Psi^{\lan \wh{Q} \ran}_{\bn} : \cU_{q, \wh{\bQ}} (\bn) \ra \sS_{m,r}(\bn)$ 
be the algebra homomorphism given in \cite[Theorem 8.1]{W16}. 
Assume that $\wh{Q}_i \not=0$ for all $1 \leq i\leq  r-1$. 
Set $\bQ'=(Q'_1,\dots, Q'_{n-1}) \in \CC^{n-1}$ as \eqref{set Qi}, 
and put $\bQ =(Q_1, Q_2,\dots, Q_{n-1}) = \bQ'_{[q]}$. 
Namely, we have 
\begin{align}
\label{Put Qi} 
Q_i = \begin{cases}
		q^{-(n_1+ \dots + n_k)} \wh{Q}_k^{-1} & \text{ if } \xi^{-1}(i) =(n_k,k) \text{ for some } k, 
		\\
		0 & \text{ otherwise.}
	\end{cases}
\end{align}
Then, we have the algebra homomorphism 
\begin{align}
\Phi^{\lan \wh{\bQ} \ran}_{\bn} :=  \Psi^{\lan \wh{\bQ} \ran}_{\bn} \circ \Om^{\lan \wh{\bQ} \ran}_{\bn} \circ \Up^{\lan \bQ' \ran} 
	: U_q (\Fsl_n^{\lan \bQ \ran}[x]) \ra \sS_{m,r}(\bn). 
\end{align}
Through the algebra homomorphism $\Phi^{\lan \wh{\bQ} \ran}_{\bn}$, 
we regard $\sS_{m,r}(\bn)$-modules as $U_q(\Fsl_n^{\lan \bQ \ran}[x])$-modules. 

\begin{prop}
\label{Prop hw Weyl}
Assume that $q$ is not a root of unity, and $n_k \geq m$ for all $k$. 
For $\la \in \vL^+_{m,r}$, 
the Weyl module $\D(\la)$ is a highest weight $U_q(\Fsl_n^{\lan \bQ \ran}[x])$-module, 
and the highest weight of $\D(\la)$ is given by 
$(\bu^{\lan Q_i \ran} (\vf_i) )_{i \in I}$, where 
\begin{align*}
\vf_i = \begin{cases}
		\dis 
		\prod_{p=1}^{\la_j^{(k)} - \la_{j+1}^{(k)}} \big( x - q^{i- 2 j + 2 \la_j^{(k)} - 2 (p-1)} \wh{Q}_{k-1}  \big) 
		\\ \dis 
		\hspace{10em} 
		\text{ if } i = \sum_{l=1}^{k-1} n_l+j  \text{ for some $k$ and } 1\leq j < n_k,
	\\
	\dis 
	q^{-\la_1^{(k+1)}} \prod_{p=1}^{\la_{n_k}^{(k)}} \big( x - q^{i- 2 n_k + 2 \la_{n_k}^{(k)} - 2 (p-1)} \wh{Q}_{k-1}  \big) 
		\quad 
		\text{ if } i = \sum_{l=1}^k n_l \text{ for some } k.
	\end{cases}
\end{align*}

\begin{proof}
By the definition of $\Phi^{\lan \wh{\bQ} \ran}_{\bn}$ together with the argument in \cite{W16}, 
we see that the Weyl module $\D(\la)$ ($\la \in \vL^+_{m,r}$) is a highest weight $U_q(\Fsl_n^{\lan \bQ \ran}[x])$-module. 
Let $v_0 \in \D(\la)$ be a highest weight vector. 
For $i \in I$, put $(j,k) = \xi^{-1}(i)$. 
Then, 
by \cite[Theorem 8.3]{W16} together with the definition of $\Phi^{\lan \wh{\bQ} \ran}_{\bn}$, 
we have 
\begin{align*}
K_i \cdot v_0 
&= \wt{\cK}_{(j,k)}^+ \wt{\cK}_{(j+1,k)}^- \cdot v_0 
=\begin{cases}
	q^{\la^{(k)}_j - \la^{(k)}_{j+1}} v_0 & \text{ if } j \not= n_k, 
	\\
	q^{\la_{n_k}^{(k)} - \la_1^{(k+1)}} v_0  & \text{ if } j=n_k, 
\end{cases}  
\\
J_{i,t} \cdot v_0 
&=q^{t i } \big( q^{-t} \cI^+_{(j,k),t} - q^t \cI^-_{(j+1,k),t} 
			- (q-q^{-1}) \sum_{z=1}^{t-1} q^{-t+ 2 z} \cI^+_{(j,k),t-z} \cI^-_{(j+1,k),z}\big) \cdot v_0 
\\
&= \begin{cases} 
	q^{ (i - 2 j) t}   (\wh{Q}_{k-1})^t \big\{ 
		q^{( 2 t -1) \la_j^{(k)}  } [\la_j^{(k)}]  
		-  q^{\la_{j+1}^{(k)}  } [\la_{j+1}^{(k)}] 
		\\ \hspace{5em}  \dis 
		- (q-q^{-1}) \sum_{z=1}^{t-1} q^{ ( 2  (t-  z) -1) \la_j^{(k)} + \la_{j+1}^{(k)} } 
		[\la_j^{(k)}]   [\la_{j+1}^{(k)}] \big\} 
		\quad \text{ if } j \not=n_k, 
	\\
	q^{(i - 2 j) t } \big\{ (\wh{Q}_{k-1})^t q^{( 2 t -1) \la_j^{(k)} } [\la_j^{(k)}]  
	-  (\wh{Q}_{k})^t q^{2 j t  + \la_{1}^{(k+1)} } [\la_{1}^{(k+1)}] 
	\\ \qquad  \dis 
	- (q-q^{-1}) \sum_{z=1}^{t-1} 
		(\wh{Q}_{k-1})^{t-z} (\wh{Q}_{k})^z q^{ 2 j z + ( 2 (t-z) -1) \la_j^{(k)} + \la_{1}^{(k+1)} } [\la_j^{(k)}] [\la_{1}^{(k+1)}] \big\} 
		\\ \hspace{27em} \text{ if } j =n_k 
	\end{cases}
\\
&= \begin{cases}
		q^{ (i - 2 j + 2 \la_j^{(k)}) t}   (\wh{Q}_{k-1})^t q^{ -( \la_j^{(k)} - \la_{j+1}^{(k)})} [\la_j^{(k)} - \la_{j+1}^{(k)}]  
		\quad \text{ if } j \not=n_k, 
	\\ 
	q^{(i-2j + 2 \la_j^{(k)})t} (\wh{Q}_{k-1})^t  q^{- \la_j^{(k)}} [\la_j^{(k)}]  
	+ \big( - q^{ \la_1^{(k+1)}  } [\la_1^{(k+1)}] \big) (q^{-i} \wh{Q}_k^{-1})^{-t}
	\\ \qquad \dis 
	+ (q-q^{-1}) \sum_{z=1}^{t-1}  \big( - q^{\la_1^{(k+1)}  } [\la_1^{(k+1)}] \big) (q^{-i} \wh{Q}_k^{-1})^{-z} 
		\\ \hspace{10em} \times 
		\big( q^{ (i - 2 j + 2  \la_j^{(k)} ) (t-z)   } (\wh{Q}_{k-1})^{t-z}  q^{-\la_j^{(k)}} [\la_j^{(k)}] \big) 
		\quad \text{ if } j=n_k. 
	\end{cases}
\end{align*}
On the other hand, for $b \in \CC$ and $c,z  \in \ZZ_{>0}$, 
we can show that 
\begin{align*}
p_z (q) (b, b  q^{-2 \cdot 1}, b q^{- 2 \cdot 2}, \dots, b q^{- 2 \cdot (c-1)}) 
= b^z q^{-c} [c]
\end{align*}
by the induction on $c$ using Lemma \ref{Lemma ptq} (\roi). 
Thus, we have 
\begin{align*}
J_{i,t} \cdot v_0 
&= \begin{cases} 
	p_t \big( b, b q^{-2 \cdot1}, \dots, b q^{- 2 \cdot (\la_j^{(k)} - \la_{j+1}^{(k)}-1)} \big) 
	& \text{ if } j \not=n_k, 
	\\ 
	p_t \big( b, b q^{-2 \cdot 1}, \dots, b q^{-2 \cdot (\la_j^{(k)}-1)} \big)   + \wt{\b} Q_i^{-t} 
		\\ \qquad \dis 
	+ (q-q^{-1}) \sum_{z=1}^{t-1} \wt{\b} Q_i^{-z} p_{t-z} \big(b, b q^{- 2 \cdot 1}, \dots, b q^{- 2 \cdot (\la_j^{(k)} -1)} \big) 
	& \text{ if } j =n_k, 
	\end{cases}
\end{align*}
where we put $b = q^{(i-2j + 2 \la_j^{(k)})} (\wh{Q}_{k-1})$ and $\wt{\b} = - q^{ \la_1^{(k+1)}  } [\la_1^{(k+1)}]$, 
and we note that $Q_i = q^{-i} \wh{Q}_k^{-1}$ if $ j=n_k$ by \eqref{Put Qi}. 
Moreover, 
applying the definition \eqref{def ptQqb} to the above equation in the case where 
$j =n_k$, 
we have 
\begin{align*}
J_{i,t} \cdot v_0 
= p_t^{\lan Q_i \ran} (q; \b) (b, bq^{-2 \cdot 1}, \dots, b q^{-2 \cdot (\la_j^{(k)}-1)} ) 
\text{ if } j =n_k, 
\end{align*}
where 
$\b = q^{-\la_1^{(k+1)}} $ 
since 
$ \wt{\b} = - q^{ \la_1^{(k+1)}  } [\la_1^{(k+1)}] = (1 - q^{2 \la_1^{(k+1)}}) (q-q^{-1})^{-1}$. 
Note that 
$i = \sum_{l=1}^{k-1} n_l +j$ since $(j,k) = \xi^{-1}(i)$, 
then we obtain the proposition. 
\end{proof}
\end{prop}





\end{document}